\documentclass{my_gOMS2e}

\usepackage{subfigure}

\usepackage[table]{xcolor}

\theoremstyle{plain}

\theoremstyle{remark}

\newtheorem{definition}{Definition}[subsection]
\theoremstyle{definition}

\graphicspath{
{figures/betas/}
{figures/pops/}
{figures/eps1/}
{figures/eps2/}
{figures/tradeoffs/}
{figures/methods/}
}

\begin{document}

\received{Received 09 May 2014; Revised 19 October 2014; Accepted 30-Nov-2014}

\title{Multiobjective approach to optimal control for a tuberculosis model}

\author{Roman Denysiuk$^{\rm a}$,
Cristiana J. Silva$^{\rm b}$,
Delfim F. M. Torres$^{\rm b}$$^{\ast}$
\thanks{$^{\ast}$ Corresponding author. Email: delfim@ua.pt \vspace{6pt}}\\\vspace{6pt}
$^{a}${\em{Algoritmi R\&D Center, University of Minho, Campus de Gualtar, 4710-057 Braga, Portugal}};
$^{b}${\em{Center for Research and Development in Mathematics and Applications (CIDMA),
Department of Mathematics, University of Aveiro, 3810-193 Aveiro, Portugal}}}

\maketitle


\begin{abstract}
Mathematical modelling can help to explain the nature and dynamics of infection transmissions,
as well as support a policy for implementing those strategies that are most likely to bring public
health and economic benefits. The paper addresses the application of optimal control strategies
in a tuberculosis model. The model consists of a system of ordinary differential equations,
which considers reinfection and post-exposure interventions. We propose a multiobjective optimization
approach to find optimal control strategies for the minimization of active infectious and persistent
latent individuals, as well as the cost associated to the implementation of the control strategies.
Optimal control strategies are investigated for different values of the model parameters.
The obtained numerical results cover a whole range of the optimal control strategies,
providing valuable information about the tuberculosis dynamics
and showing the usefulness of the proposed approach.

\bigskip

\end{abstract}



\section{Introduction}

Tuberculosis (TB) is an important international public health issue.
\emph{Mycobacterium tuberculosis} is a pathogenic bacterial species
that is the cause of most occurrences of tuberculosis. It is spread
through the air when people who have an active TB infection cough,
sneeze, or transmit respiratory fluids through the air. It typically
affects the lungs (pulmonary TB) but can affect other sites as well
(extrapulmonary TB). The classic symptoms of active TB infection are
a chronic cough with blood-tinged sputum, fever, night sweats, and weight loss.
Infection of other organs causes a wide range of symptoms. Only approximately
10\% of people infected with \emph{mycobacterium tuberculosis} develop active TB disease,
whereas approximately 90\% of infected people remain latent. Latent infected TB people
are asymptomatic and do not transmit TB, but may progress to active TB through
either endogenous reactivation or exogenous reinfection \cite{Small_Fuj_2001,Styblo_1978}.
Treatment of TB  depends on close cooperation between patients and health care providers.
One can distinguish three types of TB treatment: (i) vaccination to prevent infection;
(ii) treatment to cure active TB; (iii) treatment of latent TB to prevent
endogenous reactivation~\cite{Abu_Raddad:etal:2009,Gomes_etall_2007}.
The treatment of active infectious individuals can have different timings~\cite{Kruk_etall_2008}.
Here we consider treatment with the duration of six months \cite{MyID:305}.
One of the difficulties related to the success of these treatments is to make sure that patients complete them.
After two months, patients no longer have symptoms of the disease and feel healed,
so many of them stop taking the medicines. When the treatment is not concluded,
the patients are not cured and reactivation can occur and/or the patients
may develop resistant TB. One way to prevent patients of not completing the treatment
is based on supervision and patient support. This is one of the measures proposed
by the Direct Observation Therapy (DOT) of World Health Organization
(WHO)~\cite{WHO_treatTB_2010}. One example of treatment supervision consists
in recording each dose of anti-TB drugs on the patients treatment card~\cite{WHO_treatTB_2010}.
These measures are very expensive since the patients need to stay longer in a hospital
or specialized people are to be payed to supervise patients till they finish their treatment.
On the other hand, it is recognized that the treatment of latent TB individuals reduces
the chances of reactivation \cite{Gomes_etall_2007}.

Although the anti-TB drugs developed since 1940 have dramatically
reduced mortality rates (in clinical cases, cure rates of 90\% have
been documented)~\cite{WHO_2012}, TB remains a major health problem.
In 2012, there were 8.6 million of new TB cases and 1.3 million of TB deaths.
TB is the second leading cause of death from an infectious disease worldwide
after HIV~\cite{WHO_2012}. Mortality rates are especially high without treatment.

Optimal control theory is a branch of mathematics developed to find optimal ways to control
a dynamic system~\cite{Cesari_1983,Fleming_Rishel_1975,Pontryagin_et_all_1962}. Nowadays,
the usefulness of optimal control theory in epidemiology is
well recognized~\cite{Kong,livro_Lenhart_2007,Rodrigues_Monteiro_Torres_2010,Rodrigues_Monteiro_Torres_Zinober_2011}.
Although different optimal control problems have been recently proposed and applied to
TB~\cite{Bowong_2010,Emvudu_et_all,Haffat_et_all}, results in tuberculosis are scarce~\cite{SLenhart_2002}
especially those where optimal strategies are found with respect to different conflicting objectives.
Our goal is to show how multiobjective optimization can be used
for finding the optimal control strategies in a tuberculosis model.
Despite the clear multiobjective nature of the underlying problem,
it has been only solved in the past using single-objective approaches.
Our main contributions are: (i) to show how to generate the whole range
of the optimal strategies, (ii) to provide the analysis of the obtained results
with clear advantages with respect to available results in the literature,
and (iii) to promote multiobjective optimization in epidemiology.

The paper~\cite{SLenhart_2002} studies a mathematical model for TB based on~\cite{Castillo_Chavez_1997},
considering two classes of infected and latent individuals (infected with typical TB and with resistant strain TB).
The authors seek to reduce the number of infected and latent individuals with resistant TB. In~\cite{Emvudu_et_all},
the model considers the existence of a class called \emph{the lost to follow up individuals} and they propose optimal
control strategies for the reduction of the number of individuals in this class. In~\cite{Haffat_et_all},
the authors adapt a  model from~\cite{Castillo_Chavez_2000} where exogenous reinfection is considered
and wish to minimize the number of infectious individuals. In~\cite{Bowong_2010}, a TB model that incorporates
exogenous reinfection, chemoprophylaxis of latently infected individuals and treatment of infections is proposed.
Optimal control strategies based on chemoprophylaxis of latently infected individuals and treatment of infectious
individuals are analyzed for the reduction of the number of active infected individuals. A TB model, which considers
reinfection and post-exposure interventions, was proposed in~\cite{Gomes_etall_2007}. The importance of considering
reinfection and post-exposure interventions was previously justified in \cite{Bandera_2001,Caminero_2001,van_Rie_1999}.
In~\cite{SilvaTorres2012,SilvaTorres2013}, a TB model from~\cite{Gomes_etall_2007} is extended by adding two controls
and two real parameters associated with controls. Optimal strategies are found by minimizing a cost functional that includes
the number of active TB infectious and persistent latent individuals as well as the cost of the measures for treatments.
An overview of existing TB modelling studies is presented in~\cite{Houben_eta_ll_2014} (see also \cite{Cohen:etal:2006}).
These works also identify high-priority areas and challenges for future modelling efforts, which are:
(i) the difficult diagnosis and high mortality of TB-HIV; (ii) the high risk of disease progression;
(iii) TB health systems in high HIV prevalence settings; (iv) uncertainty in the natural progression of TB-HIV;
and (v) combined interventions for TB-HIV.

In this paper we use the model suggested in~\cite{SilvaTorres2013} and propose a multiobjective approach
to find optimal control strategies. This approach reflects the intrinsic nature of an underlying
decision-making problem. It avoids the use of additional parameters needed to formulate the cost
functional and allows to obtain the whole range of optimal solutions. We also investigate the effect
of different model parameters. Section~\ref{sec:model} presents the mathematical model for tuberculosis
with controls. The optimal control problem is then formulated in Section~\ref{sec:problem}, while
in Section~\ref{sec:results} we present and discuss optimal control strategies obtained
by numerical simulations, considering several variations of some of the model parameters.
Section~\ref{sec:5} compares our approach to find the optimal control
strategies, based on the $\epsilon$-constraint method,
with the goal attainment and Chebyshev methods.
We end with Section~\ref{sec:conclusions} of conclusions
and some directions of future work.


\section{Tuberculosis model with controls}
\label{sec:model}

In the following, we consider a TB model taken from~\cite{SilvaTorres2013}. The model without controls
is based on~\cite{Gomes_etall_2007} and considers reinfection and post-exposure interventions,
consisting of a system of nonlinear ordinary differential equations representing population dynamics.
In the model, the population is divided into five categories:
\begin{quote}
\begin{tabular}{lcl}
$S(t)$ & -- & susceptible; \\
$L_1(t)$ & -- & early latent, \textrm{i.e.}, individuals recently infected (less than two years) but\\
&& not infectious; \\
$I(t)$ & -- & infected, \textrm{i.e.}, individuals who have active tuberculosis and are infectious; \\
$L_2(t)$  & -- & persistent latent, \textrm{i.e.}, individuals who were infected and remain latent; \\
$R(t)$ & -- & recovered, \textrm{i.e.}, individuals who were previously infected and treated. \\
\end{tabular}
\end{quote}
It is assumed that at birth all individuals are equally susceptible and differentiate as they experience
infection and respective therapy~\cite{Gomes_etall_2007}. The total population, $N$, is assumed to be constant,
so, $N = S(t) + L_1(t) + I(t) + L_2(t) + R(t)$. This way, it is assumed that the rates of birth and death,
$\mu$, are equal (corresponding to a mean life time of 70 years~\cite{Gomes_etall_2007})
and there are no disease-related deaths.

The model includes control variables representing prevention and treatment measures,
which are continuously implemented during a considered period of disease treatment:
\begin{quote}
\begin{tabular}{lcl}
$u_1(t)$ &--& represents the effort that prevents the failure of treatment in active \\
&& TB infectious individuals, $I$, e.g., supervising the patients, helping them \\
&& to take the TB medications regularly and to complete the TB treatment;\\
$u_2(t)$ &--& represents the fraction of persistent latent individuals, $L_2$, that is identified \\
&& and put under treatment.\\
\end{tabular}
\end{quote}

According to~\cite{SilvaTorres2013}, the tuberculosis is modeled
by the nonlinear time-varying state equations
\begin{equation}
\label{tb:model}
\small
\begin{cases}
\dot{S}(t) = \mu N - \frac{\beta}{N} I(t) S(t) - \mu S(t)\\
\dot{L_1}(t) = \frac{\beta}{N} I(t)\left( S(t) + \sigma L_2(t)
+ \sigma_R R(t)\right) - (\delta + \tau_1 + \mu)L_1(t)\\
\dot{I}(t) = \phi \delta L_1(t) + \omega L_2(t) + \omega_R R(t)
- (\tau_0 + \epsilon_1 u_1(t) + \mu) I(t)\\
\dot{L_2}(t) = (1 - \phi) \delta L_1(t) - \sigma \frac{\beta}{N} I(t) L_2(t)
- (\omega + \epsilon_2 u_2(t) + \tau_2 + \mu)L_2(t)\\
\dot{R}(t) = (\tau_0 + \epsilon_1 u_1(t) )I(t) +  \tau_1 L_1(t)
+ \left(\tau_2 +\epsilon_2 u_2(t)\right) L_2(t)
- \sigma_R \frac{\beta}{N} I(t) R(t) - \left(\omega_R + \mu\right) R(t)
\end{cases}
\end{equation}
with the initial conditions
\begin{equation}
\label{tb:initial}
S(0) = \frac{76}{120}N, \quad L_1(0) = \frac{37}{120}N,
\quad I(0) = \frac{4}{120}N, \quad L_2(0) = \frac{2}{120}N, \quad R(0) = \frac{1}{120}N.
\end{equation}

\begin{table}
\centering
\small
\begin{tabular}{|l|l|l|}
\hline
Symbol & Description  & Value \\
\hline
$\beta$ & Transmission coefficient & $75, 100, 150, 175$ \\
$\mu$ & Death and birth rate &  $1/70 \, yr^{-1}$\\
$\delta$ & Rate at which individuals leave $L_1$ & $12 \, yr^{-1}$\\
$\phi$ & Proportion of individuals going to $I$ & $0.05$\\
$\omega$ & Rate of endogenous reactivation for persistent latent infections & $0.0002 \, yr^{-1}$\\
$\omega_R$ & Rate of endogenous reactivation for treated individuals & $0.00002 \, yr^{-1}$\\
$\sigma$ & Factor reducing the risk of infection as a result of acquired  & \\
& immunity to a previous infection for $L_2$ & $0.25$ \\
$\sigma_R$ & Rate of exogenous reinfection of treated patients & 0.25 \\
$\tau_0$ & Rate of recovery under treatment of active TB &  $2 \, yr^{-1}$\\
$\tau_1$ & Rate of recovery under treatment of latent individuals $L_1$ &  $2 \, yr^{-1}$\\
$\tau_2$ & Rate of recovery under treatment of latent individuals $L_2$ &  $1 \, yr^{-1}$\\
$N$ & Total population & $30000, 40000, 60000$ \\
$\epsilon_1$ & Efficacy of treatment of active TB $I$ & $0.25, 0.5, 0.75$ \\
$\epsilon_2$ & Efficacy of treatment of latent TB $L_2$ & $0.25, 0.5, 0.75$ \\
$T$ & Total simulation duration & $5 \, yr$ \\
\hline
\end{tabular}
\caption{Model parameters.}\label{parameters}
\end{table}
\noindent The values of the model parameters presented in the control system~\eqref{tb:model}
are given in Table~\ref{parameters}. The values of the rates $\delta$, $\phi$, $\omega$,
$\omega_R$, $\sigma$ and $\tau_0$ are taken from~\cite{Gomes_etall_2007}
and the references cited therein. The parameter $\delta$ denotes the rate
at which individuals leave $L_1$ compartment; $\phi$ is the proportion
of individuals going to compartment $I$; $\omega$ is the rate of endogenous reactivation
for persistent latent infections (untreated latent infections); $\omega_R$ is the rate
of endogenous reactivation for treated individuals (for those who have undergone
a therapeutic intervention). The parameter $\sigma$ is the factor that reduces the risk of infection,
as a result of acquired immunity to a previous infection, for persistent latent individuals,
\textrm{i.e.}, this factor affects the rate of exogenous reinfection of untreated individuals;
while $\sigma_R$ represents the same parameter factor but for treated patients.

The parameter $\tau_0$ is the rate of recovery under treatment of active TB,
assuming an average duration of infectiousness of six months.
The parameters $\tau_1$ and $\tau_2$ apply to latent individuals $L_1$ and $L_2$, respectively.
They are the rates at which chemotherapy or a post-exposure vaccine is applied \cite{Abu_Raddad:etal:2009}.
We consider that the rate of recovery of early latent individuals under post-exposure interventions,
$\tau_1$, is equal to the rate of recovery under treatment of active TB, $\tau_0$, and greater
than the rate of recovery of persistent latent individuals under post-exposure interventions, $\tau_2$.

The parameters $\epsilon_i \in (0, 1)$, $i=1, 2$,
measure the effectiveness of the controls $u_i$, $i=1, 2$, respectively, \textrm{i.e.},
these parameters measure the efficacy of treatment interventions for active and persistent
latent TB individuals, respectively.

According to~\cite{Styblo_1991}, the risk of developing disease after infection is much higher
in the first five years following infection, and decline exponentially after that. Therefore,
as done in~\cite{SilvaTorres2013}, we consider here the total simulation duration of $T=5$.


\section{Problem formulation}
\label{sec:problem}

The main goal of our study is to find the most effective ways of applying the controls in~\eqref{tb:model},
aimed at restriction of the tuberculosis epidemic. In this section, we present two formulations
of an optimal control problem: (i) the first one, based on optimal control theory;
(ii) the second one, based on multiobjective optimization.


\subsection{Optimal control problem}

The aim is to find the optimal values $u_1^*$ and $u_2^*$ of the controls $u_1$ and $u_2$,
such that the associated state trajectories $S^*$, $L_1^*$, $I^*$, $L_2^*$, $R^*$
are solution of the system \eqref{tb:model} in the time interval $[0, T]$,
with initial conditions~\eqref{tb:initial}, and minimize an objective functional.
Here the objective functional considers the number of active TB infectious individuals $I$,
the number of persistent latent individuals $L_2$, and the implementation cost
of the strategies associated to the controls $u_i$, $i=1, 2$. The controls
are bounded between $0$ and $1$. When the controls vanish, no extra measures
are implemented for the reduction of $I$ and $L_2$; when the controls take the maximum value $1$,
the magnitude of the implemented measures, associated to $u_1$ and $u_2$,
take the value of the effectiveness of the controls, $\epsilon_1$ and $\epsilon_2$, respectively.

Consider the state system of ordinary differential equations~\eqref{tb:model}
and the set of admissible control functions given by
\[
\Omega = \{ (u_1(\cdot),u_2(\cdot)) \in (L^{\infty}(0,T))^2 \, | \,
0 \leq u_1(t),u_2(t) \leq 1, \forall t \in [0,T] \}.
\]
The objective functional is defined by
\begin{equation}
\label{oc:problem}
J(u_1(\cdot),u_2(\cdot))
= \int_{0}^{T}\left[ I(t) + L_2(t) + W_1u_1^2(t) + W_2u_2^2(t) \right]dt,
\end{equation}
where the constants $W_1$ and $W_2$ are a measure of the relative cost
of the interventions associated to the controls $u_1$ and $u_2$, respectively.
We consider the optimal control problem of determining
$\left(S^*(\cdot), L_1^*(\cdot), I^*(\cdot), L_2^*(\cdot), R^*(\cdot)\right)$,
associated to an admissible control pair
$\left(u_1^*(\cdot), u_2^*(\cdot) \right) \in \Omega$ on the time interval $[0, T]$,
satisfying \eqref{tb:model}, the initial conditions~\eqref{tb:initial},
and minimizing the cost function \eqref{oc:problem}, \textrm{i.e.},
\[
J(u_1^*(\cdot), u_2^*(\cdot)) = \min_{\Omega} J(u_1(\cdot), u_2(\cdot)) \, .
\]


\subsection{Multiobjective optimization}

Without loss of generality, a multiobjective optimization problem
with $m$ objectives and $n$ decision variables can be formulated as follows:
\begin{equation}\label{mop}
\begin{array}{rl}
\text{minimize:} & \boldsymbol{f}(\boldsymbol{x})
=(f_1(\boldsymbol{x}),f_2(\boldsymbol{x}),\ldots,f_m(\boldsymbol{x}))^{\text{T}}\\
\text{subject to:} & \boldsymbol{x} \in \Omega,
\end{array}
\end{equation}
where $\boldsymbol{x}$ is the decision vector, $\Omega \subseteq \mathbb{R}^n$
is the feasible decision space, and $\boldsymbol{f}(\boldsymbol{x})$ is the
objective vector defined in the objective space $\mathbb{R}^m$.

When several objectives are simultaneously optimized, there
is no natural ordering in the objective space. The objective space is partially ordered.
In such a scenario, solutions are compared on the basis of the Pareto dominance relation.

\begin{definition}[Pareto dominance]
For two solutions $\boldsymbol{a}$ and $\boldsymbol{b}$ from $\Omega$,
a solution $\boldsymbol{a}$ is said to dominate a solution $\boldsymbol{b}$
(denoted by $\boldsymbol{a} \prec \boldsymbol{b}$) if
\[
\forall i \in \{ 1, \ldots ,m\} : f_i (\boldsymbol{a}) \le f_i (\boldsymbol{b})
\, \wedge \, \exists j \in \{ 1, \ldots ,m\} :f_j (\boldsymbol{a}) < f_j (\boldsymbol{b}).
\]
\end{definition}

The presence of multiple conflicting objectives gives
rise to a set of optimal solutions, generally known as the Pareto optimal set.
The concepts of optimality for multiobjective optimization are defined as follows.

\begin{definition}[Pareto optimality]
A solution $\boldsymbol{x}^{\ast} \in \Omega$ is Pareto optimal if and only if
\[
\nexists \boldsymbol{y} \in \Omega :\boldsymbol{y} \prec \boldsymbol{x}^{\ast}.
\]
\end{definition}

\begin{definition}[Pareto optimal set]
For a multiobjective optimization problem~\eqref{mop},
the Pareto optimal set (or Pareto set, for short) is defined as
$$
\mathcal{PS}^{\ast} = \{ \boldsymbol{x}^{\ast} \in \Omega \, |
\, \nexists \boldsymbol{y} \in \Omega : \boldsymbol{y} \prec \boldsymbol{x}^{\ast} \}.
$$
\end{definition}

\begin{definition}[Pareto optimal front]
For a multiobjective optimization problem~\eqref{mop}
and the Pareto optimal set $\mathcal{PS}^{\ast}$, the Pareto
optimal front (or Pareto front, for short) is defined as
$$
\mathcal{PF}^{\ast} = \{ \boldsymbol{f}(\boldsymbol{x}^{\ast})
\in \mathbb{R}^m \, | \, \boldsymbol{x}^{\ast} \in \mathcal{PS}^{\ast} \}.
$$
\end{definition}


\subsection{Proposed approach}

An approach based on optimal control theory allows to obtain a single optimal solution
for the cost functional \eqref{oc:problem}, which is defined from some decision
maker's perspective using the constants $W_1$ and $W_2$. The most straightforward disadvantage
is that only a limited amount of information about the choice of the optimal strategy
can be presented to the decision maker. Moreover, the choice of proper values
of $W_1$ and $W_2$ is not straightforward, generally being not an easy task.

In our approach, we decompose the cost functional shown in~\eqref{oc:problem}
into two components representing different aspects taken into consideration
when dealing with tuberculosis. Then, we use a multiobjective optimization
method to simultaneously optimize the defined objectives. Thus,
the multiobjective optimization problem is defined as:
\begin{equation}
\label{mo:problem}
\begin{array}{rl}
\text{minimize} & f_1(u_1(\cdot),u_2(\cdot))
= \displaystyle \int_{0}^{T} \left[ I(t) + L_2(t) \right] dt, \\[0.2 cm]
& f_2(u_1(\cdot),u_2(\cdot)) = \displaystyle \int_{0}^{T} \left[ u_1^2(t) + u_2^2(t) \right] dt,
\end{array}
\end{equation}
where $f_1$ represents the number of active infected and latent individuals,
and $f_2$ represents the cost associated to the implementation of the control policies.


\section{Numerical experiments}
\label{sec:results}

We now present and discuss the numerical results for the optimal controls
using the multiobjective optimization approach. We also investigate the effects
of different parameters values on the optimal control strategies
and dynamics of the tuberculosis model.


\subsection{$\epsilon$-Constraint method}

The $\epsilon$-constraint method was introduced in \cite{Haimes1971}.
In this method, one of the objective functions is selected to be minimized,
whereas all the other functions are converted into constrains by setting an upper bound
to each of them. The problem to be solved is of the following form:
\begin{equation}
\label{method:eps}
\begin{array}{rlll}
\underset{\boldsymbol{x} \in \Omega}{\text{minimize:}} & f_l(\boldsymbol{x}) & \\
\text{subject to:} & f_i(\boldsymbol{x})\leq \epsilon_i,
& \forall i \in \{1,\ldots,m\} \wedge i \neq l.
\end{array}
\end{equation}
In the above formulation, the $l$th objective is minimized, and the parameter
$\epsilon_i$ represents an upper bound of the value of $f_i$. The $\epsilon$-constraint
method is able to obtain solutions in convex and nonconvex regions of the Pareto optimal front.
When all the objective functions in the problem~\eqref{mop} are convex,
as it happens in our study, the problem \eqref{method:eps} is also convex and it has a unique solution.
The unique solution of the problem~\eqref{method:eps} is Pareto optimal for any given
upper bound vector $\boldsymbol{\epsilon} = \{\epsilon_1, \ldots, \epsilon_{m-1}\}$.
For a proof see \cite{MiettinenBook}.


\subsection{Experimental setup}
\label{sec:setup}

The system~\eqref{tb:model} is numerically integrated using the fourth-order Runge--Kutta method.
The controls are discretized using $60$ equally spaced time intervals over the period of $[0,T]$.
The integrals in~\eqref{mo:problem} are calculated using the trapezoidal rule.
Using the formulation~\eqref{method:eps}, we minimize the first
objective in \eqref{mo:problem} setting $f_2$ as the constraint bounded
by the values of $\epsilon$, which are selected by dividing the range of
$f_2$ into 100 evenly distributed intervals. The range of $f_2$ is calculated as
$f_2^{\min} = 0$ for $u_1(\cdot), u_2(\cdot) \equiv 0$
and $f_2^{\max} = 10$ for $u_1(\cdot),u_2(\cdot) \equiv 1$.
To solve the problems with different $\epsilon$,
we use the MATLAB\textsuperscript{\textregistered} built-in function \texttt{fmincon}
with a sequential quadratic programming algorithm,
setting the maximum number of function evaluations to $20,000$.


\subsection{Experimental results}

\begin{figure}
\centering
\includegraphics[width=0.67\textwidth]{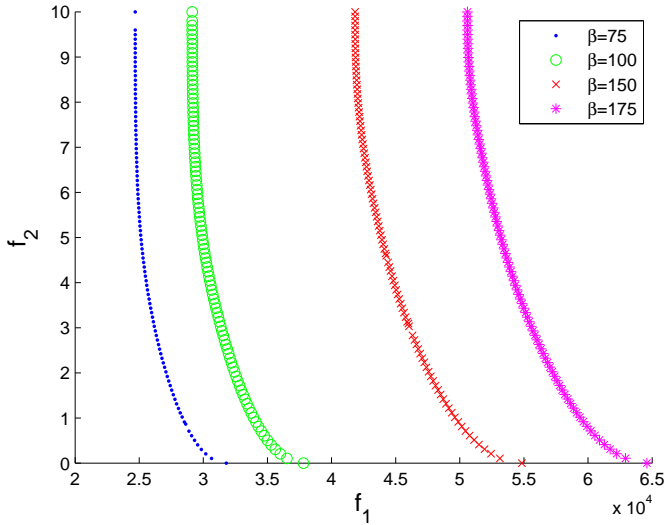}
\caption{Trade-off curves for different values of $\beta$.}
\label{tb:tradeoffs:betas}
\end{figure}
\begin{figure}
\centering
\subfigure[$I$ for $f_2=0$]{\epsfig{file=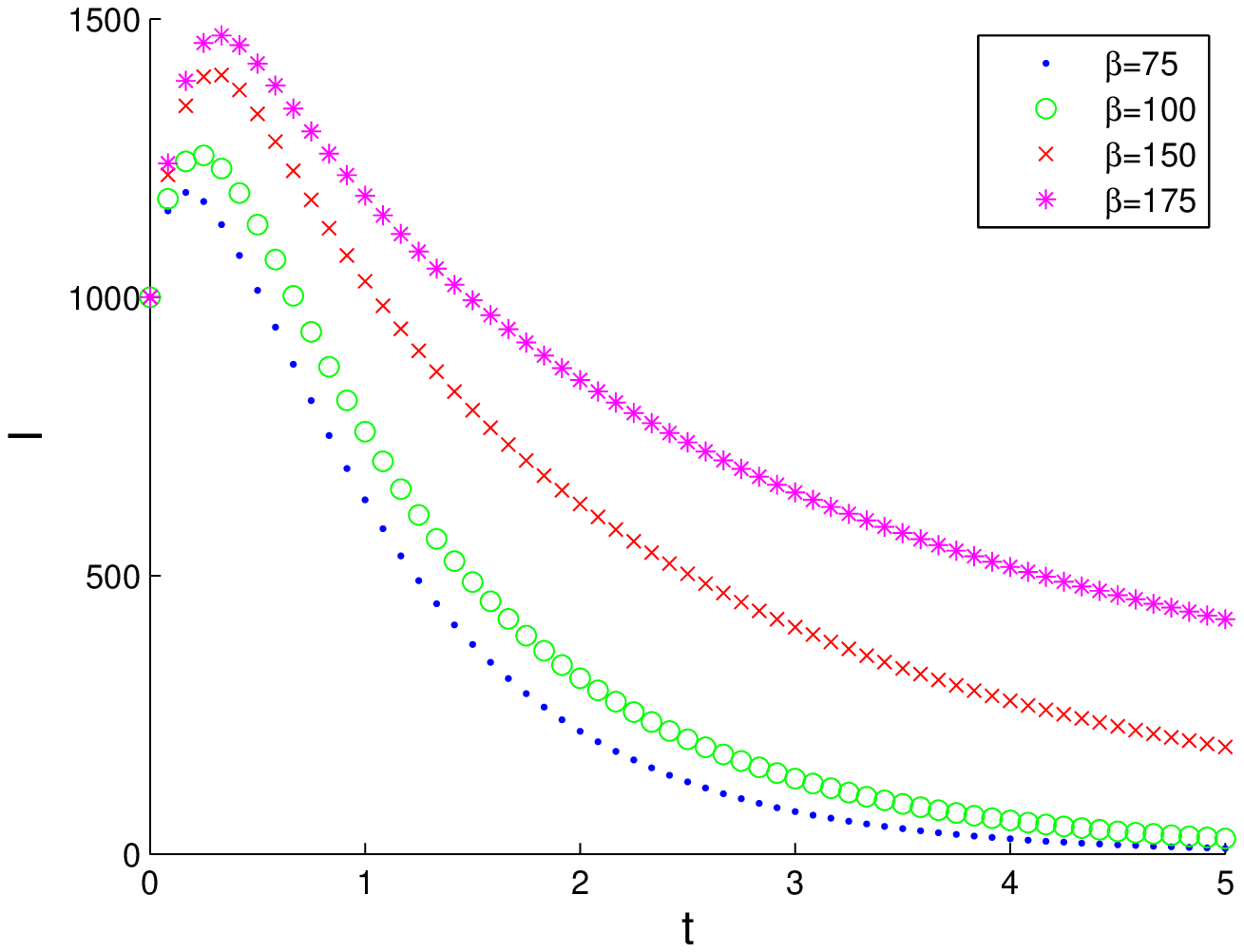,width=0.333\linewidth}\label{tb:sol:betas:a}}%
\subfigure[$L_2$ for $f_2=0$]{\epsfig{file=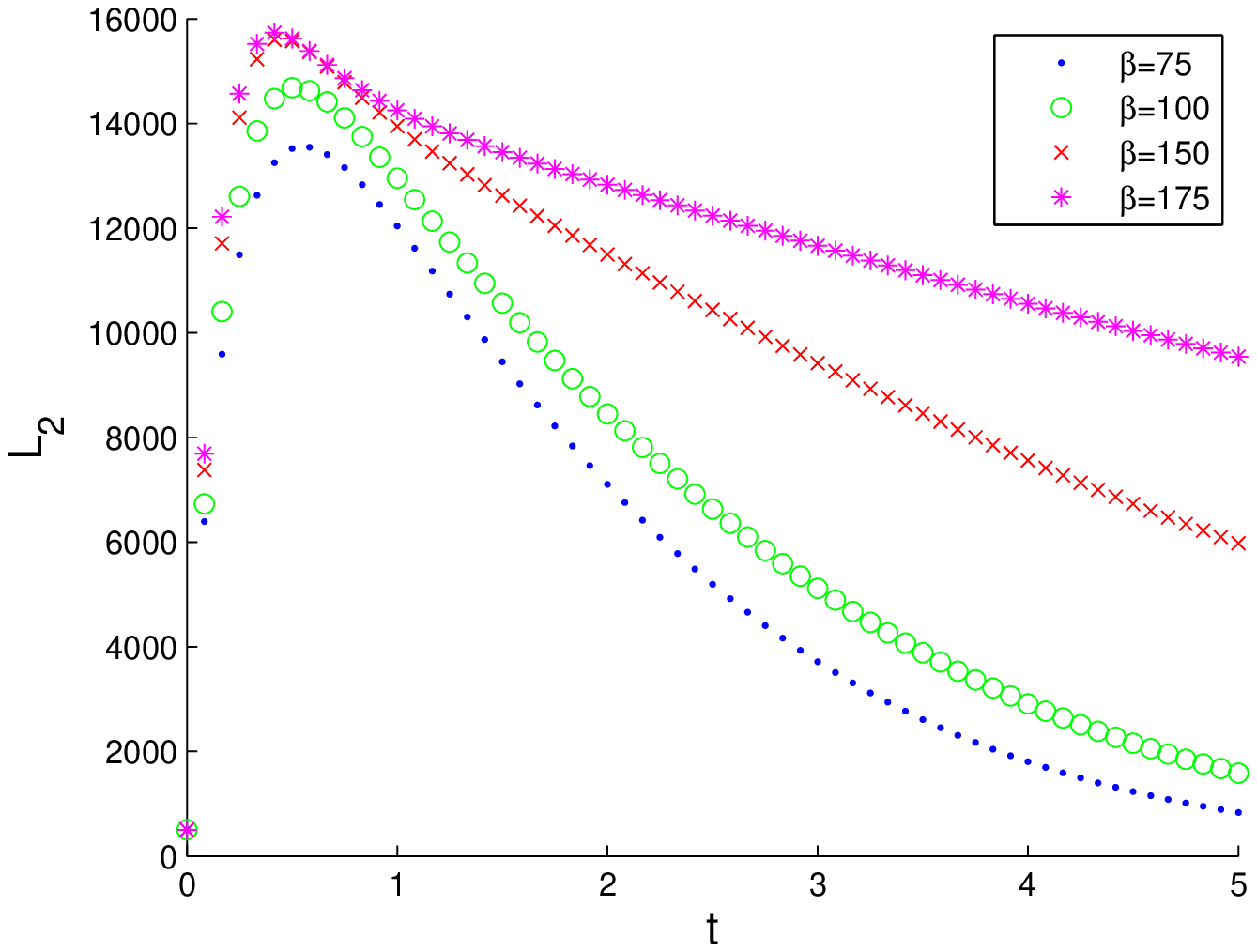,width=0.333\linewidth}\label{tb:sol:betas:b}}
\subfigure[$u_1$ for $f_2=2.5$]{\epsfig{file=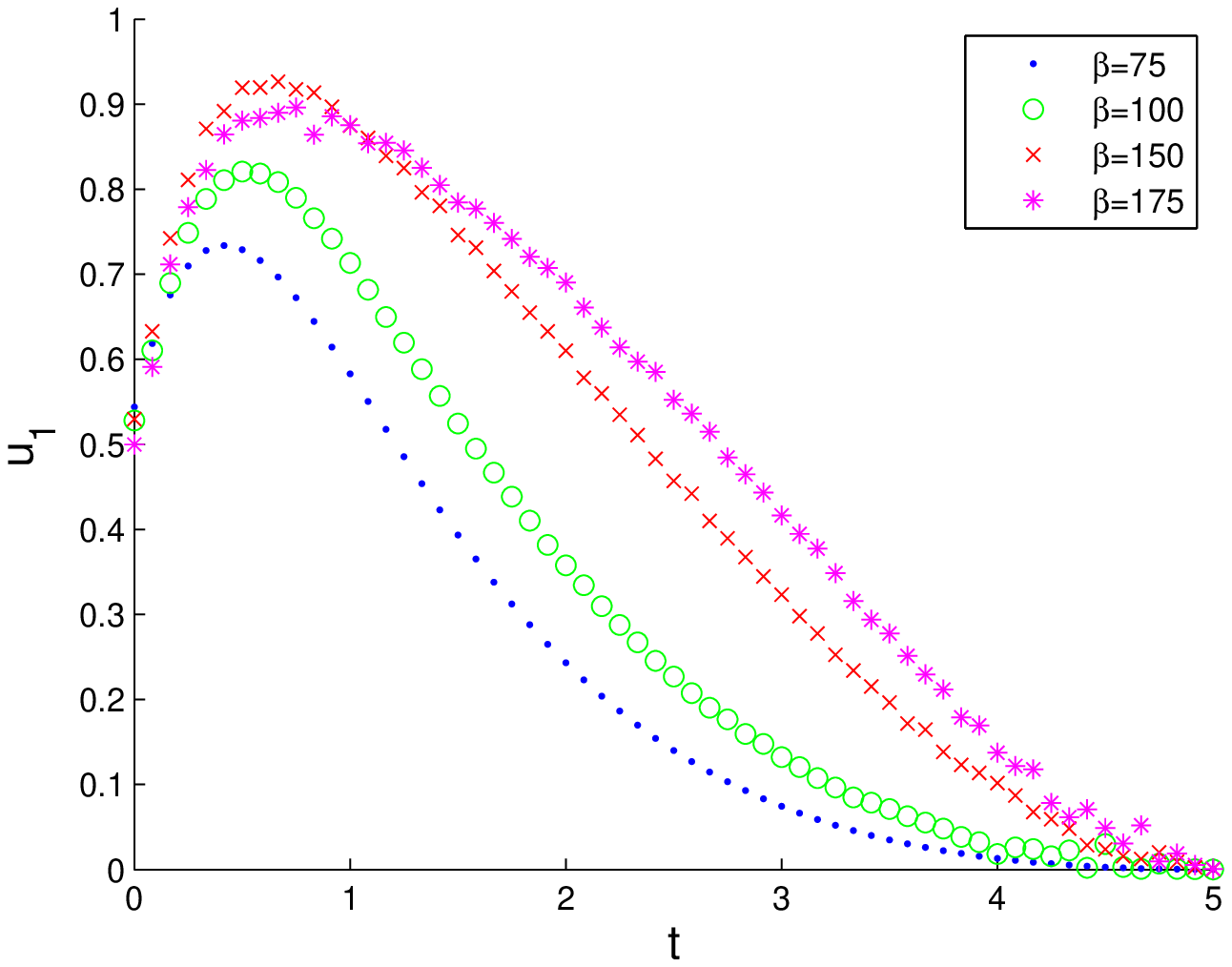,width=0.333\linewidth}\label{tb:sol:betas:c}}%
\subfigure[$u_2$ for $f_2=2.5$]{\epsfig{file=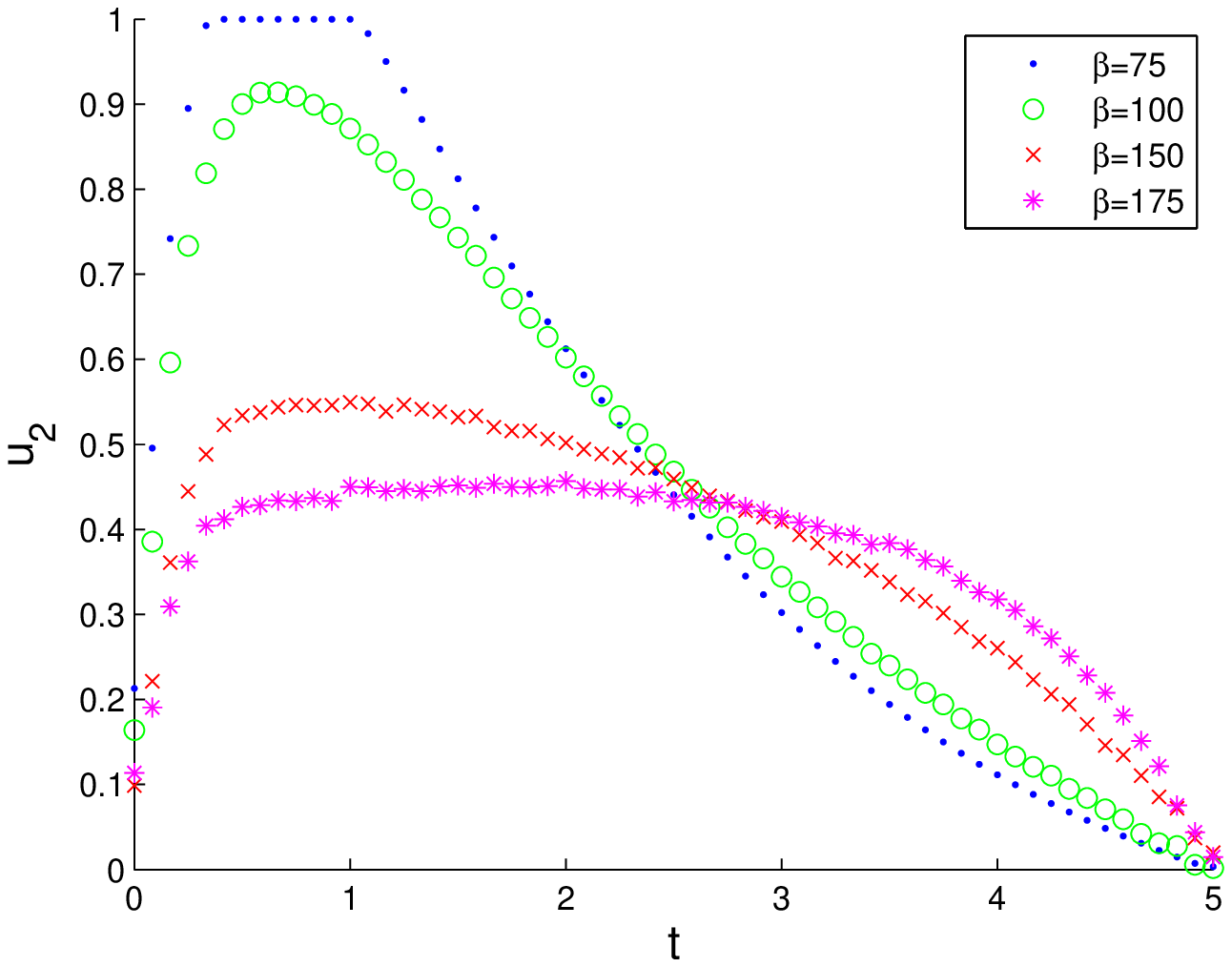,width=0.333\linewidth}\label{tb:sol:betas:d}}%
\subfigure[$I+L_2$ for $f_2=2.5$]{\epsfig{file=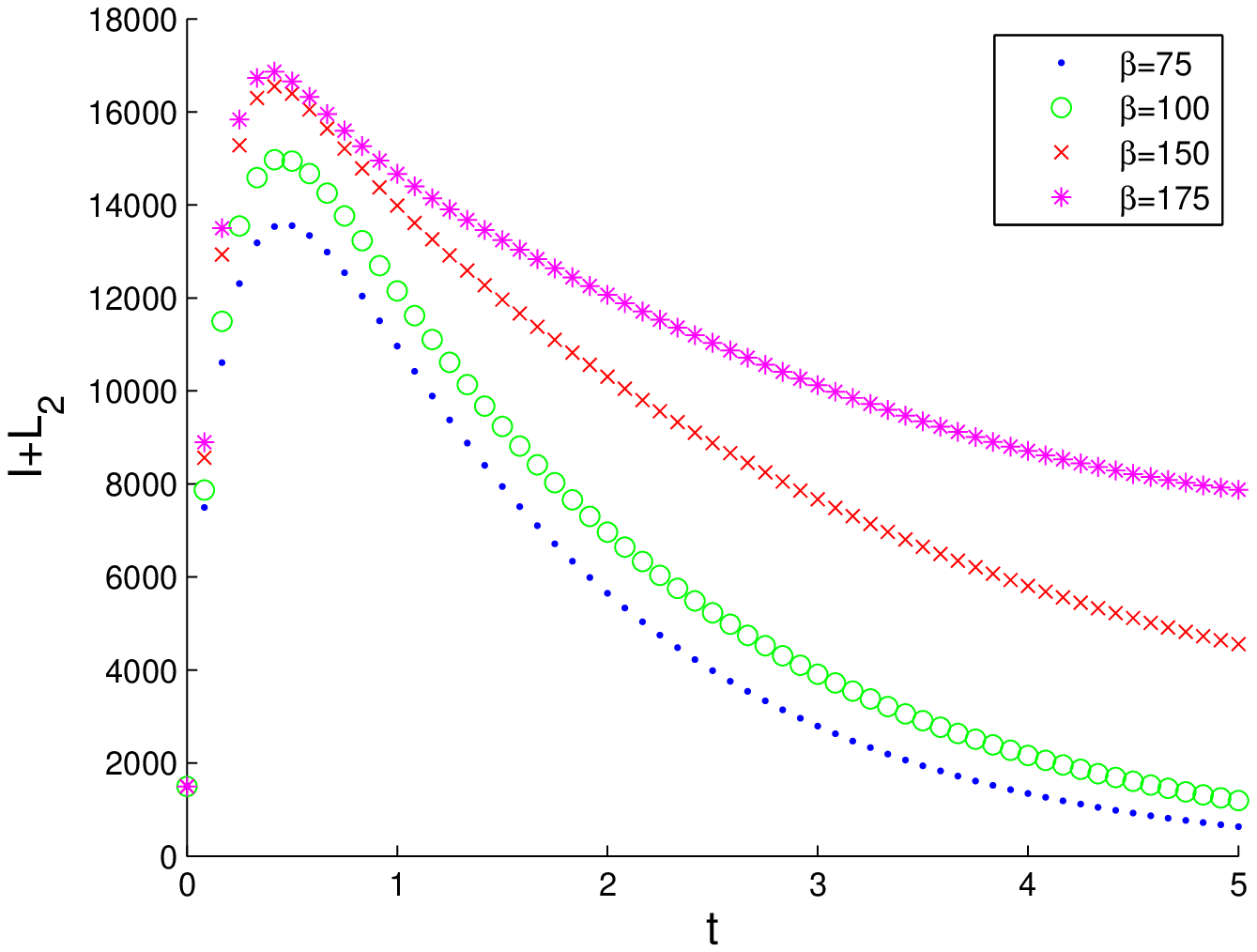,width=0.333\linewidth}\label{tb:sol:betas:e}}
\subfigure[$u_1$ for $f_2=5$]{\epsfig{file=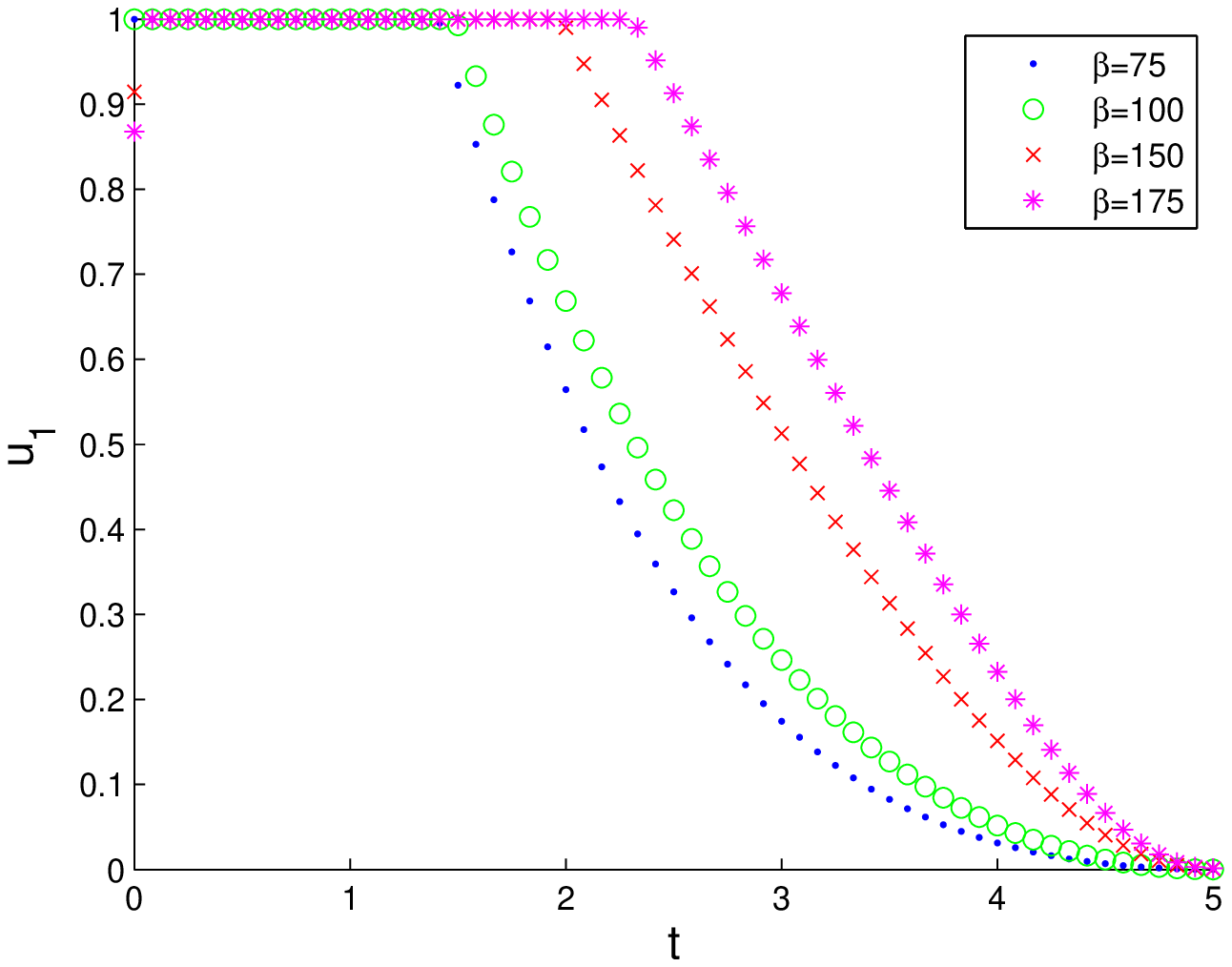,width=0.333\linewidth}\label{tb:sol:betas:f}}%
\subfigure[$u_2$ for $f_2=5$]{\epsfig{file=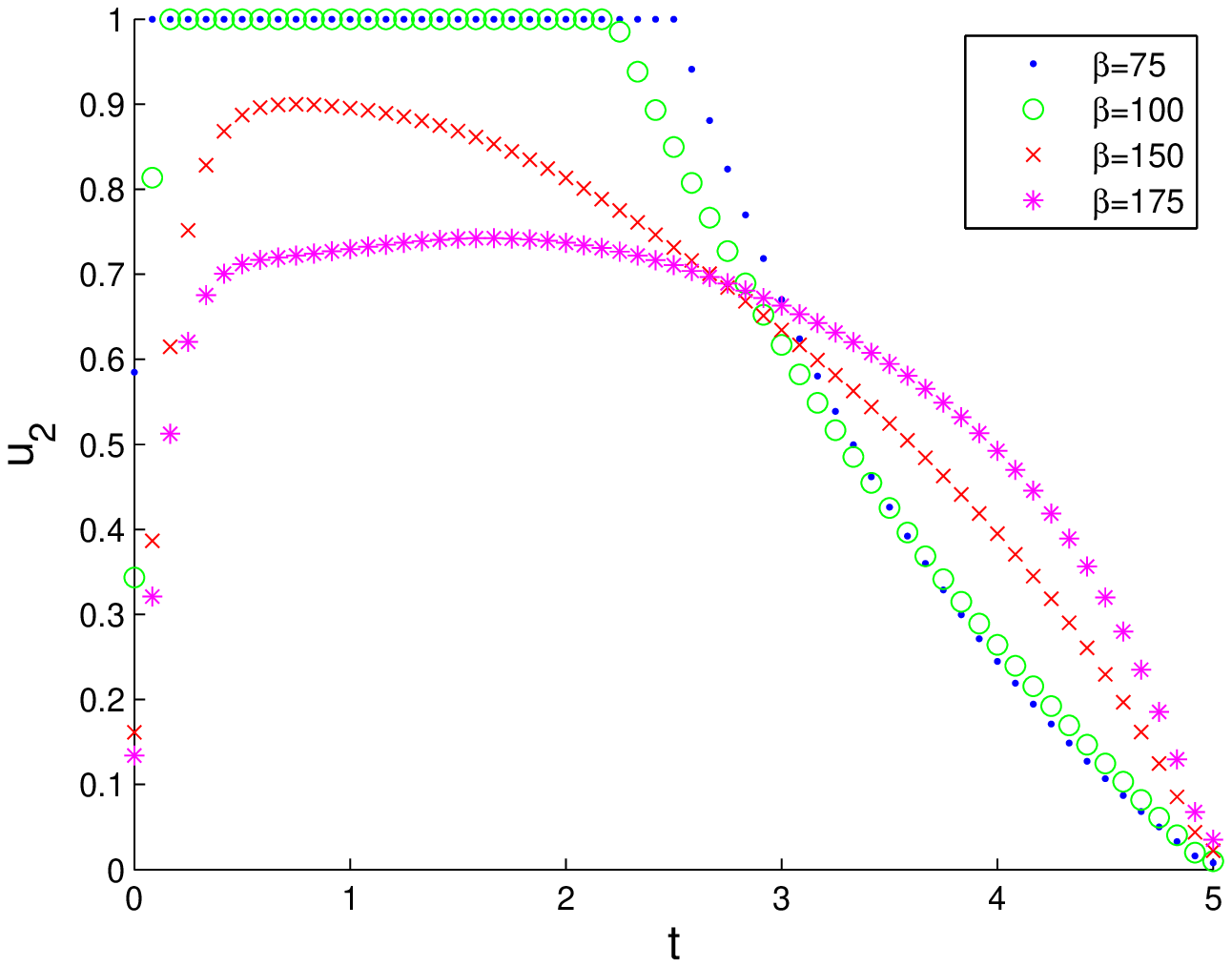,width=0.333\linewidth}\label{tb:sol:betas:g}}%
\subfigure[$I+L_2$ for $f_2=5$]{\epsfig{file=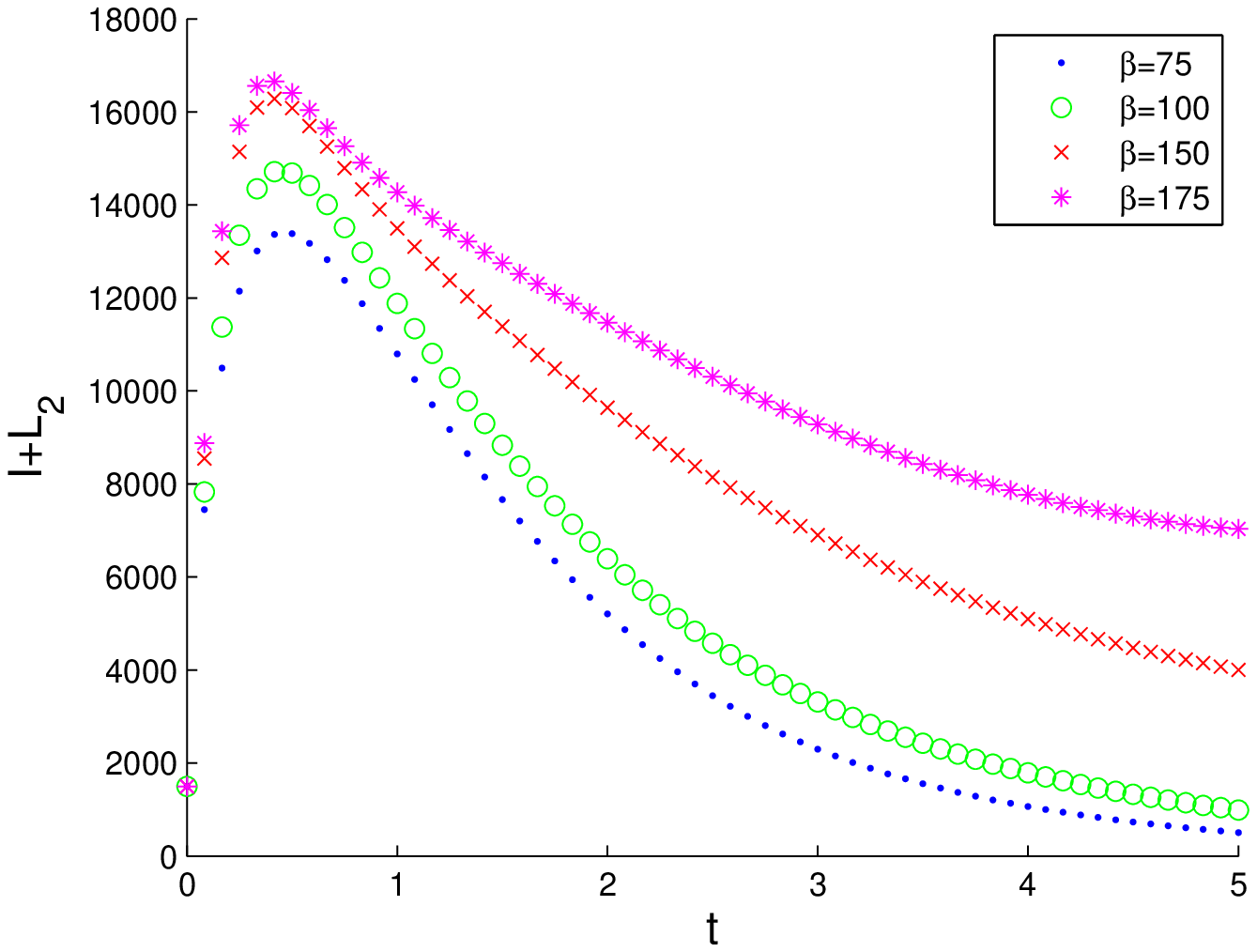,width=0.333\linewidth}\label{tb:sol:betas:h}}
\subfigure[$u_1$ for $f_2=7.5$]{\epsfig{file=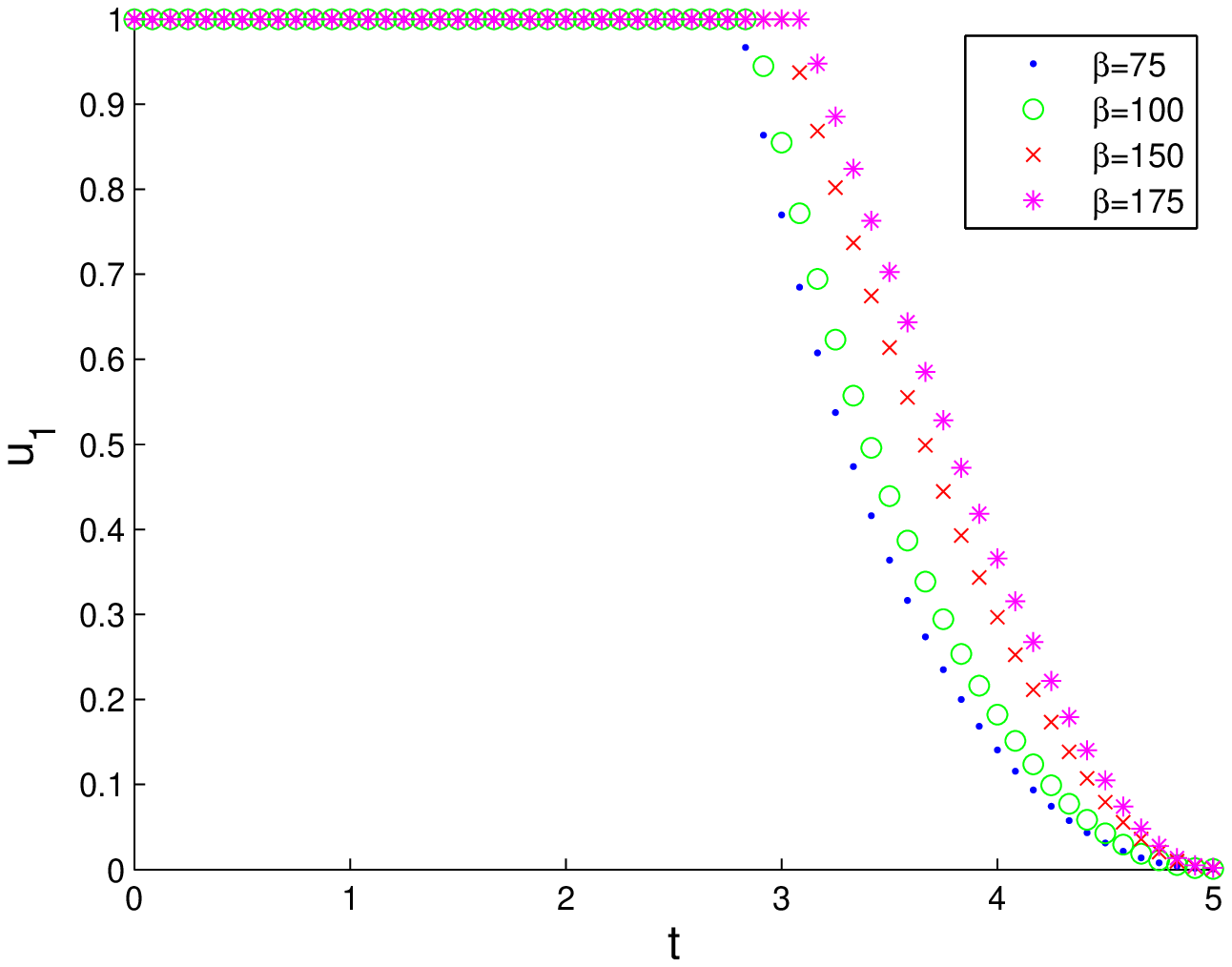,width=0.333\linewidth}\label{tb:sol:betas:i}}%
\subfigure[$u_2$ for $f_2=7.5$]{\epsfig{file=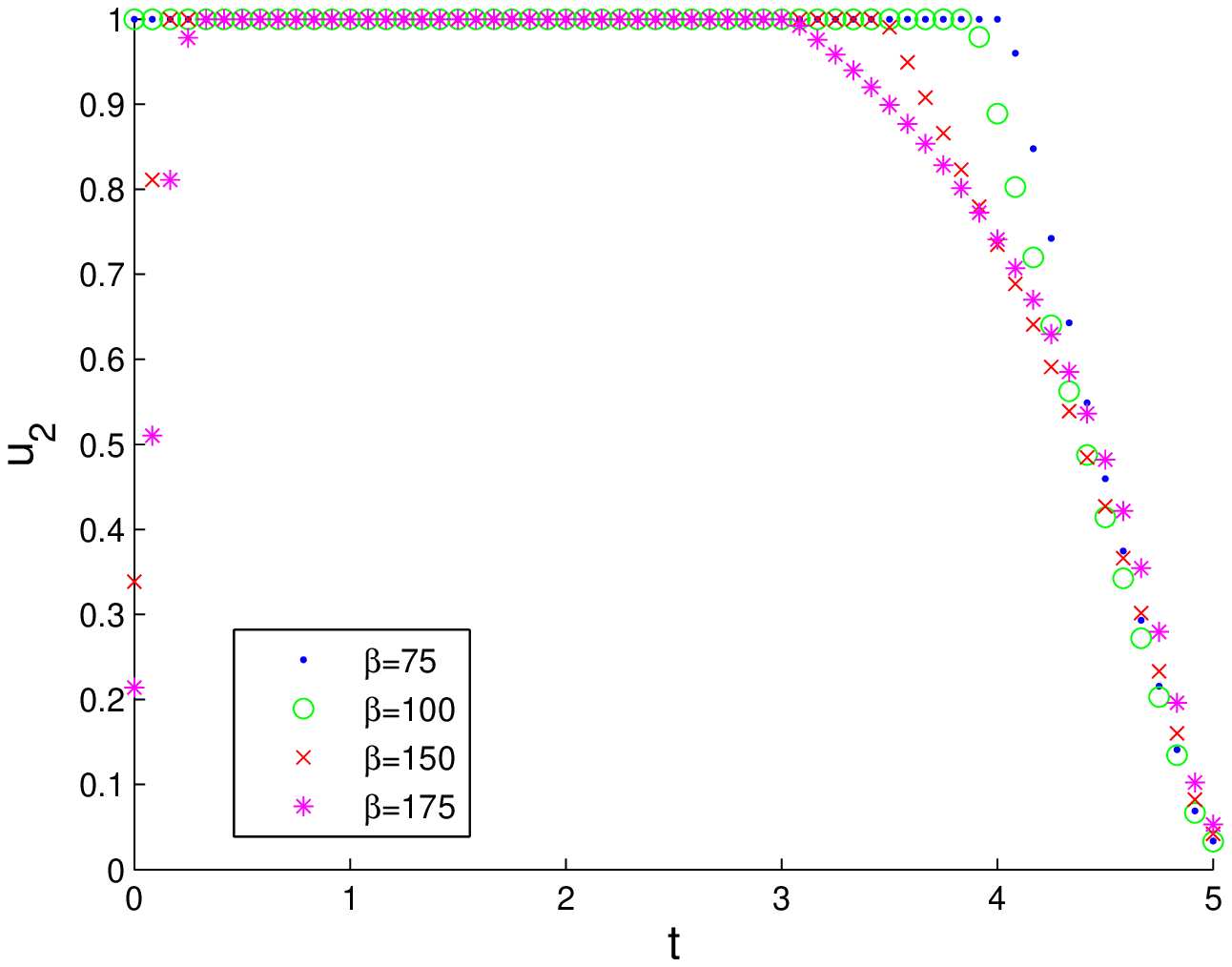,width=0.333\linewidth}\label{tb:sol:betas:j}}%
\subfigure[$I+L_2$ for $f_2=7.5$]{\epsfig{file=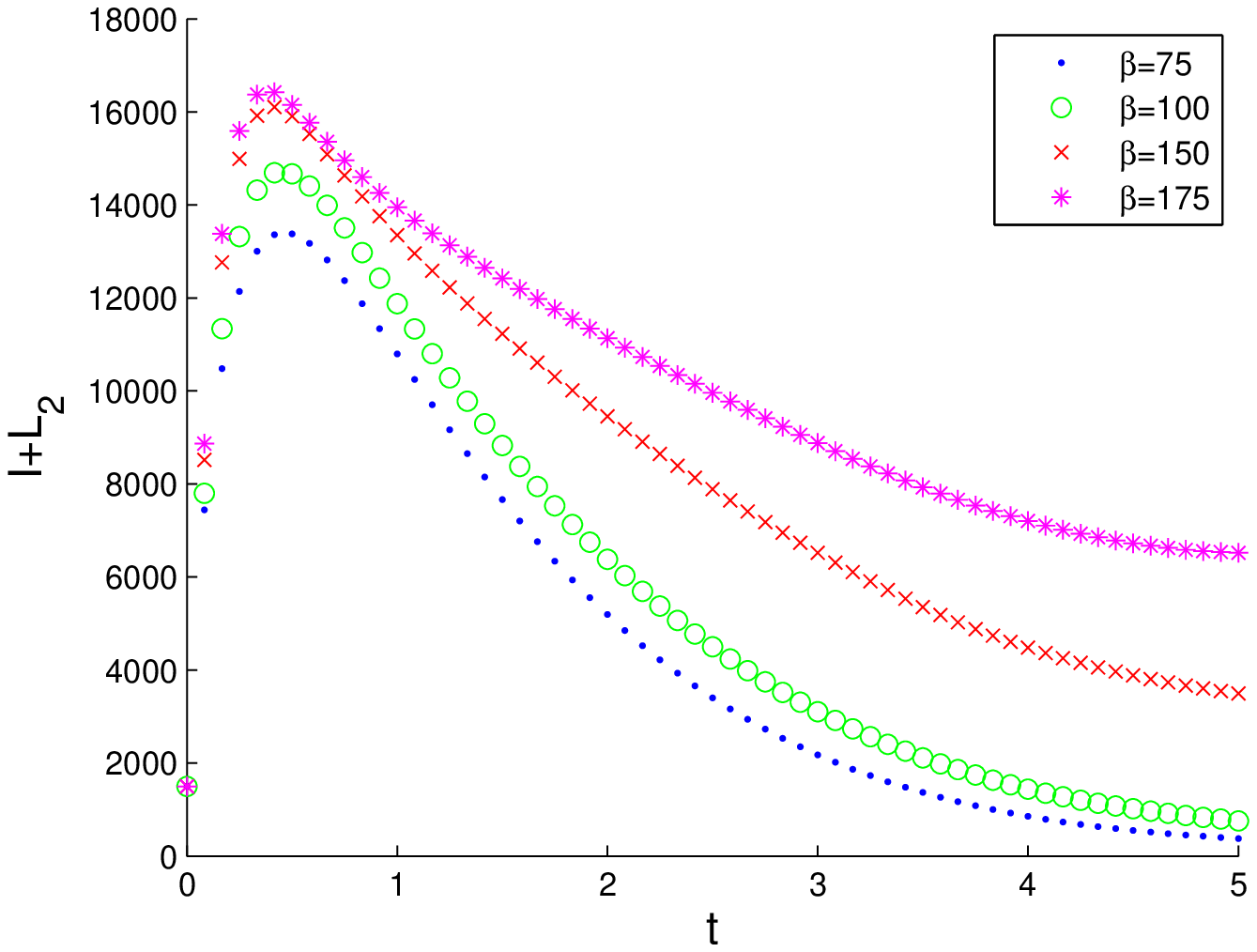,width=0.333\linewidth}\label{tb:sol:betas:k}}
\subfigure[$I$ for $f_2=10$]{\epsfig{file=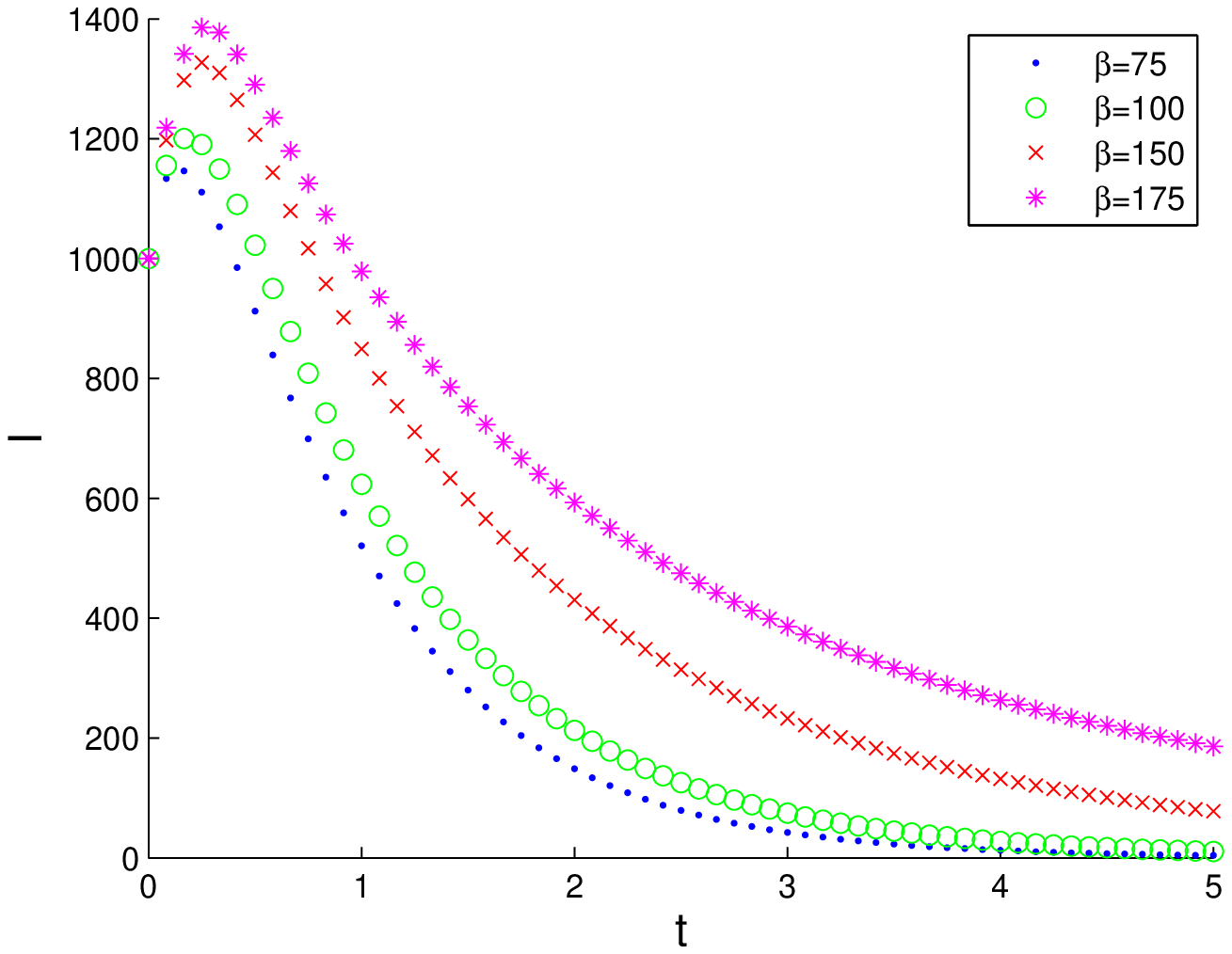,width=0.333\linewidth}\label{tb:sol:betas:l}}%
\subfigure[$L_2$ for $f_2=10$]{\epsfig{file=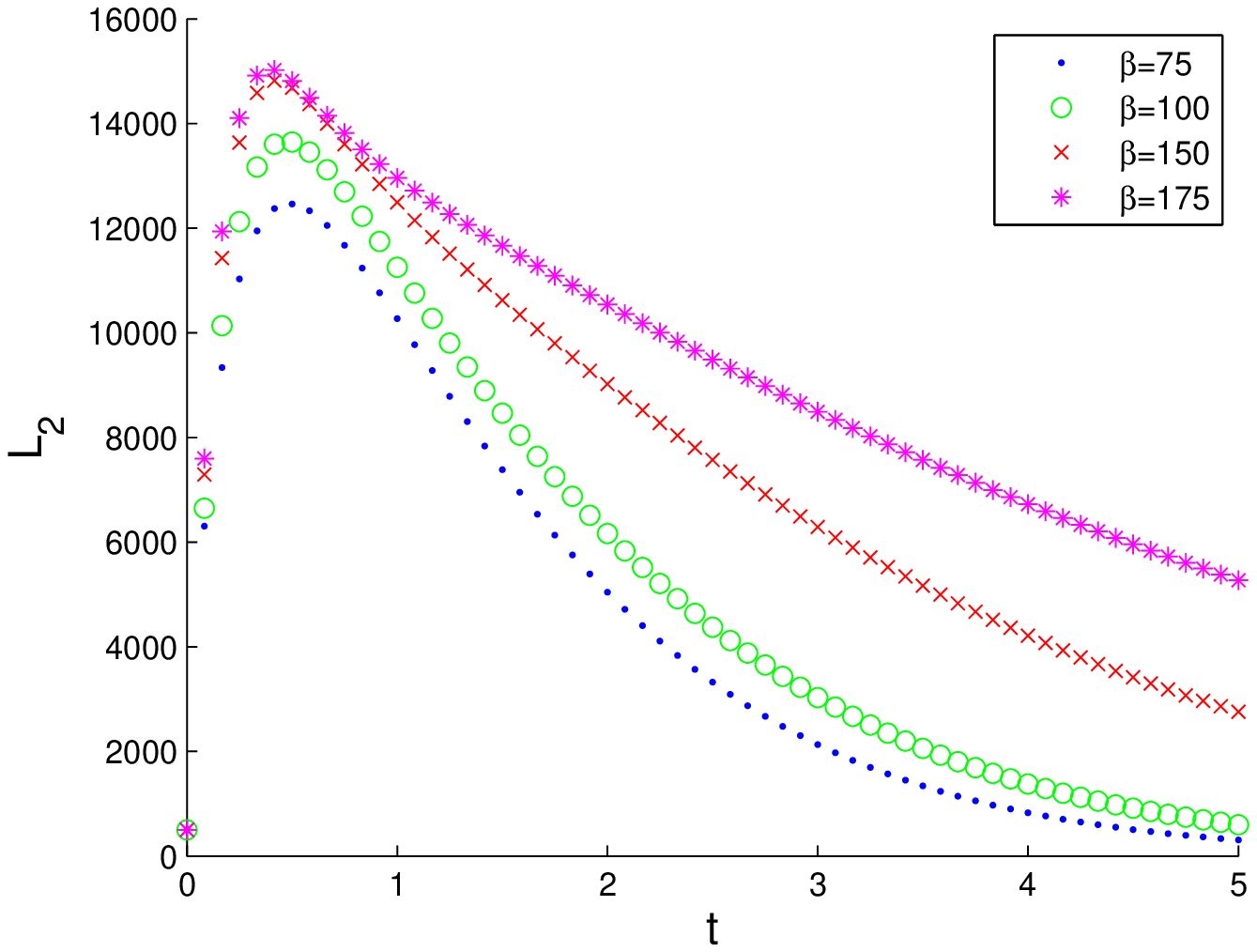,width=0.333\linewidth}\label{tb:sol:betas:m}}
\caption{Changes of controls, infectious and persistent latent individuals for different values of $\beta$.}
\label{tb:sol:betas}
\end{figure}
\begin{figure}
\centering
\includegraphics[width=0.67\textwidth]{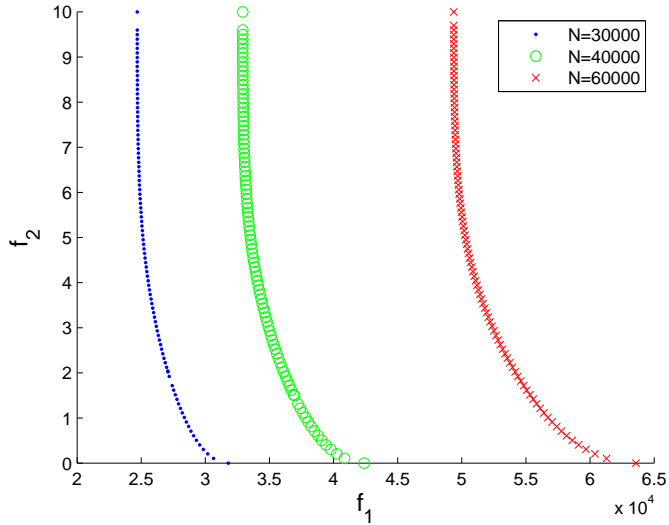}
\caption{Trade-off curves for different values of $N$.}
\label{tb:tradeoffs:pops}
\end{figure}
\begin{figure}
\centering
\subfigure[$I/N$ for $f_2=0$]{\epsfig{file=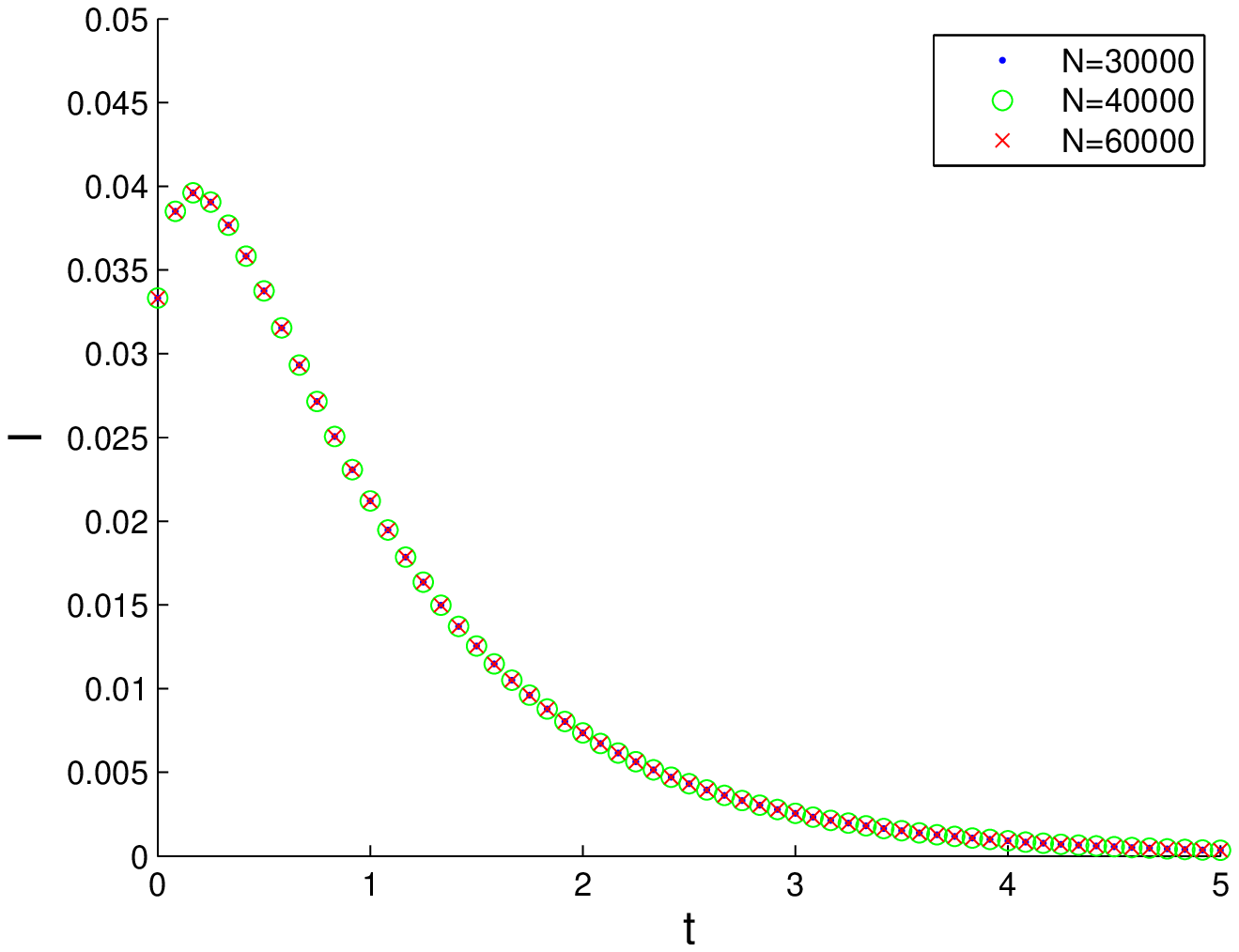,width=0.333\linewidth}\label{tb:sol:pops:a}}%
\subfigure[$L_2/N$ for $f_2=0$]{\epsfig{file=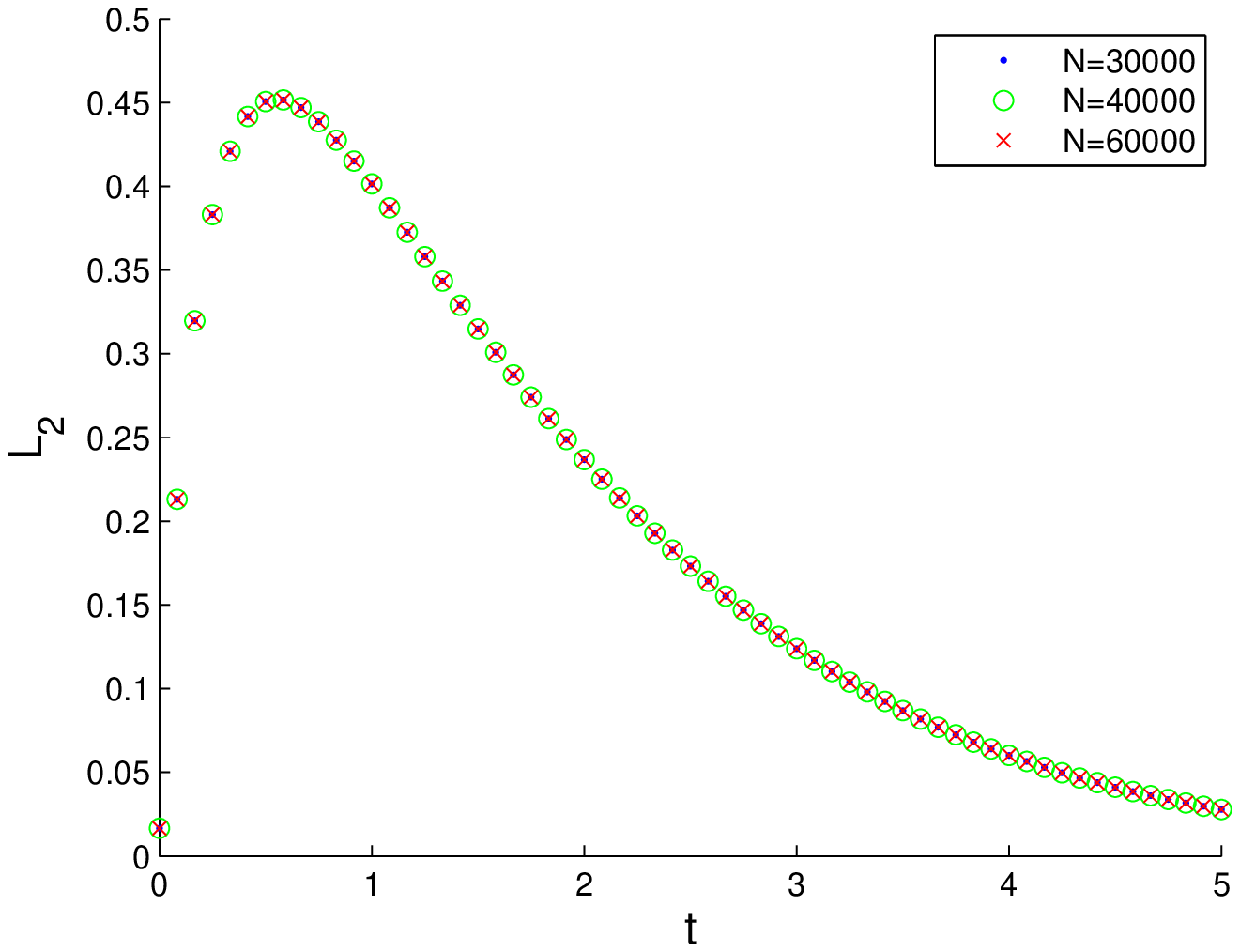,width=0.333\linewidth}\label{tb:sol:pops:b}}
\subfigure[$u_1$ for $f_2=2.5$]{\epsfig{file=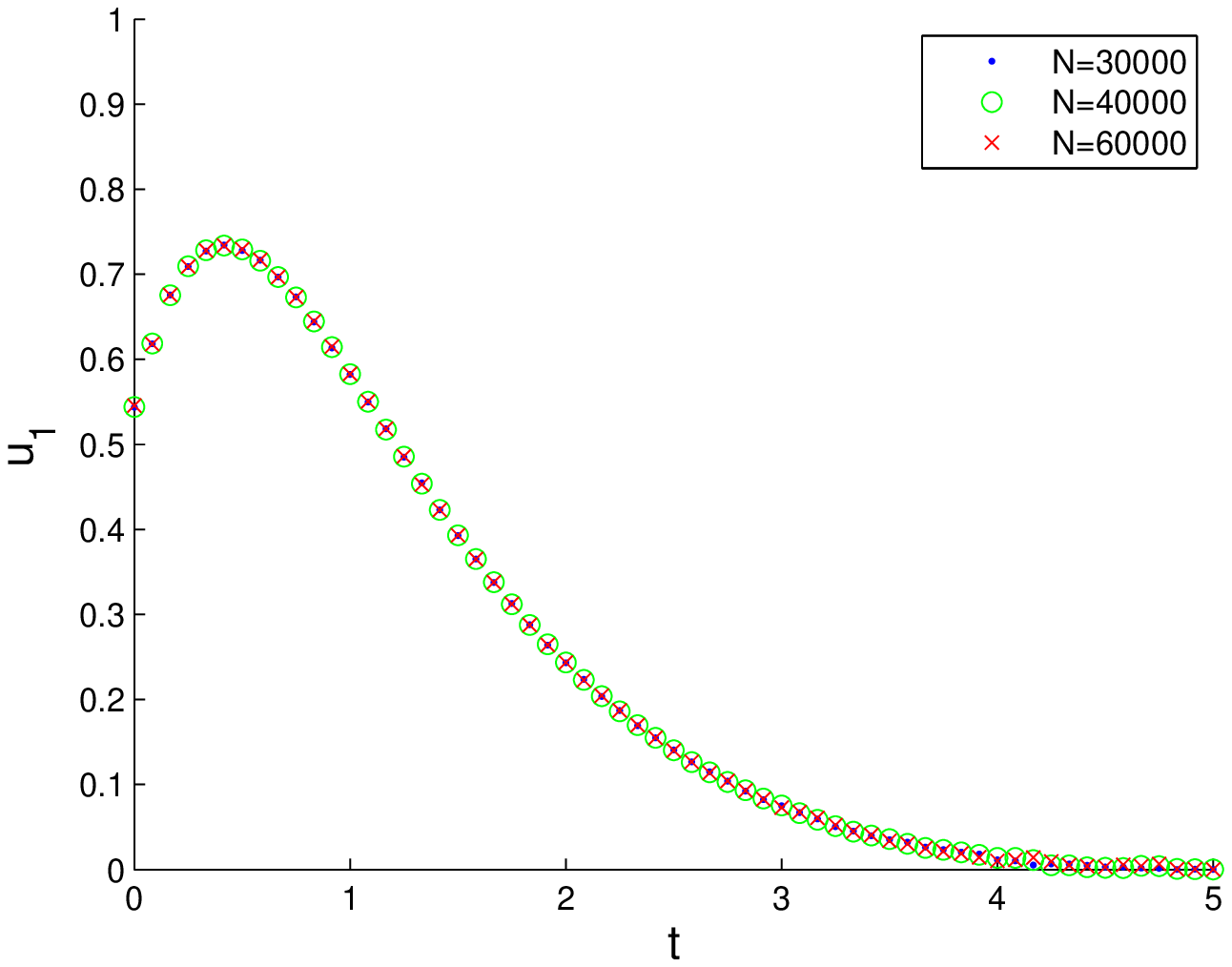,width=0.333\linewidth}\label{tb:sol:pops:c}}%
\subfigure[$u_2$ for $f_2=2.5$]{\epsfig{file=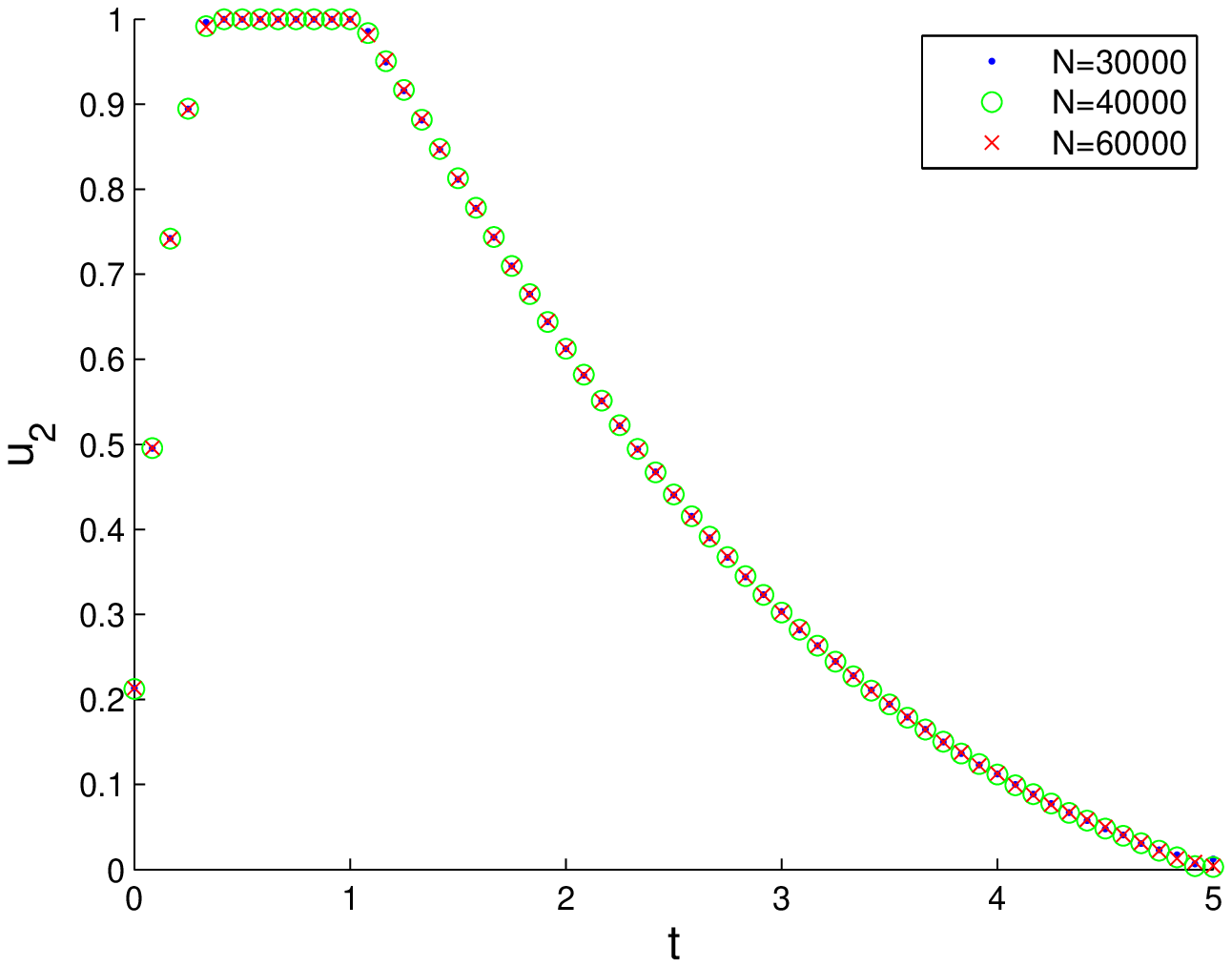,width=0.333\linewidth}\label{tb:sol:pops:d}}%
\subfigure[$(I+L_2)/N$ for $f_2=2.5$]{\epsfig{file=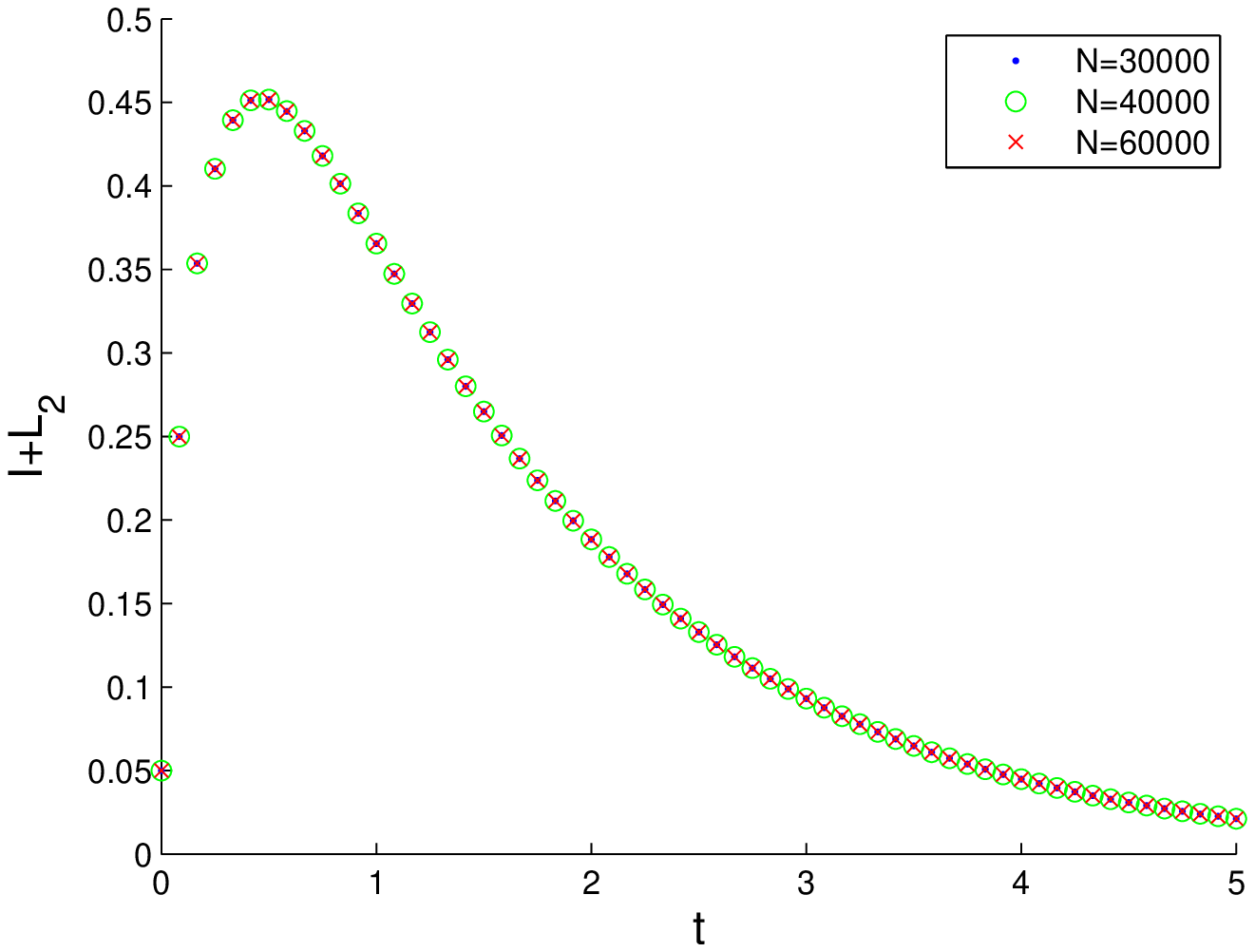,width=0.333\linewidth}\label{tb:sol:pops:e}}
\subfigure[$u_1$ for $f_2=5$]{\epsfig{file=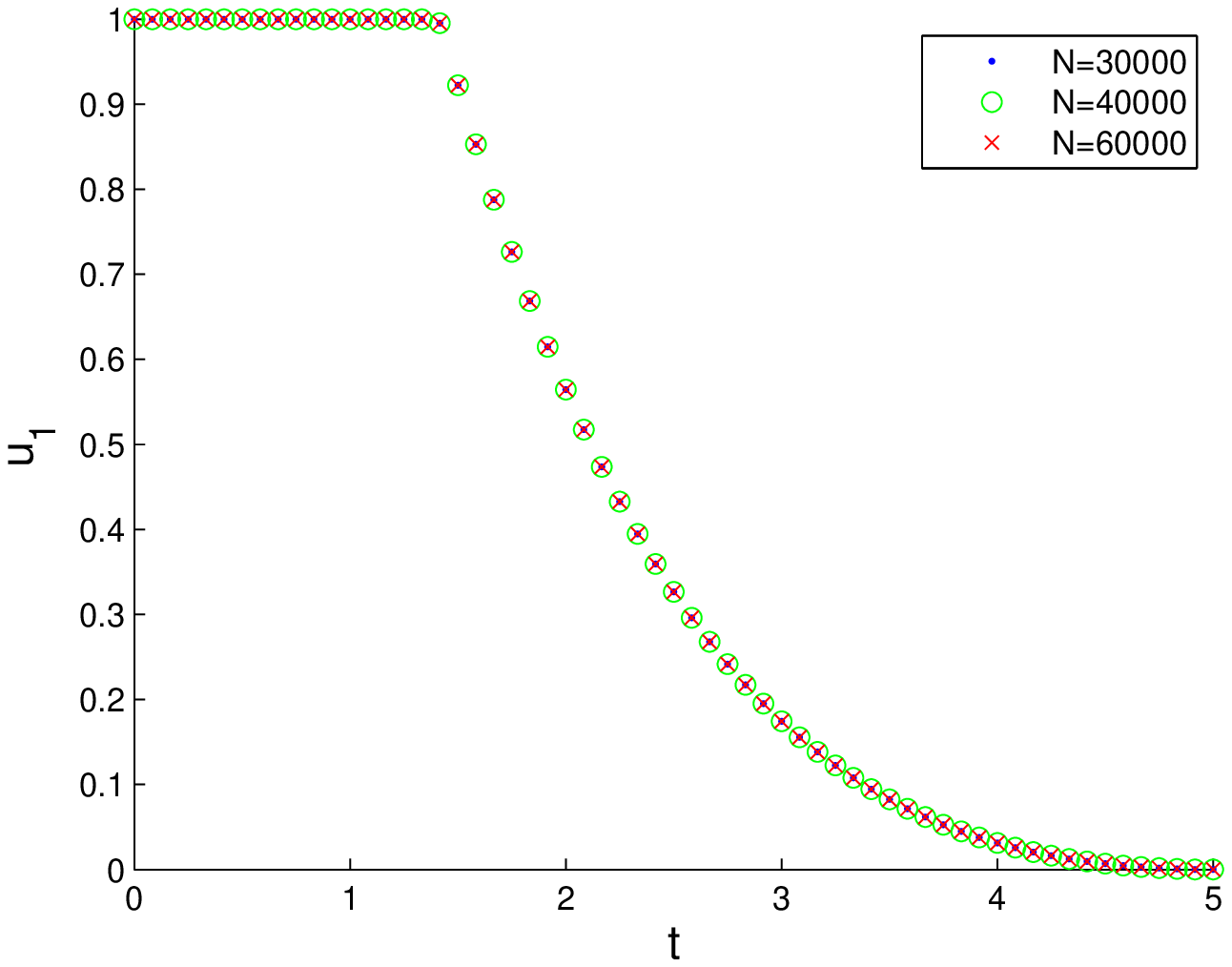,width=0.333\linewidth}\label{tb:sol:pops:f}}%
\subfigure[$u_2$ for $f_2=5$]{\epsfig{file=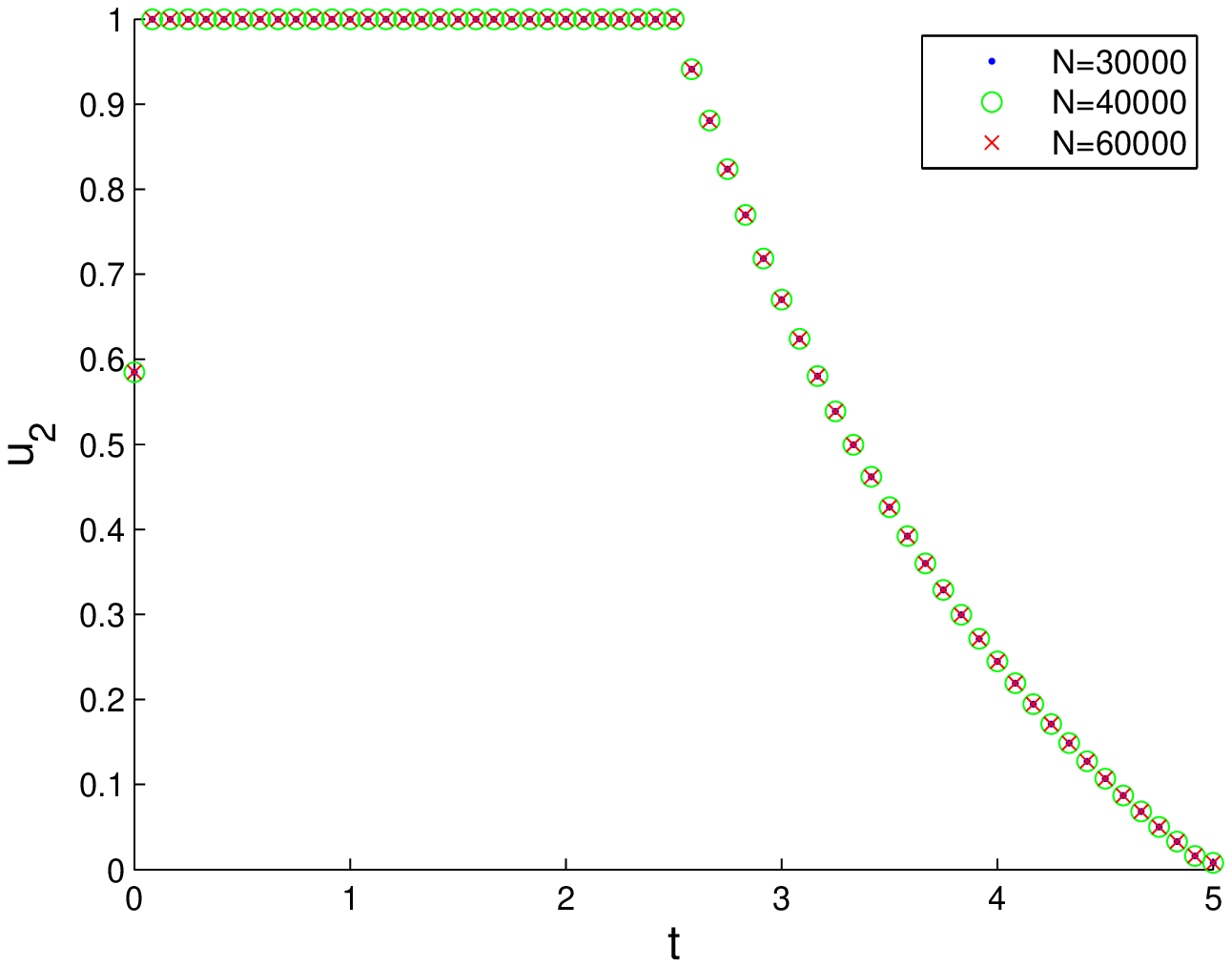,width=0.333\linewidth}\label{tb:sol:pops:g}}%
\subfigure[$(I+L_2)/N$ for $f_2=5$]{\epsfig{file=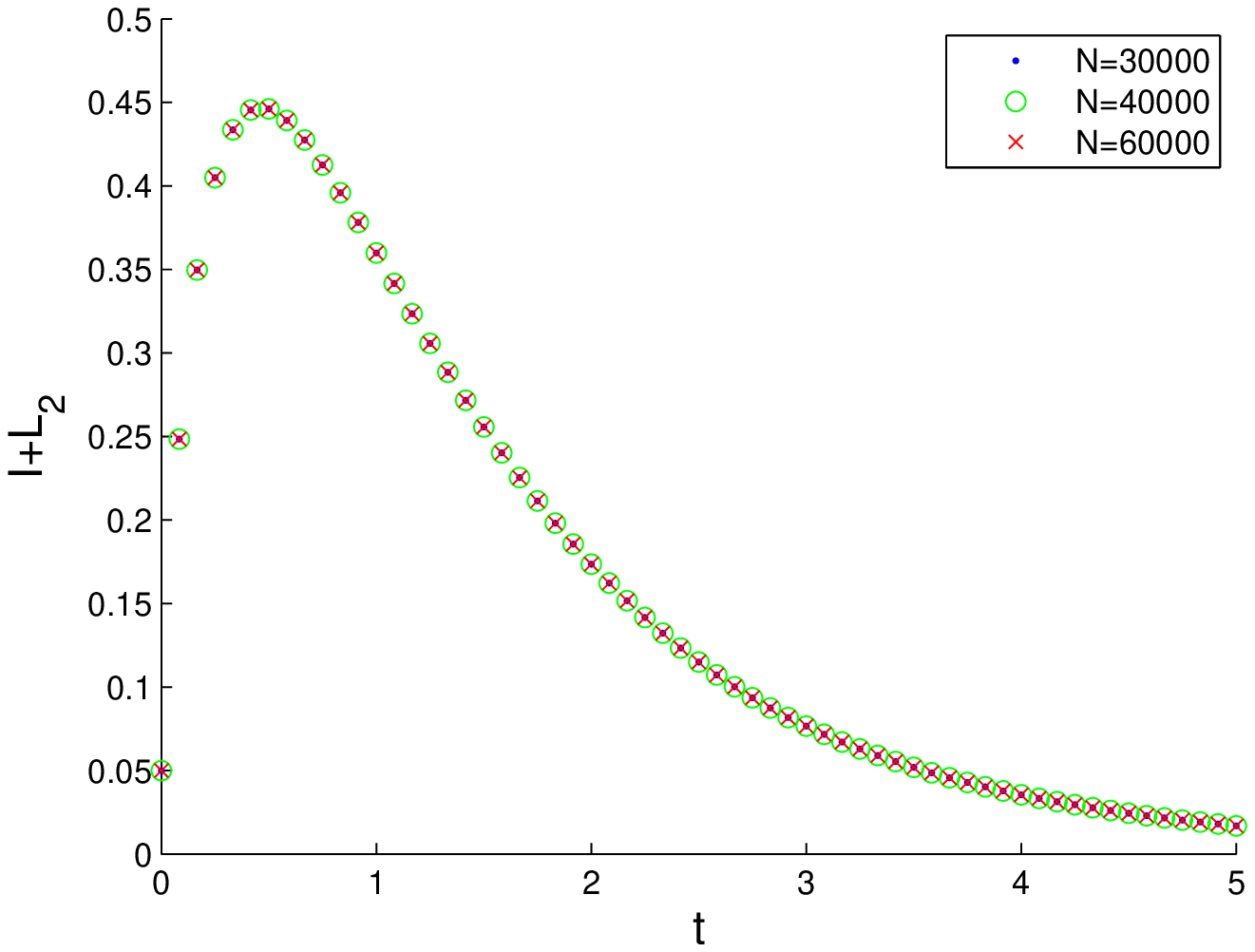,width=0.333\linewidth}\label{tb:sol:pops:h}}
\subfigure[$u_1$ for $f_2=7.5$]{\epsfig{file=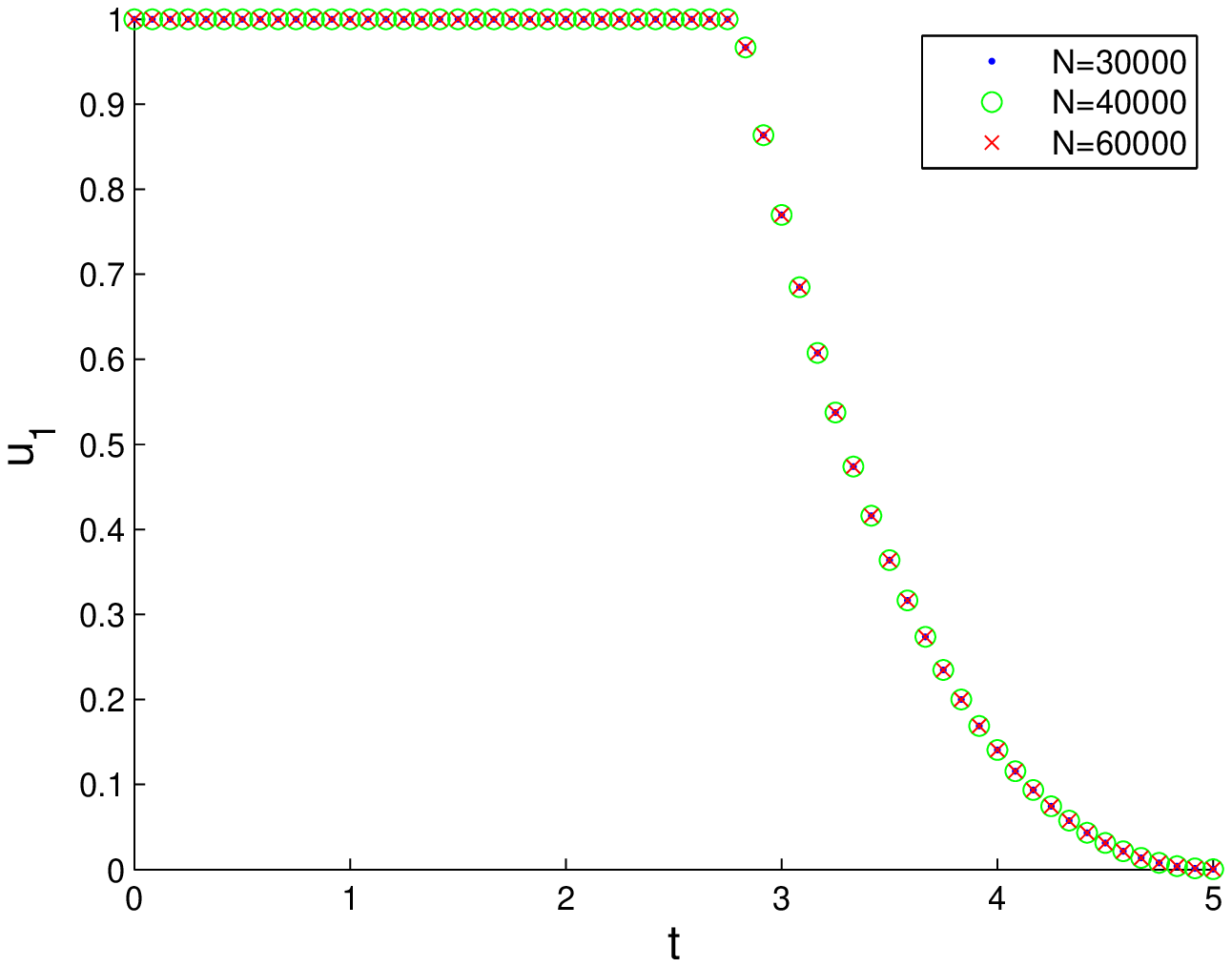,width=0.333\linewidth}\label{tb:sol:pops:i}}%
\subfigure[$u_2$ for $f_2=7.5$]{\epsfig{file=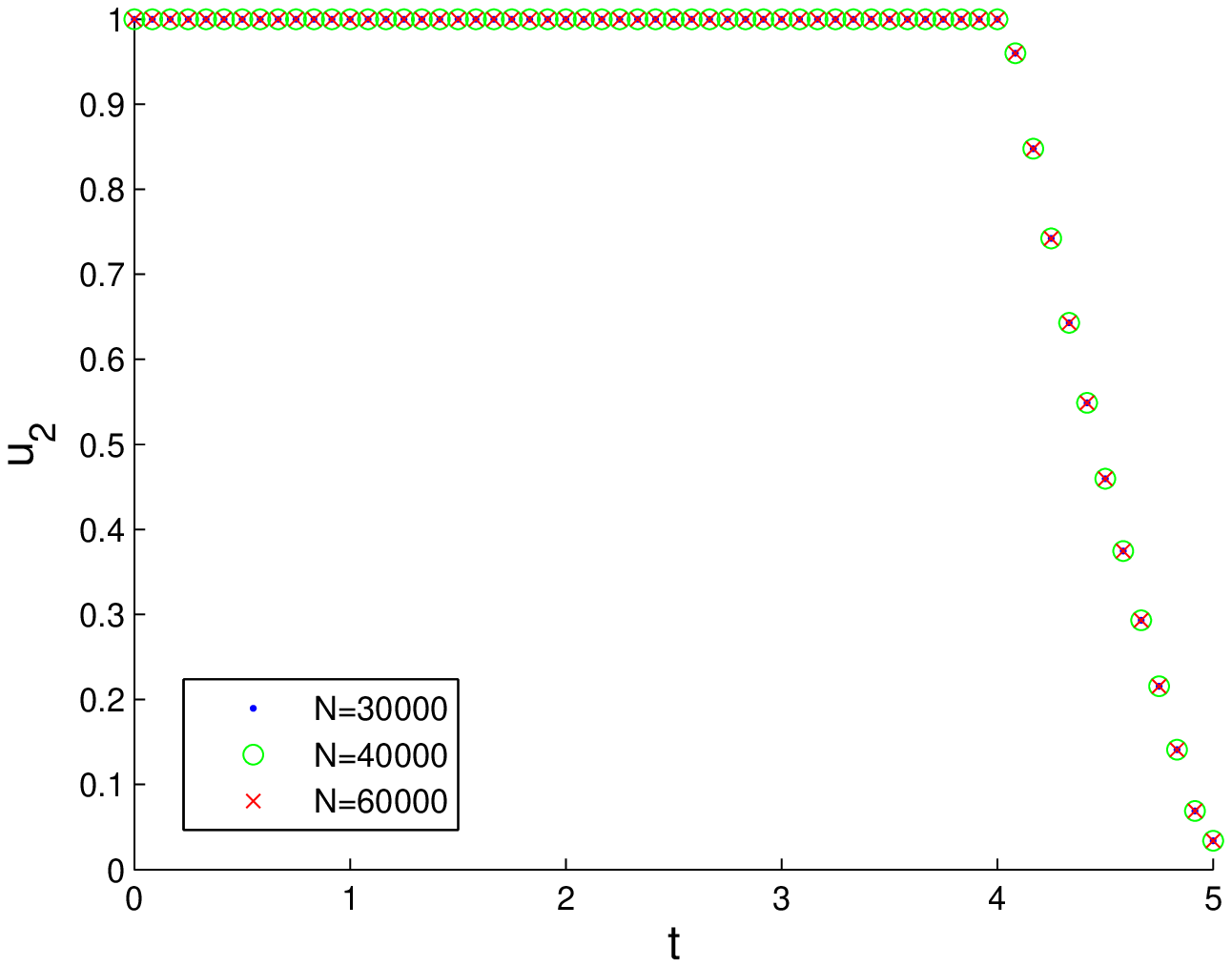,width=0.333\linewidth}\label{tb:sol:pops:j}}%
\subfigure[$(I+L_2)/N$ for $f_2=7.5$]{\epsfig{file=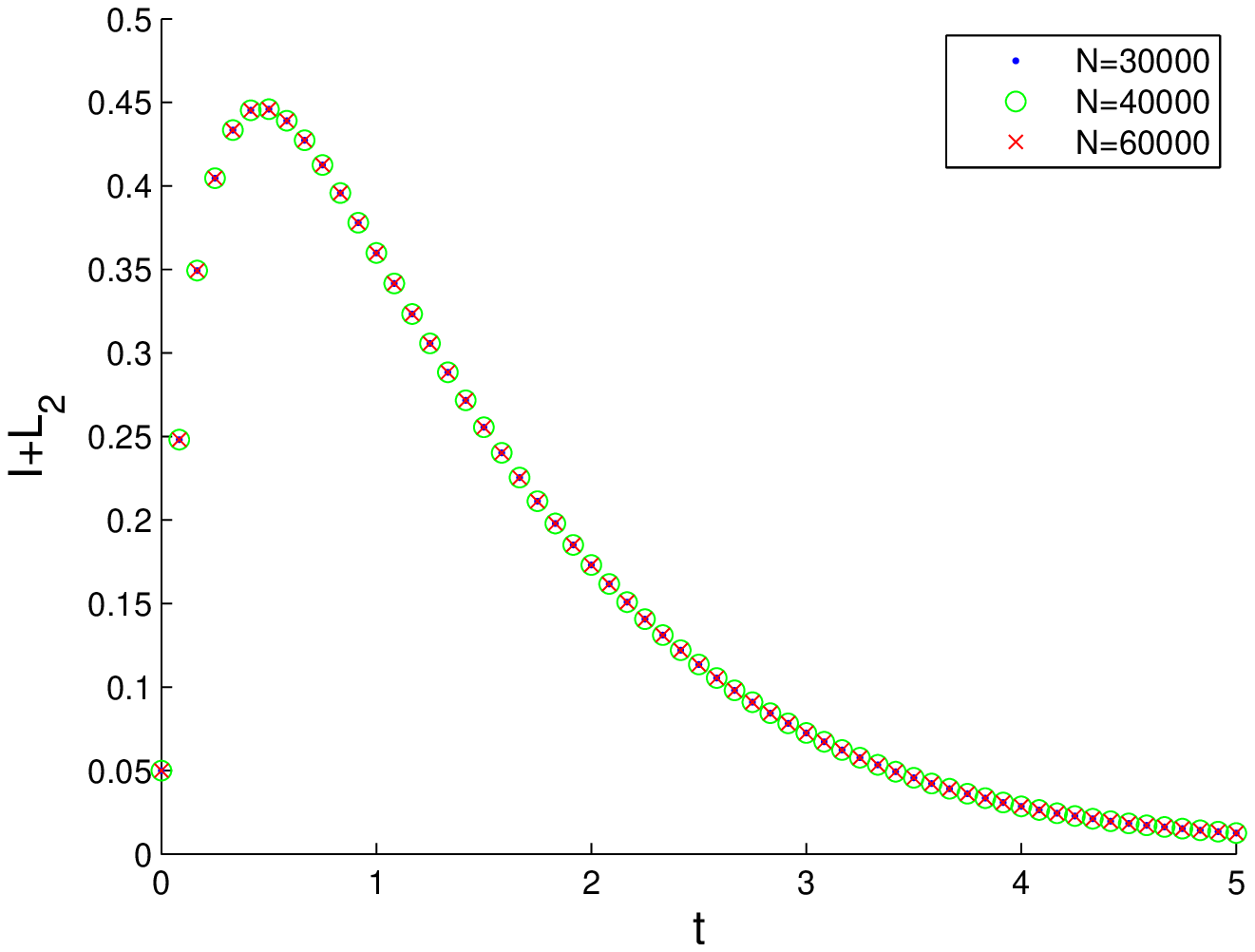,width=0.333\linewidth}\label{tb:sol:pops:k}}
\subfigure[$I/N$ for $f_2=10$]{\epsfig{file=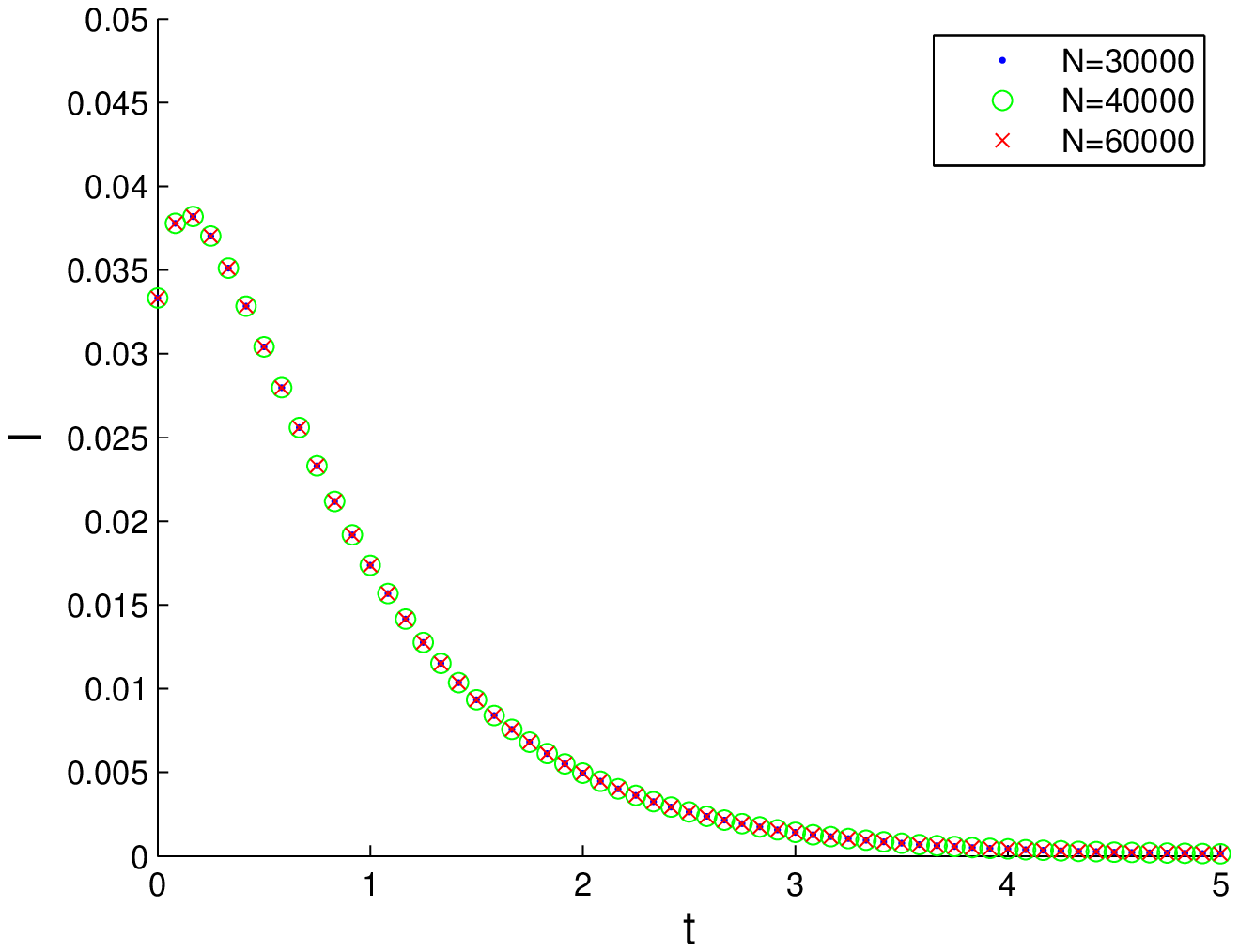,width=0.333\linewidth}\label{tb:sol:pops:l}}%
\subfigure[$L_2/N$ for $f_2=10$]{\epsfig{file=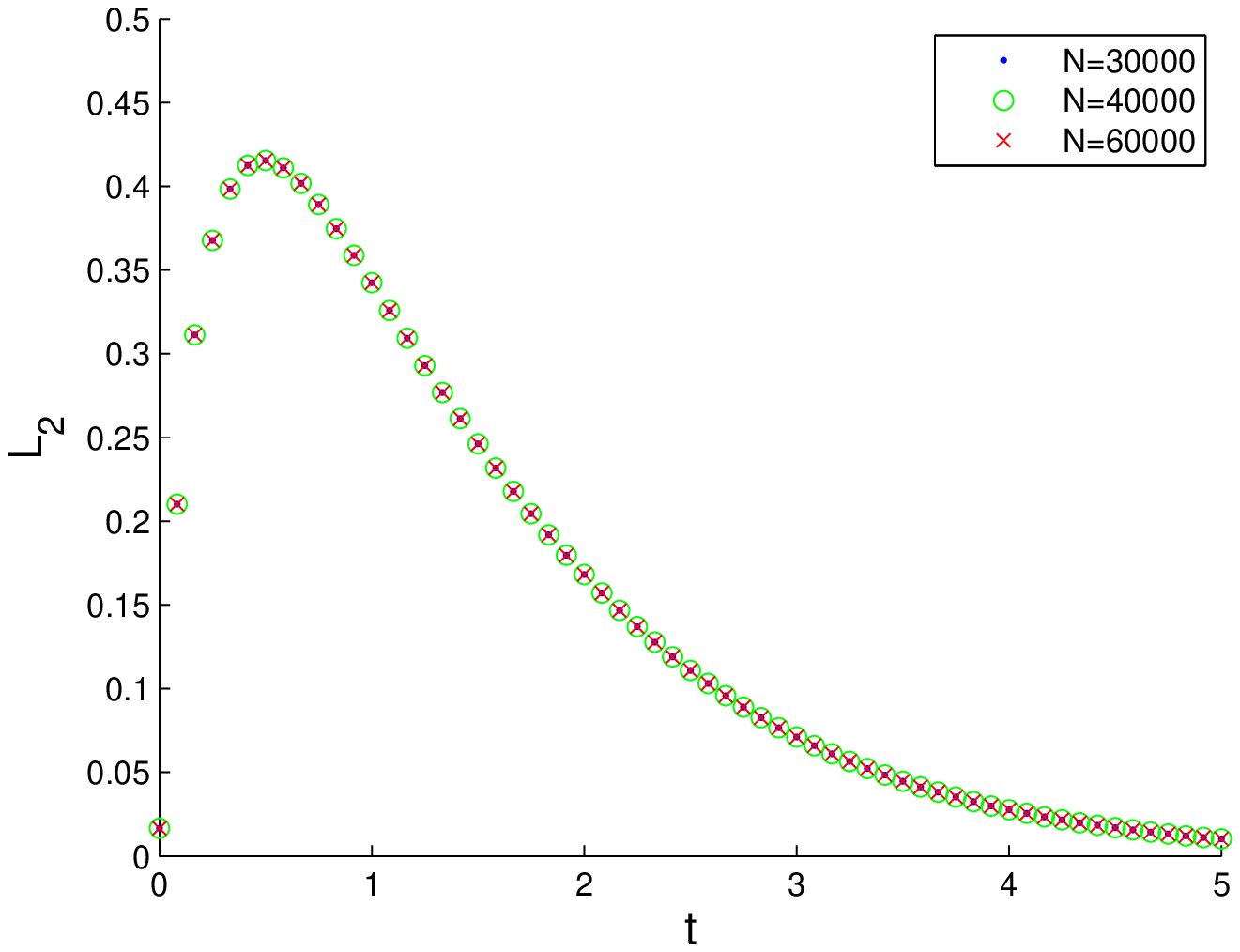,width=0.333\linewidth}\label{tb:sol:pops:m}}
\caption{Changes of controls, infectious and persistent latent individuals for different values of $N$.}
\label{tb:sol:pops}
\end{figure}
\begin{figure}
\centering
\includegraphics[width=0.67\textwidth]{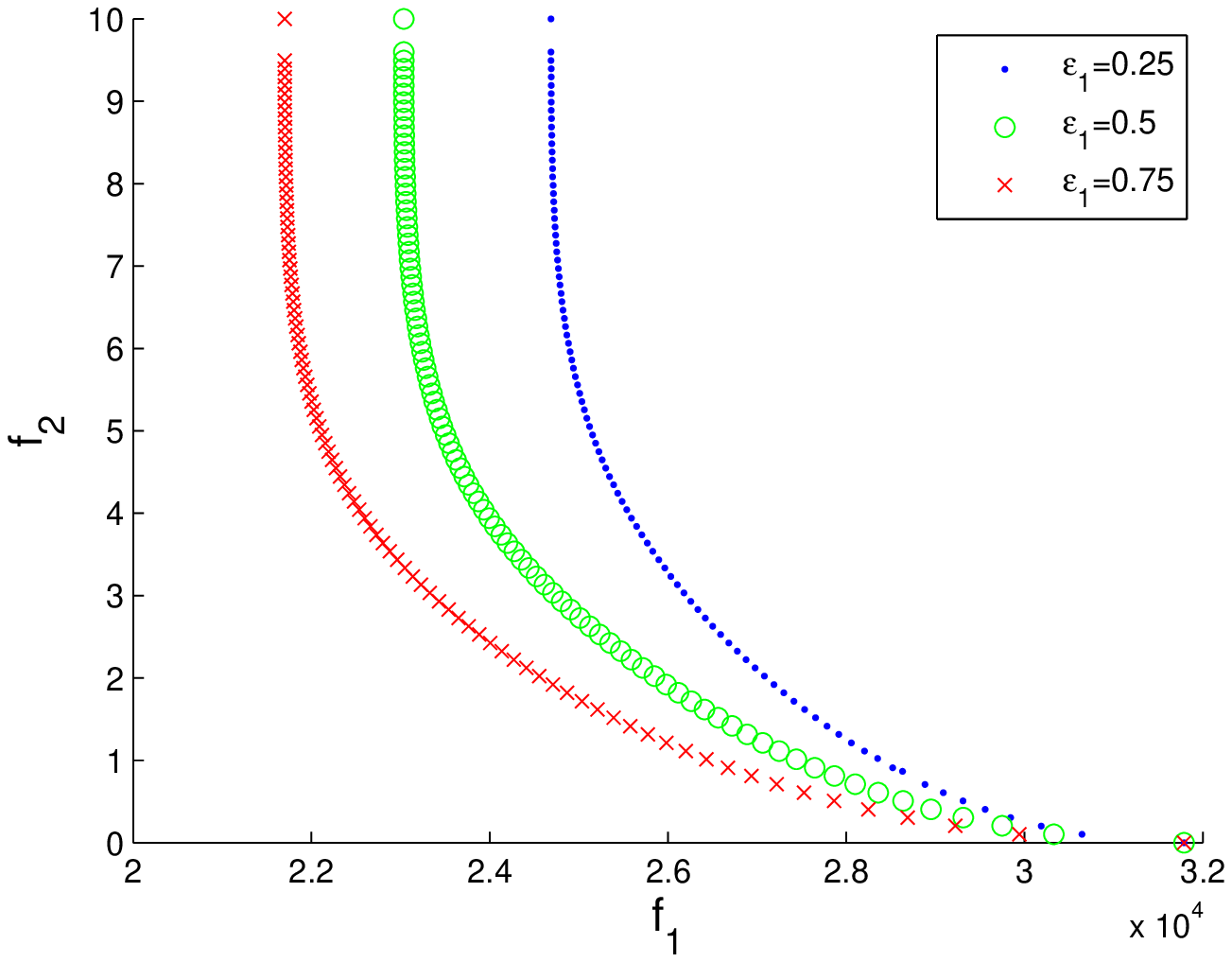}
\caption{Trade-off curves for different values of $\epsilon_1$.}
\label{tb:tradeoffs:eps1s}
\end{figure}
\begin{figure}
\centering
\subfigure[$I$ for $f_2=0$]{\epsfig{file=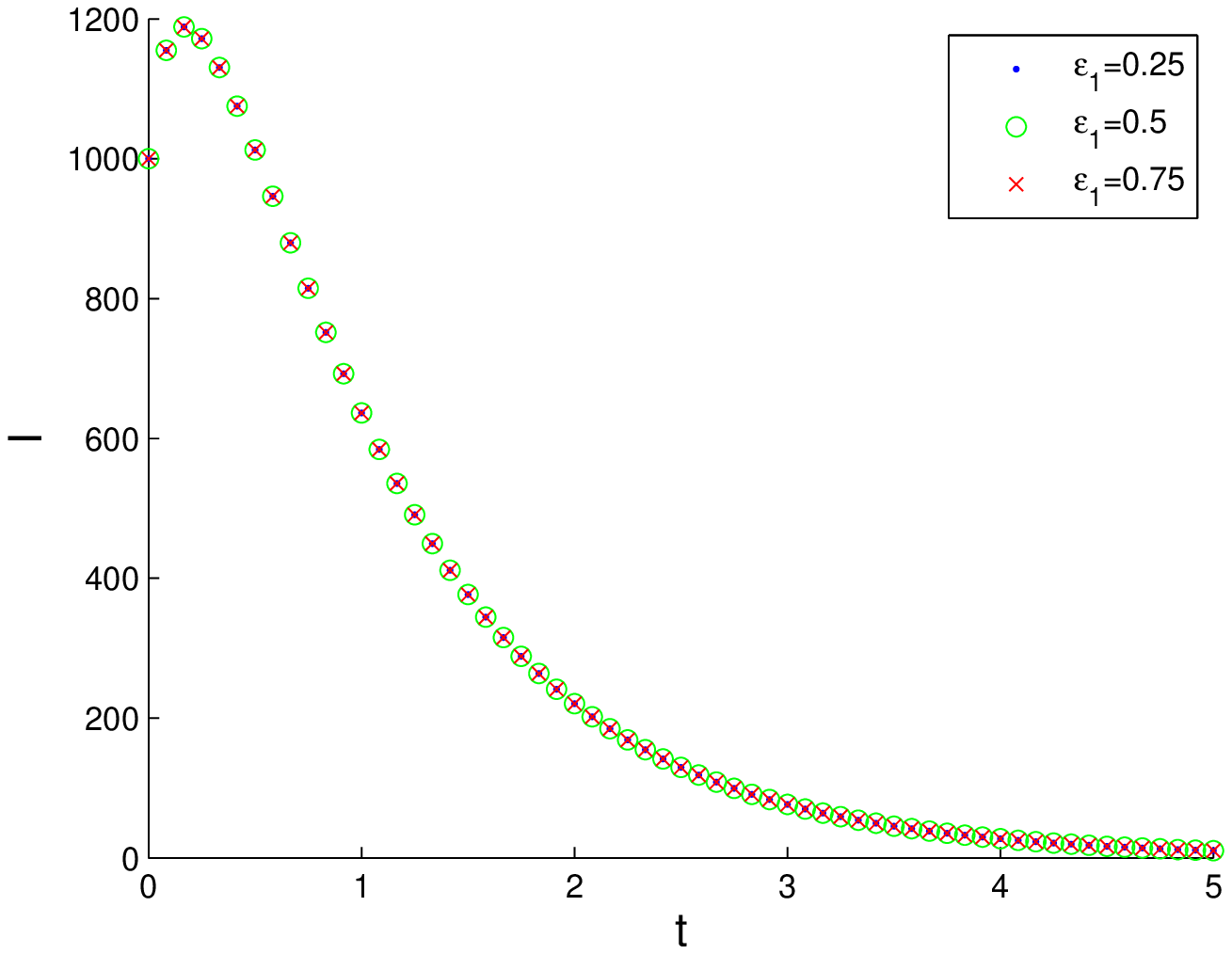,width=0.333\linewidth}\label{tb:sol:eps1s:a}}%
\subfigure[$L_2$ for $f_2=0$]{\epsfig{file=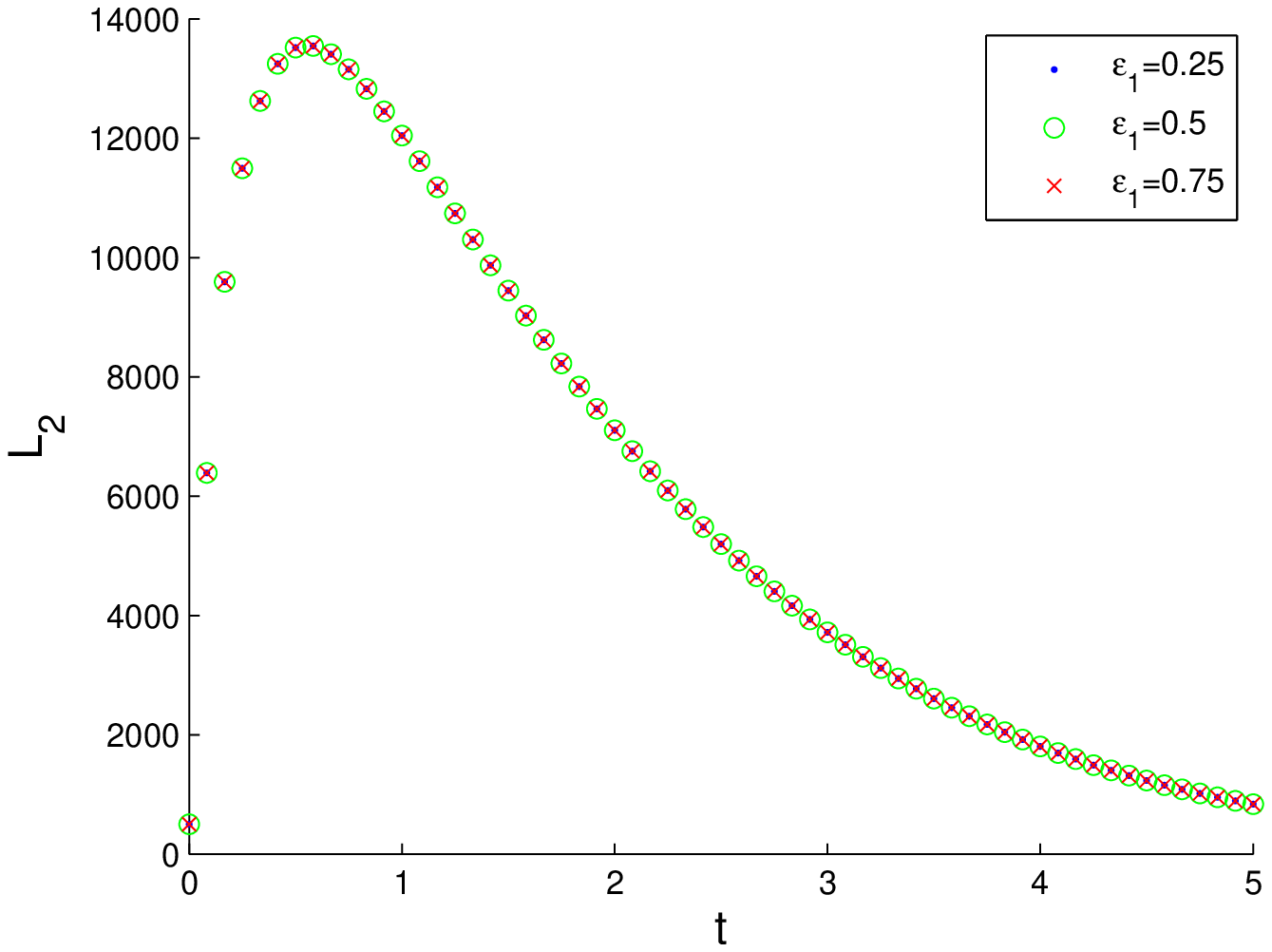,width=0.333\linewidth}\label{tb:sol:eps1s:b}}
\subfigure[$u_1$ for $f_2=2.5$]{\epsfig{file=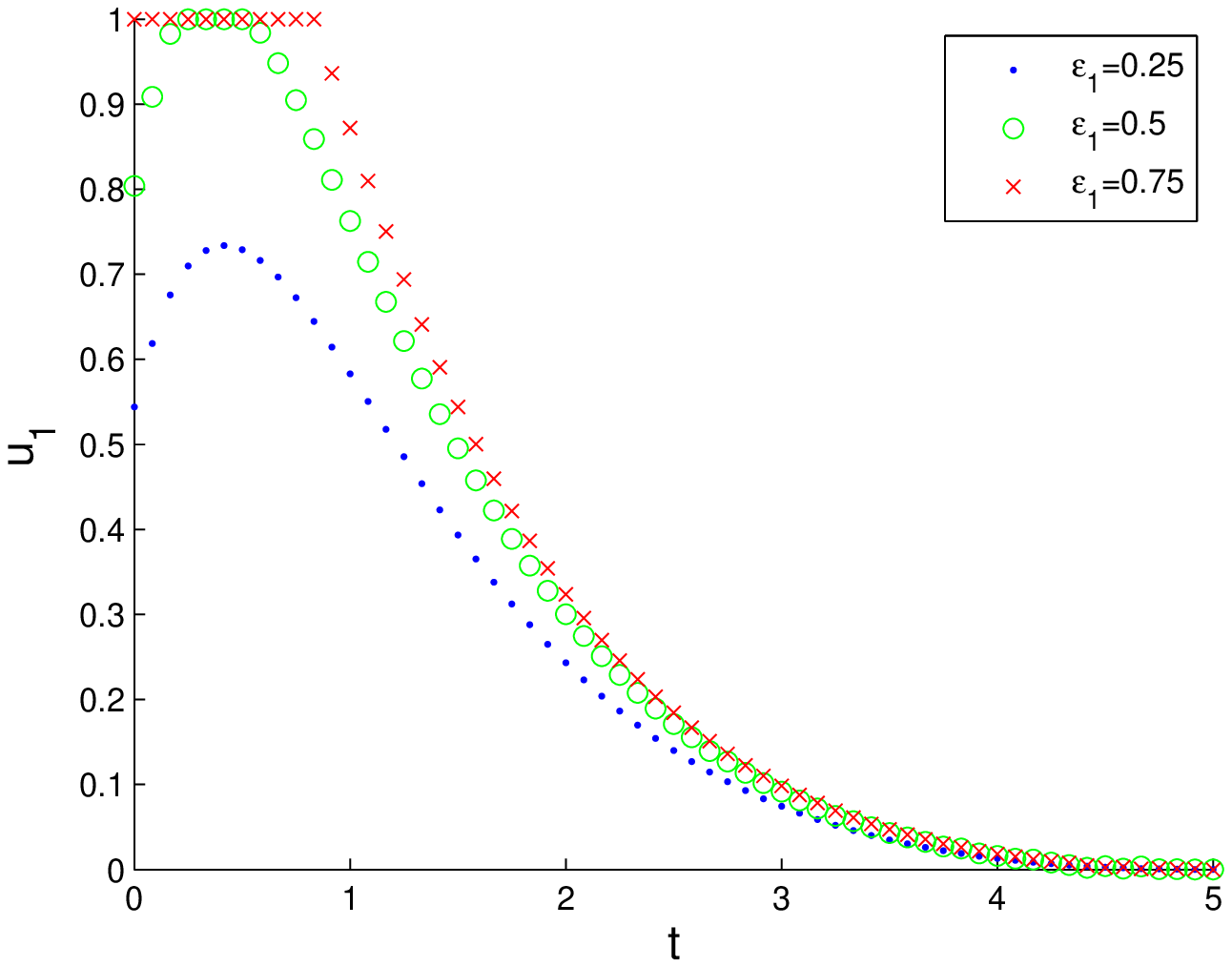,width=0.333\linewidth}\label{tb:sol:eps1s:c}}%
\subfigure[$u_2$ for $f_2=2.5$]{\epsfig{file=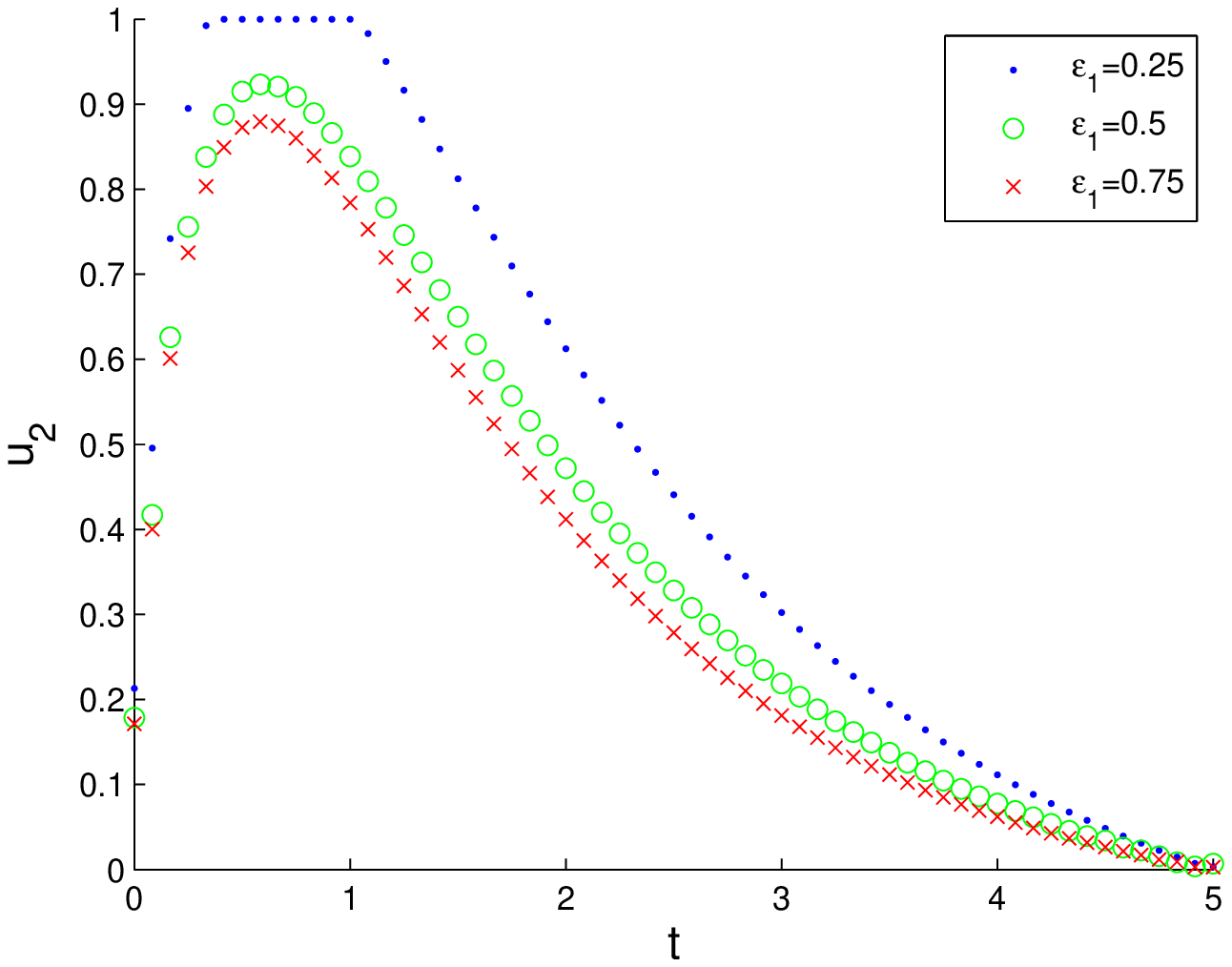,width=0.333\linewidth}\label{tb:sol:eps1s:d}}%
\subfigure[$I+L_2$ for $f_2=2.5$]{\epsfig{file=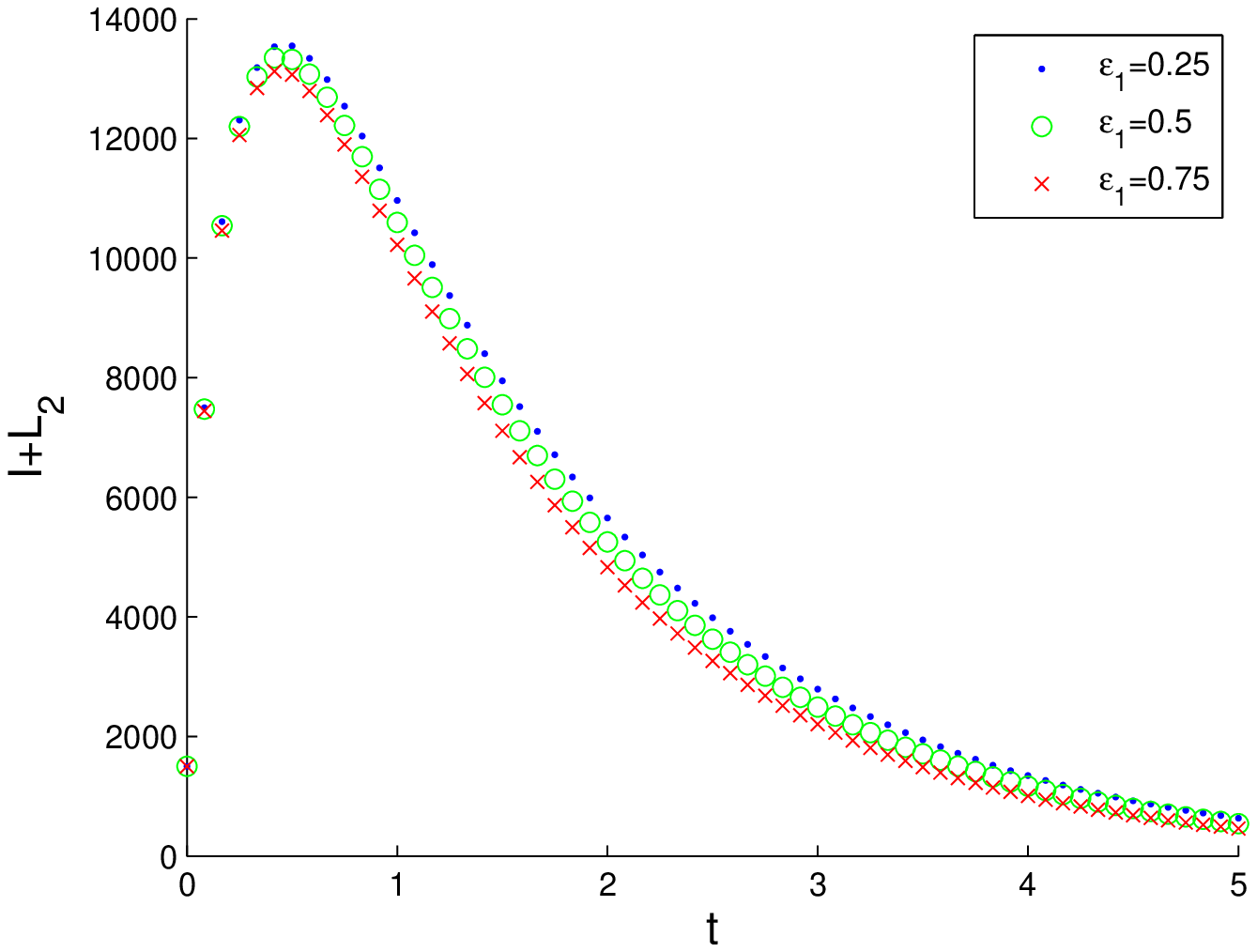,width=0.333\linewidth}\label{tb:sol:eps1s:e}}
\subfigure[$u_1$ for $f_2=5$]{\epsfig{file=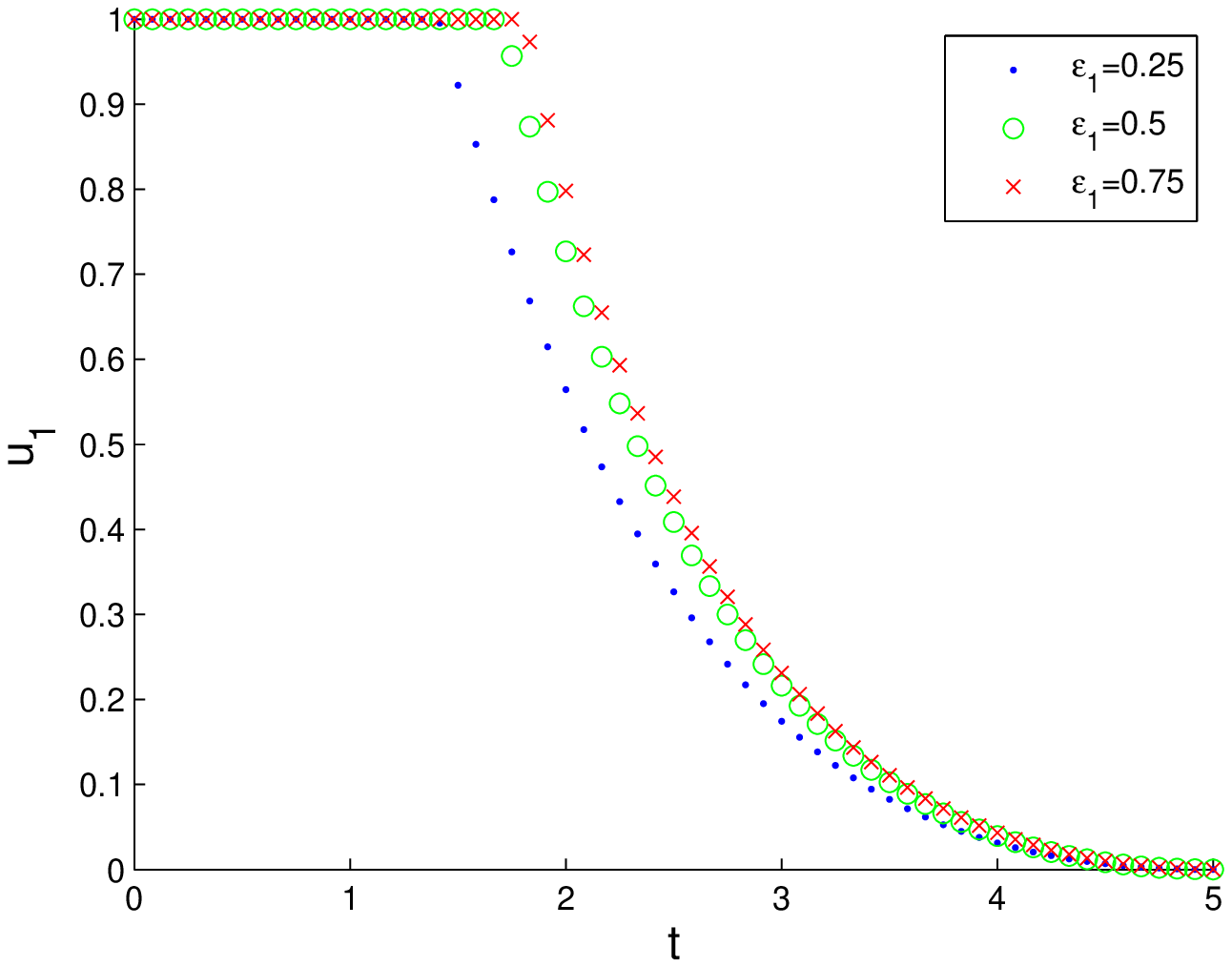,width=0.333\linewidth}\label{tb:sol:eps1s:f}}%
\subfigure[$u_2$ for $f_2=5$]{\epsfig{file=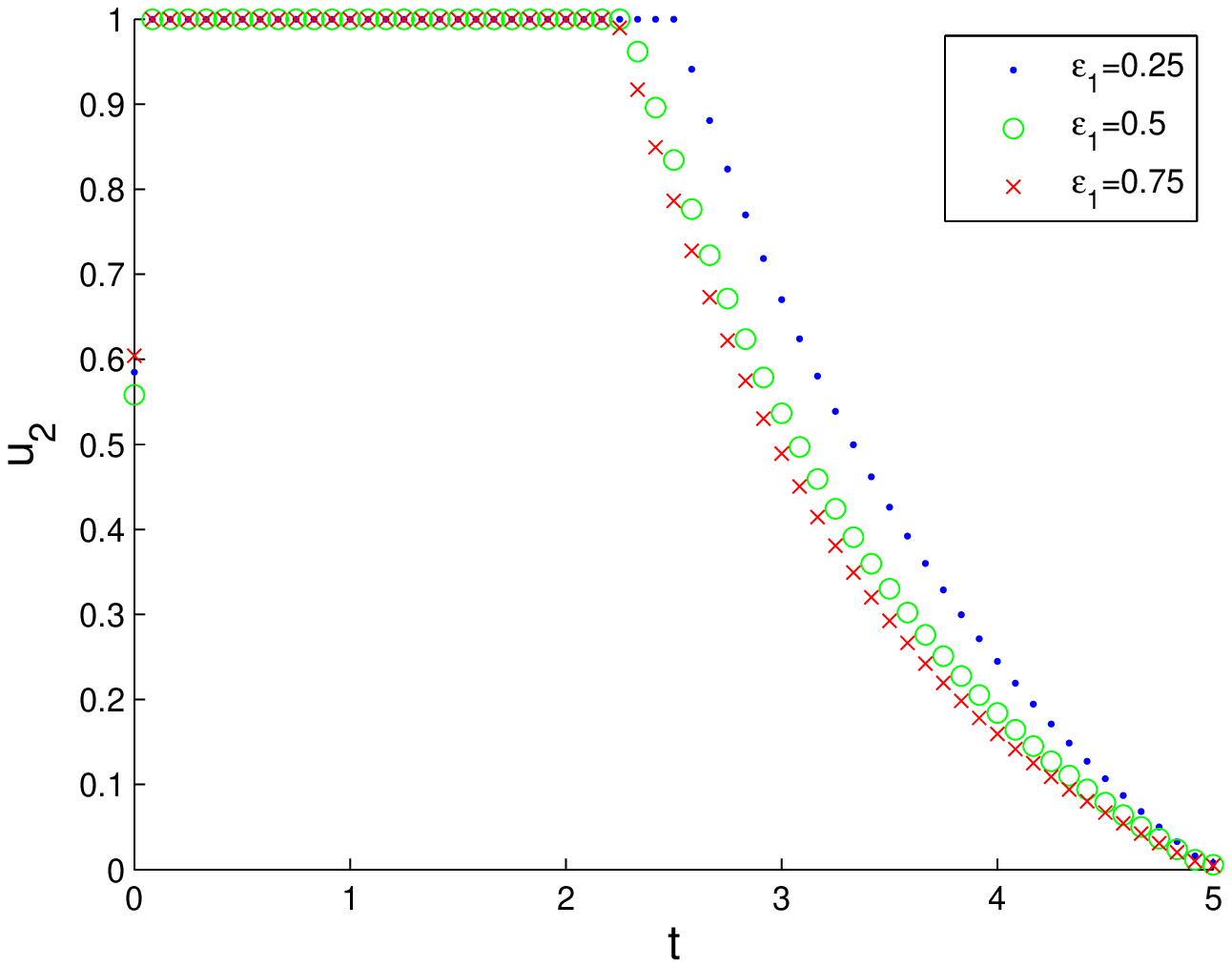,width=0.333\linewidth}\label{tb:sol:eps1s:g}}%
\subfigure[$I+L_2$ for $f_2=5$]{\epsfig{file=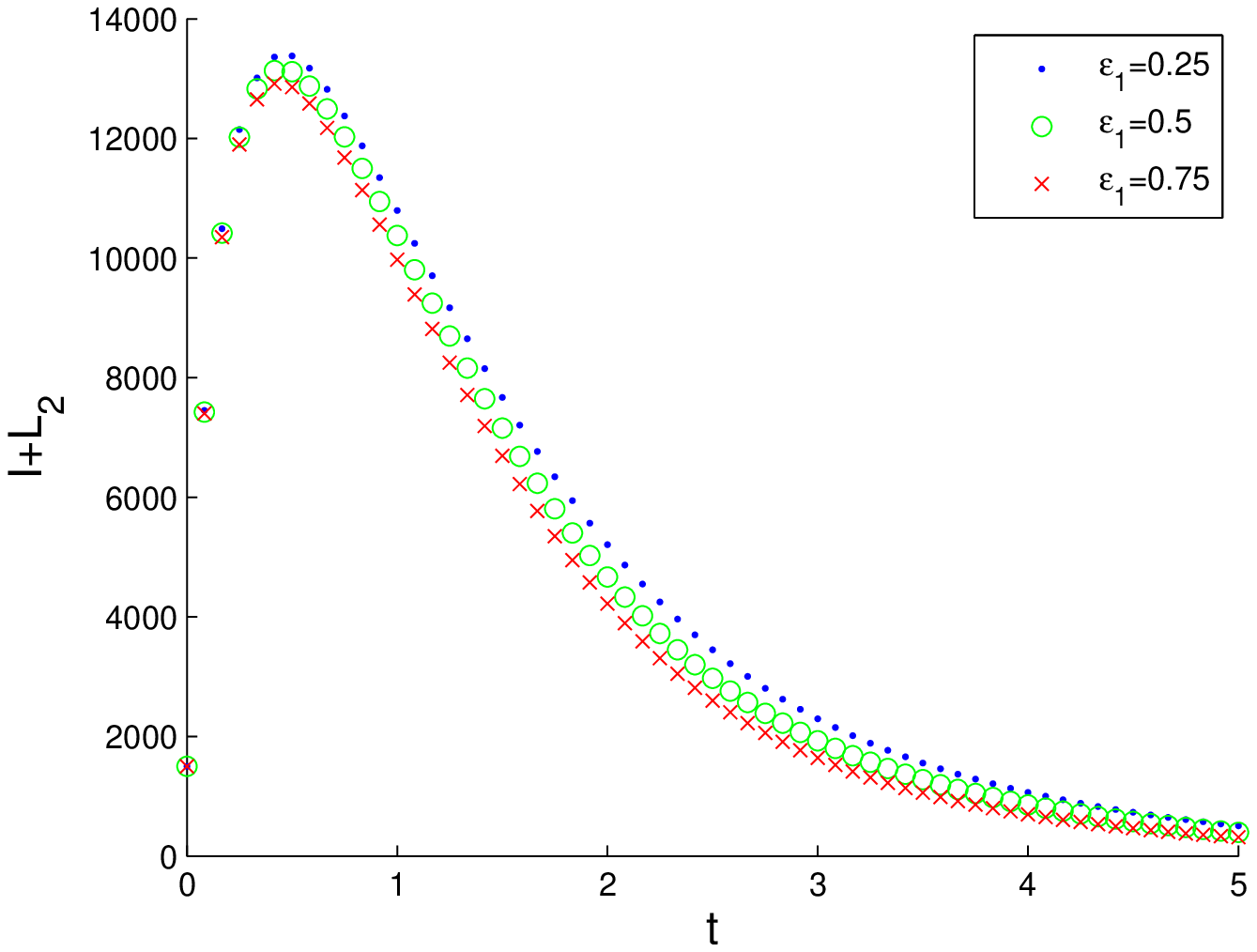,width=0.333\linewidth}\label{tb:sol:eps1s:h}}
\subfigure[$u_1$ for $f_2=7.5$]{\epsfig{file=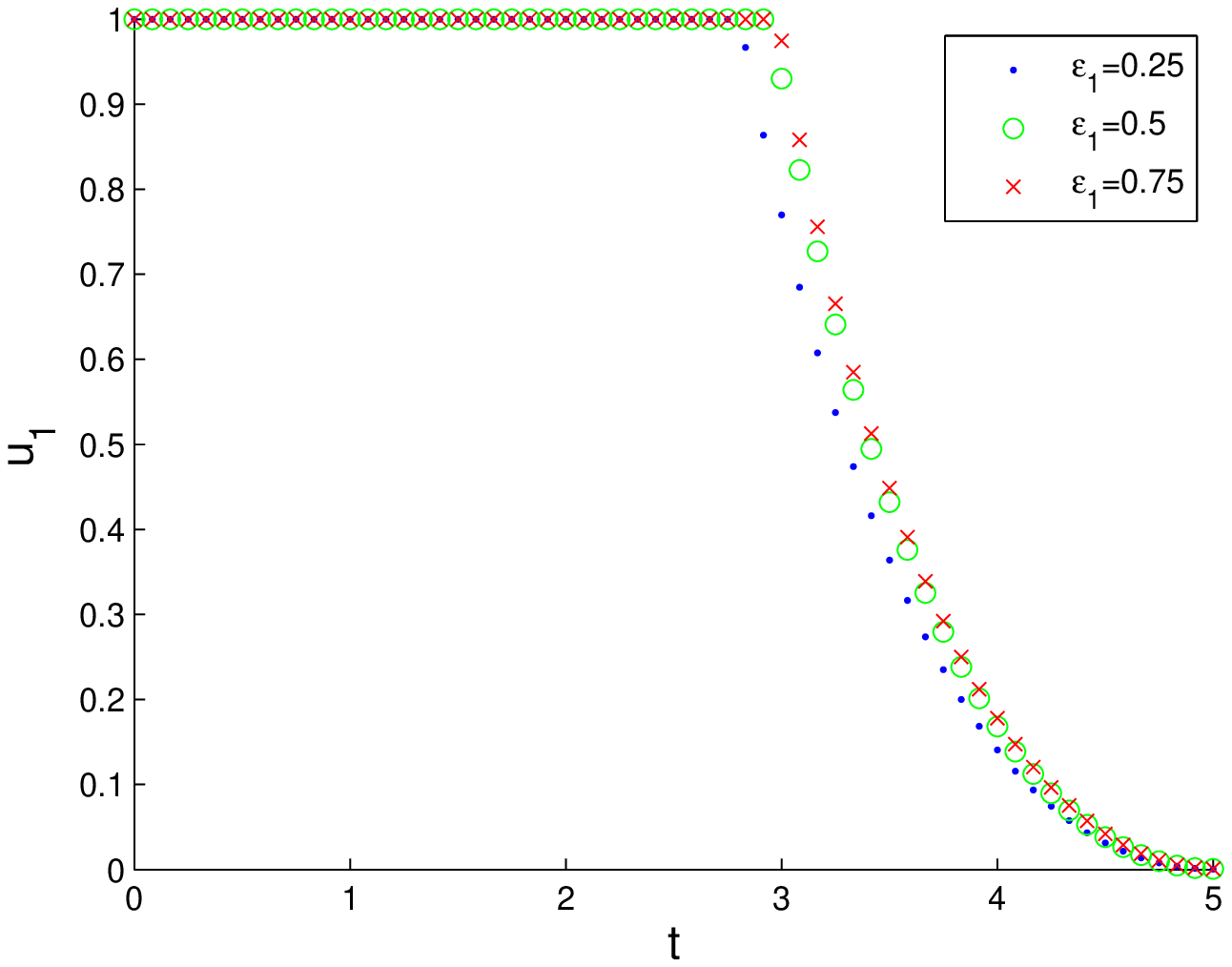,width=0.333\linewidth}\label{tb:sol:eps1s:i}}%
\subfigure[$u_2$ for $f_2=7.5$]{\epsfig{file=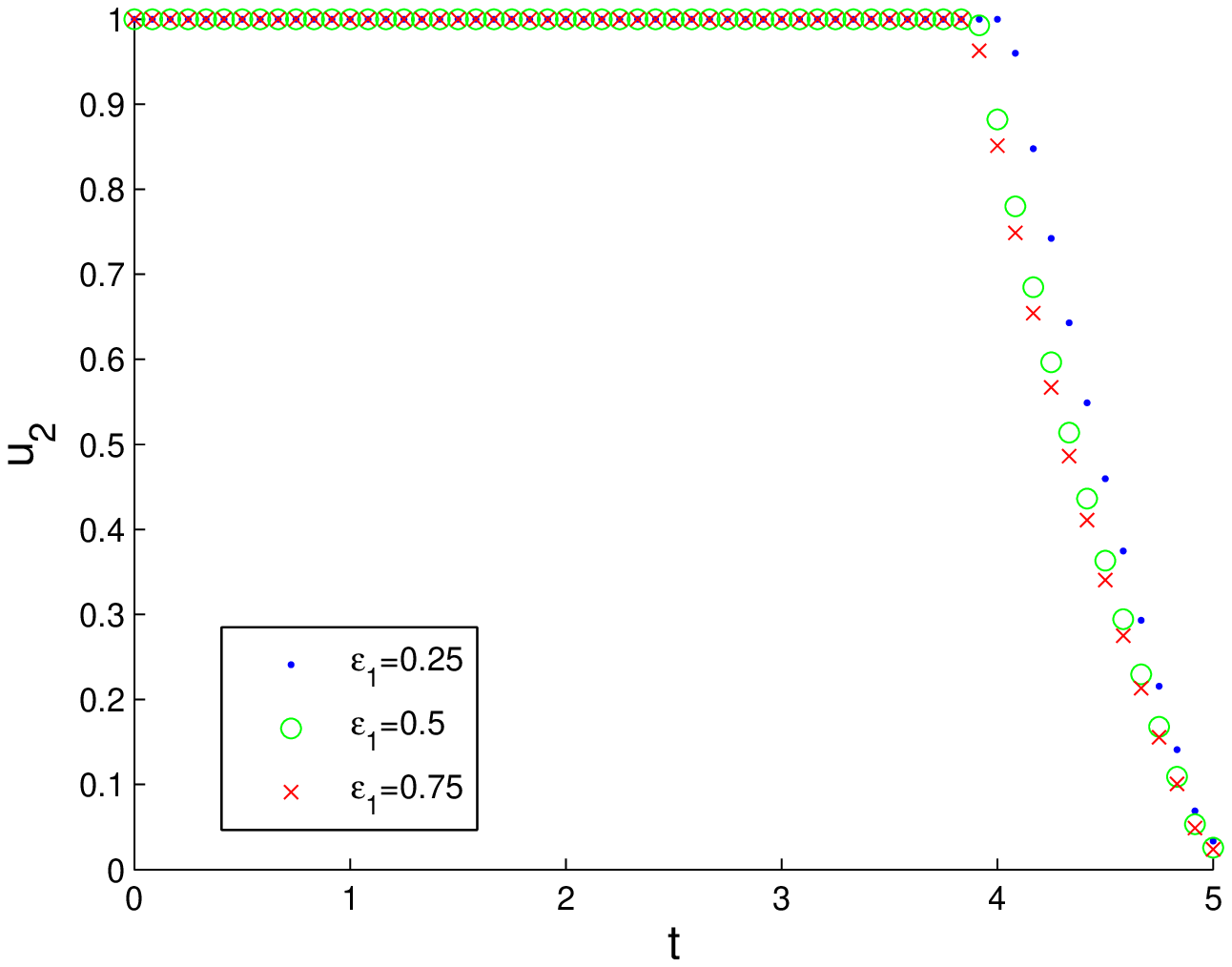,width=0.333\linewidth}\label{tb:sol:eps1s:j}}%
\subfigure[$I+L_2$ for $f_2=7.5$]{\epsfig{file=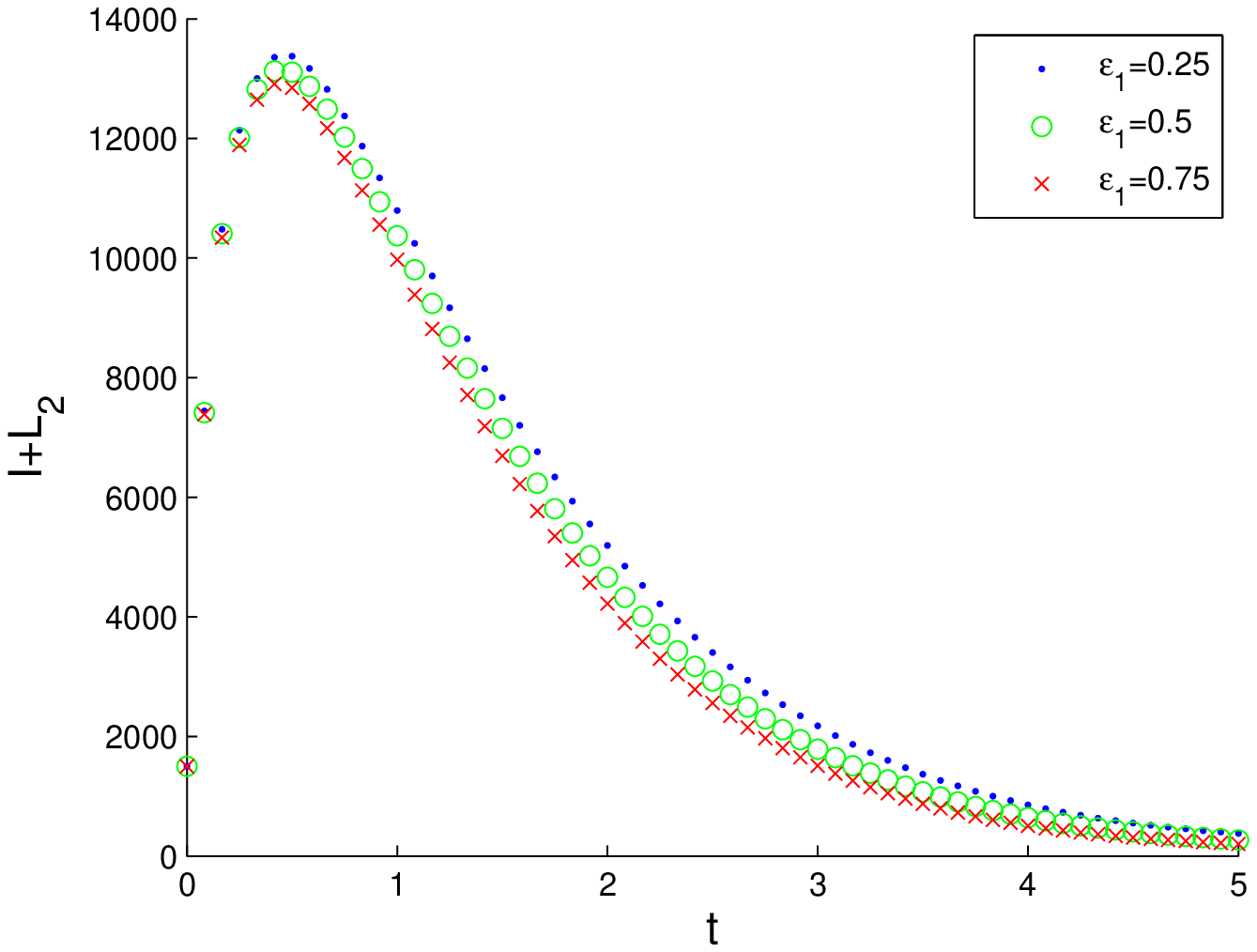,width=0.333\linewidth}\label{tb:sol:eps1s:k}}
\subfigure[$I$ for $f_2=10$]{\epsfig{file=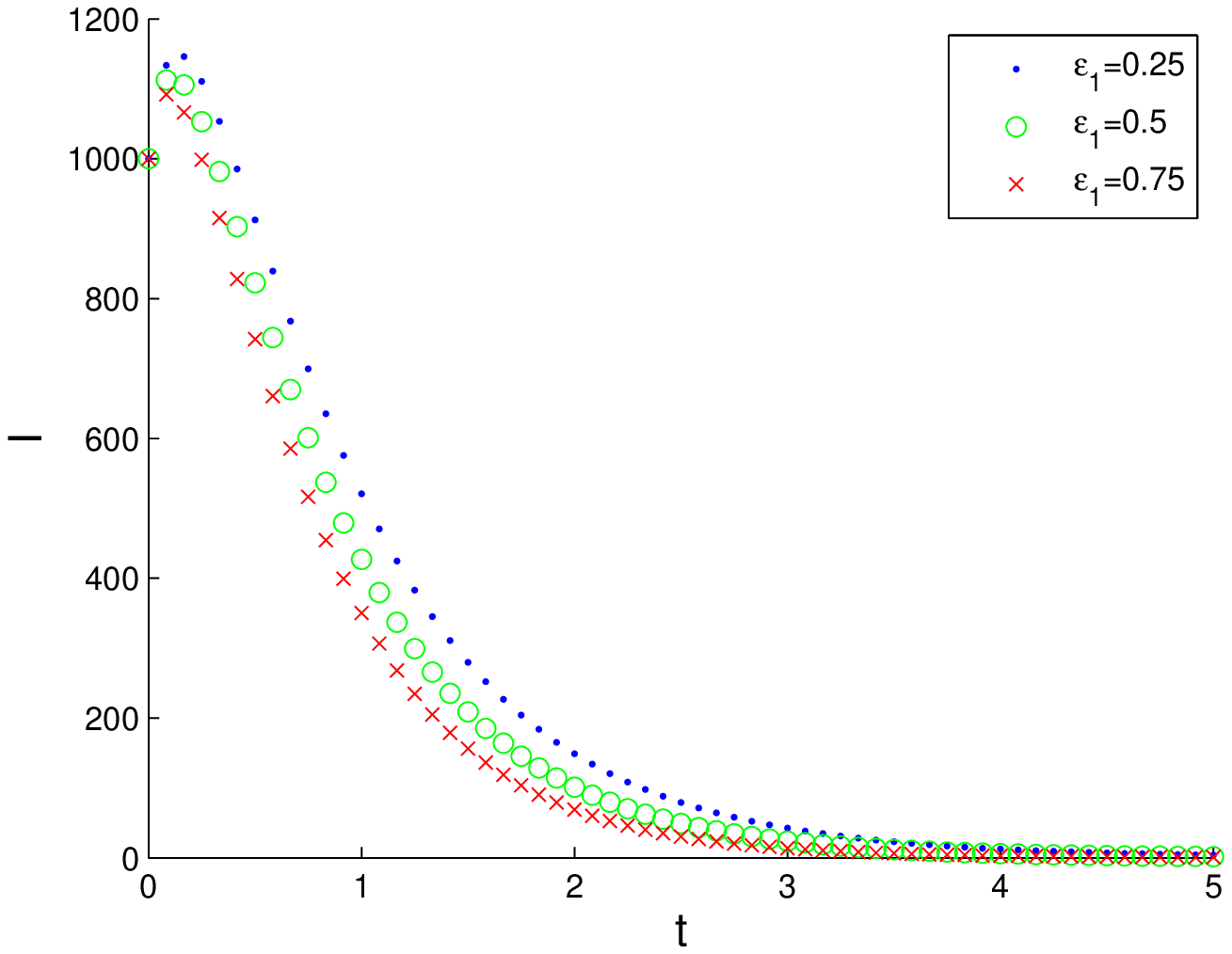,width=0.333\linewidth}\label{tb:sol:eps1s:l}}%
\subfigure[$L_2$ for $f_2=10$]{\epsfig{file=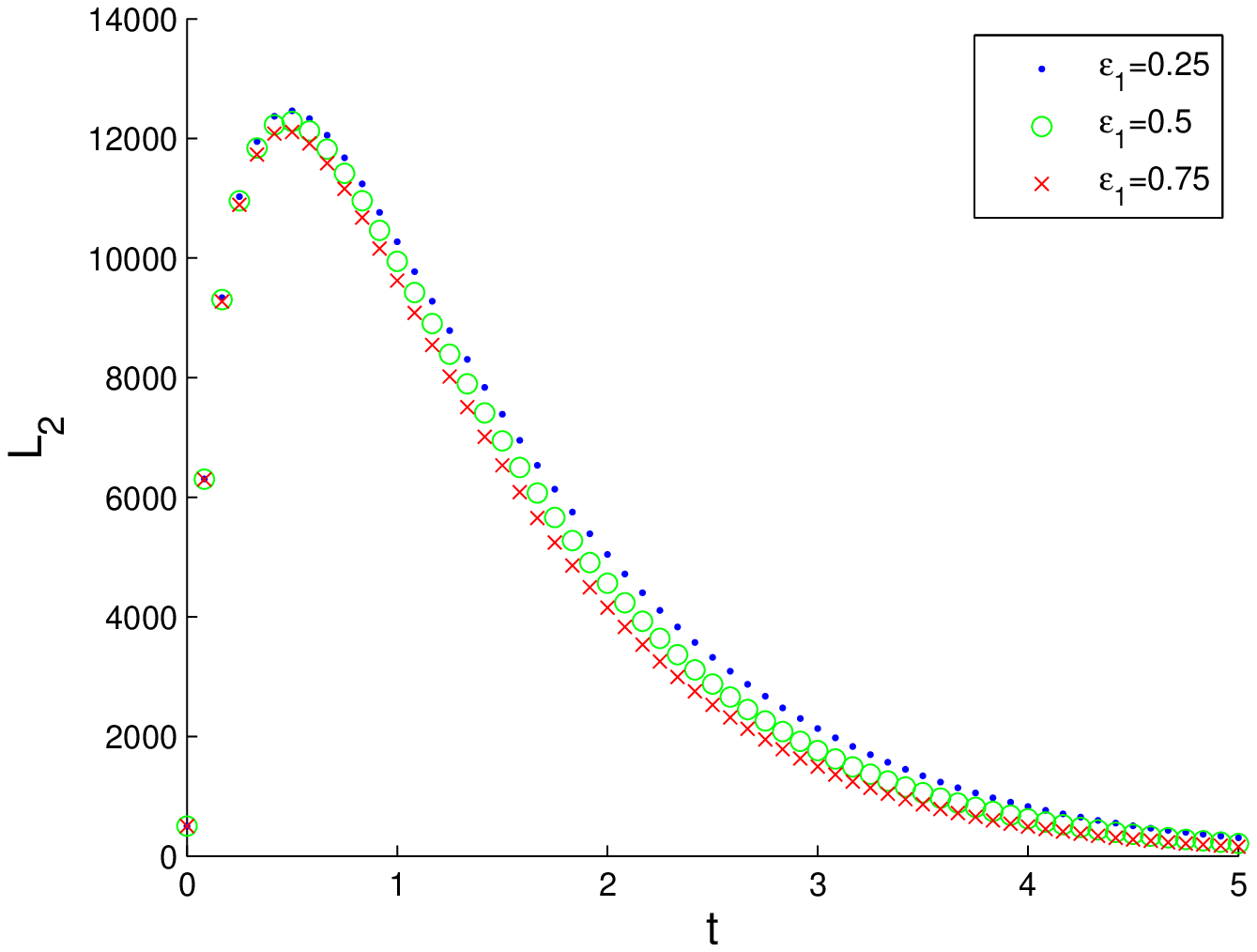,width=0.333\linewidth}\label{tb:sol:eps1s:m}}
\caption{Changes of controls, infectious and persistent latent individuals for different values of $\epsilon_1$.}
\label{tb:sol:eps1s}
\end{figure}
\begin{figure}
\centering
\includegraphics[width=0.67\textwidth]{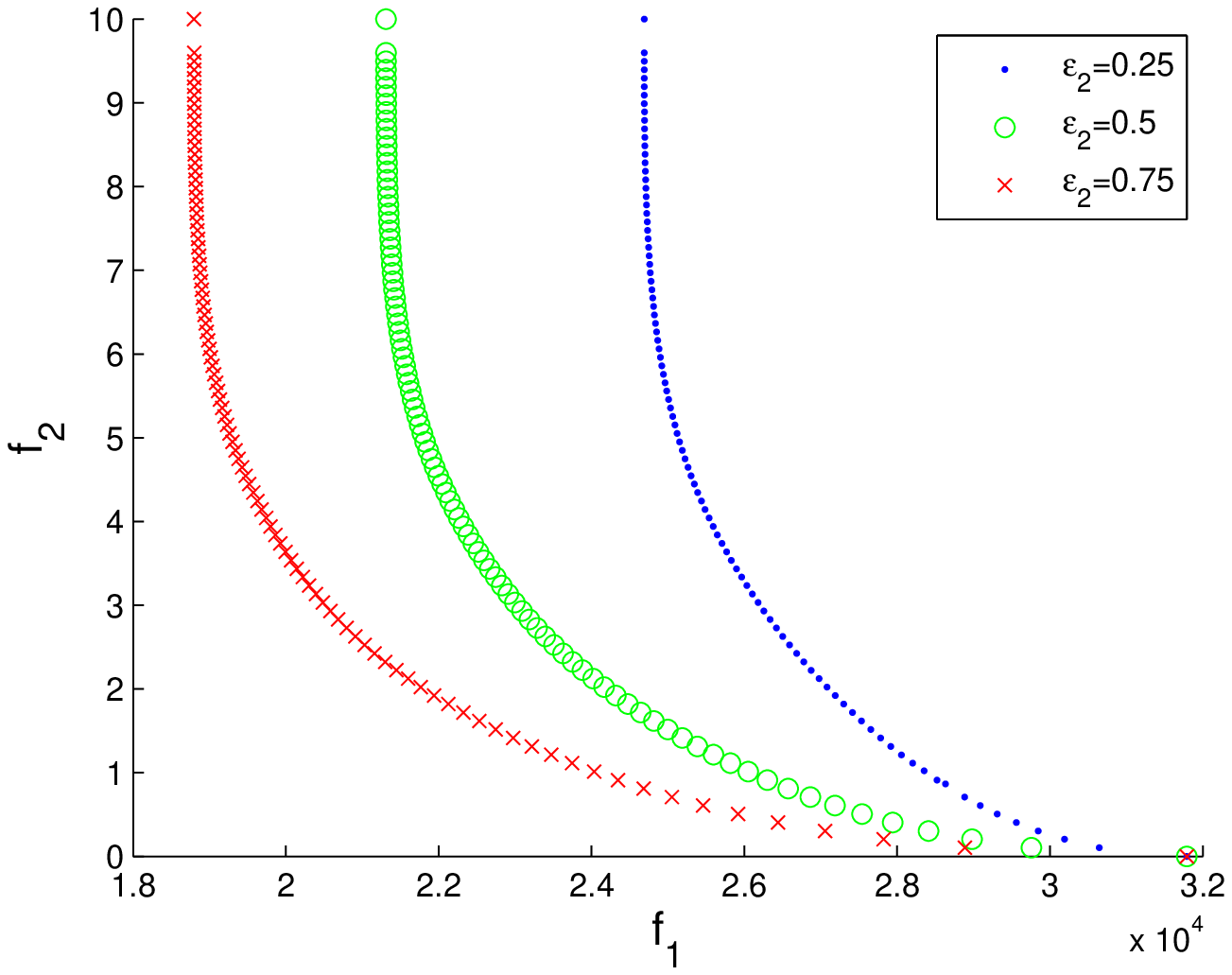}
\caption{Trade-off curves for different values of $\epsilon_2$.}
\label{tb:tradeoffs:eps2s}
\end{figure}
\begin{figure}
\centering
\subfigure[$I$ for $f_2=0$]{\epsfig{file=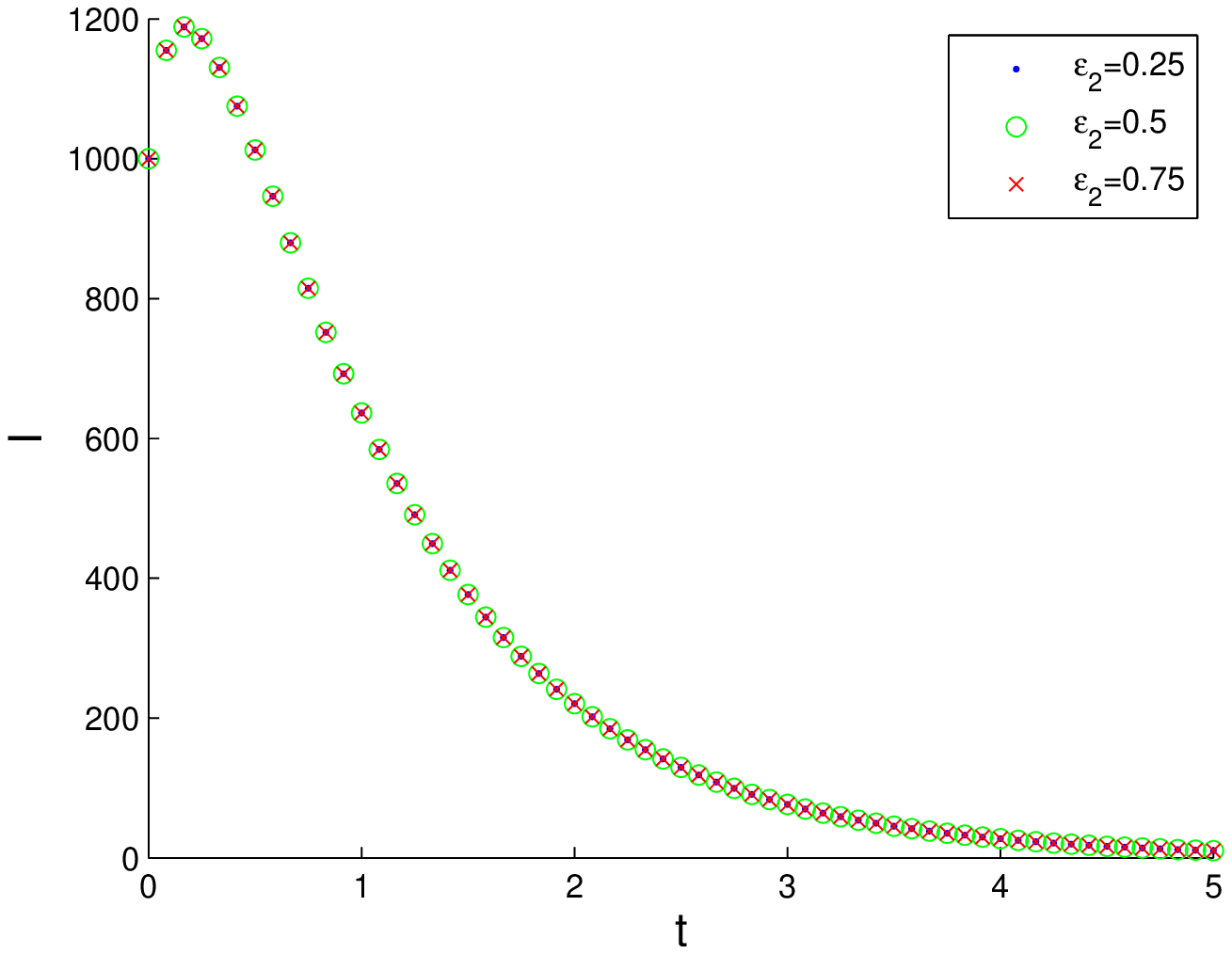,width=0.333\linewidth}\label{tb:sol:eps2s:a}}%
\subfigure[$L_2$ for $f_2=0$]{\epsfig{file=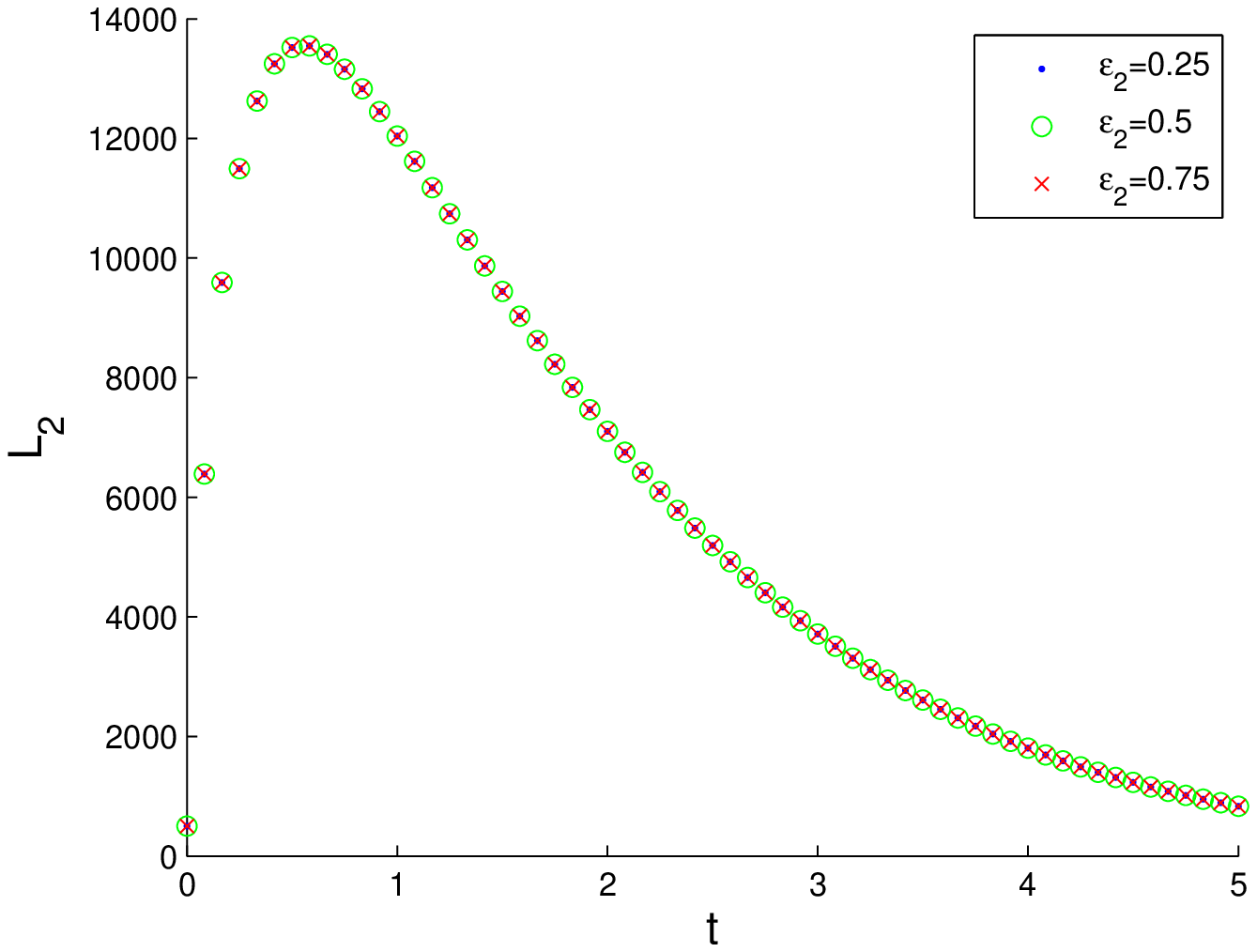,width=0.333\linewidth}\label{tb:sol:eps2s:b}}
\subfigure[$u_1$ for $f_2=2.5$]{\epsfig{file=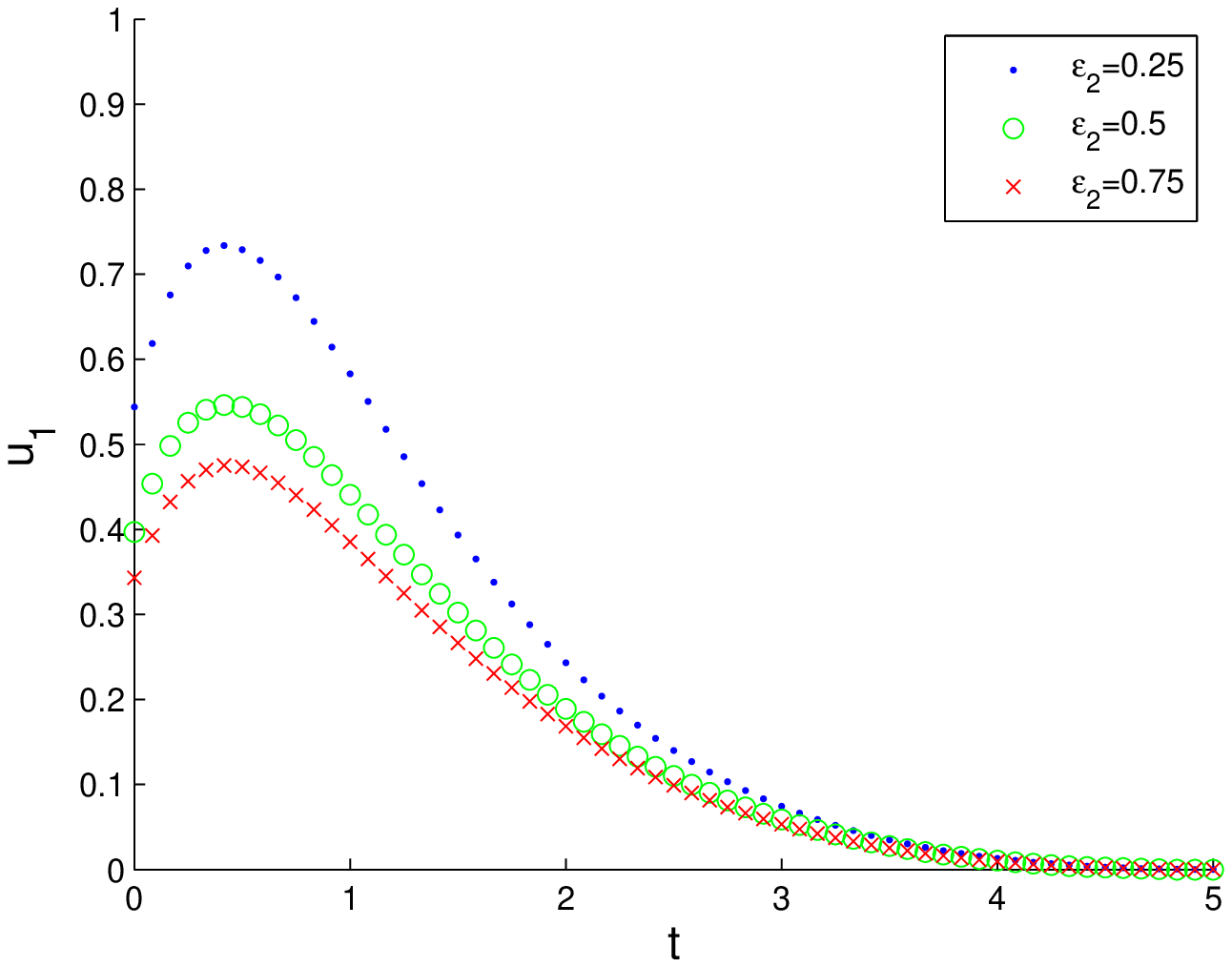,width=0.333\linewidth}\label{tb:sol:eps2s:c}}%
\subfigure[$u_2$ for $f_2=2.5$]{\epsfig{file=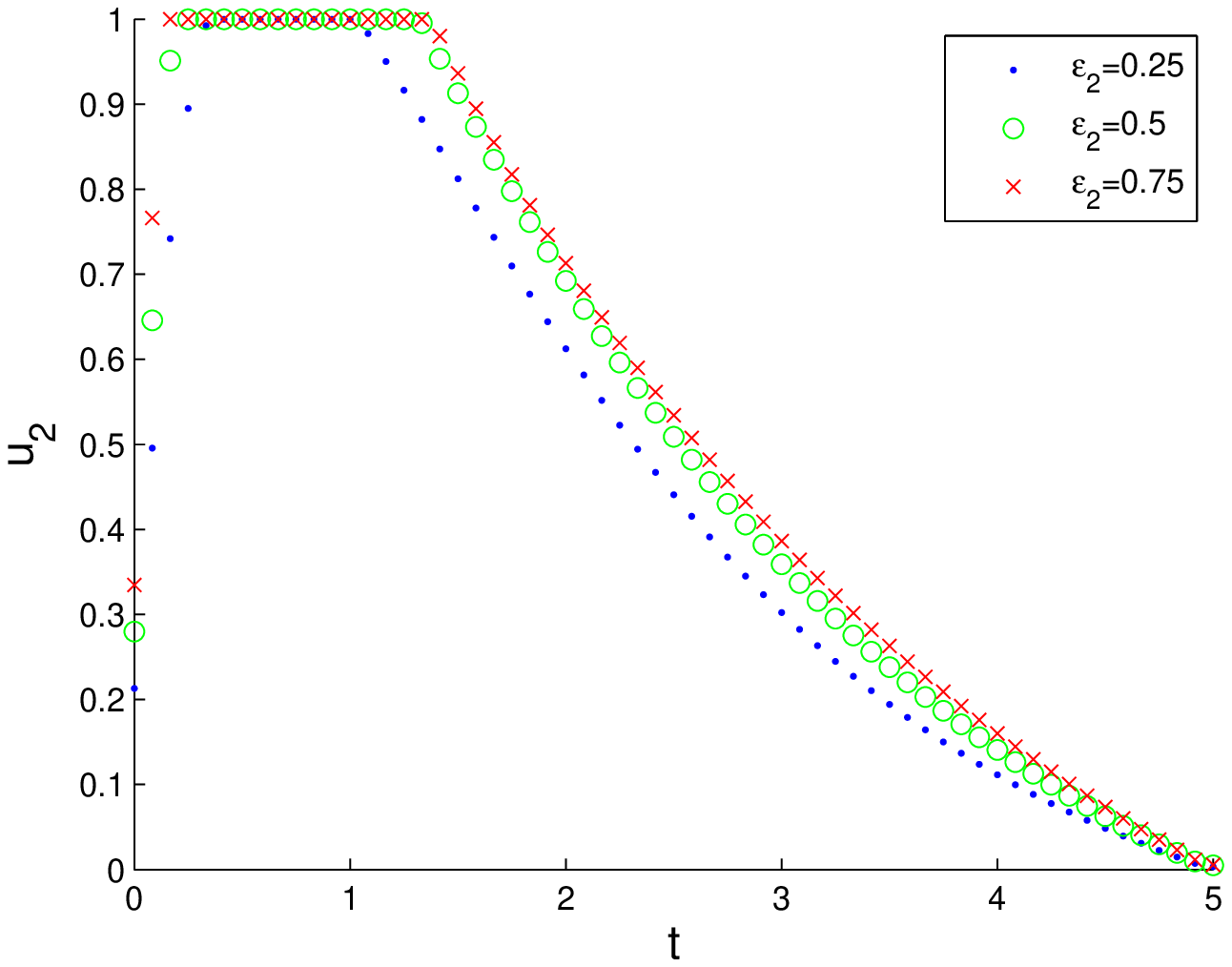,width=0.333\linewidth}\label{tb:sol:eps2s:d}}%
\subfigure[$I+L_2$ for $f_2=2.5$]{\epsfig{file=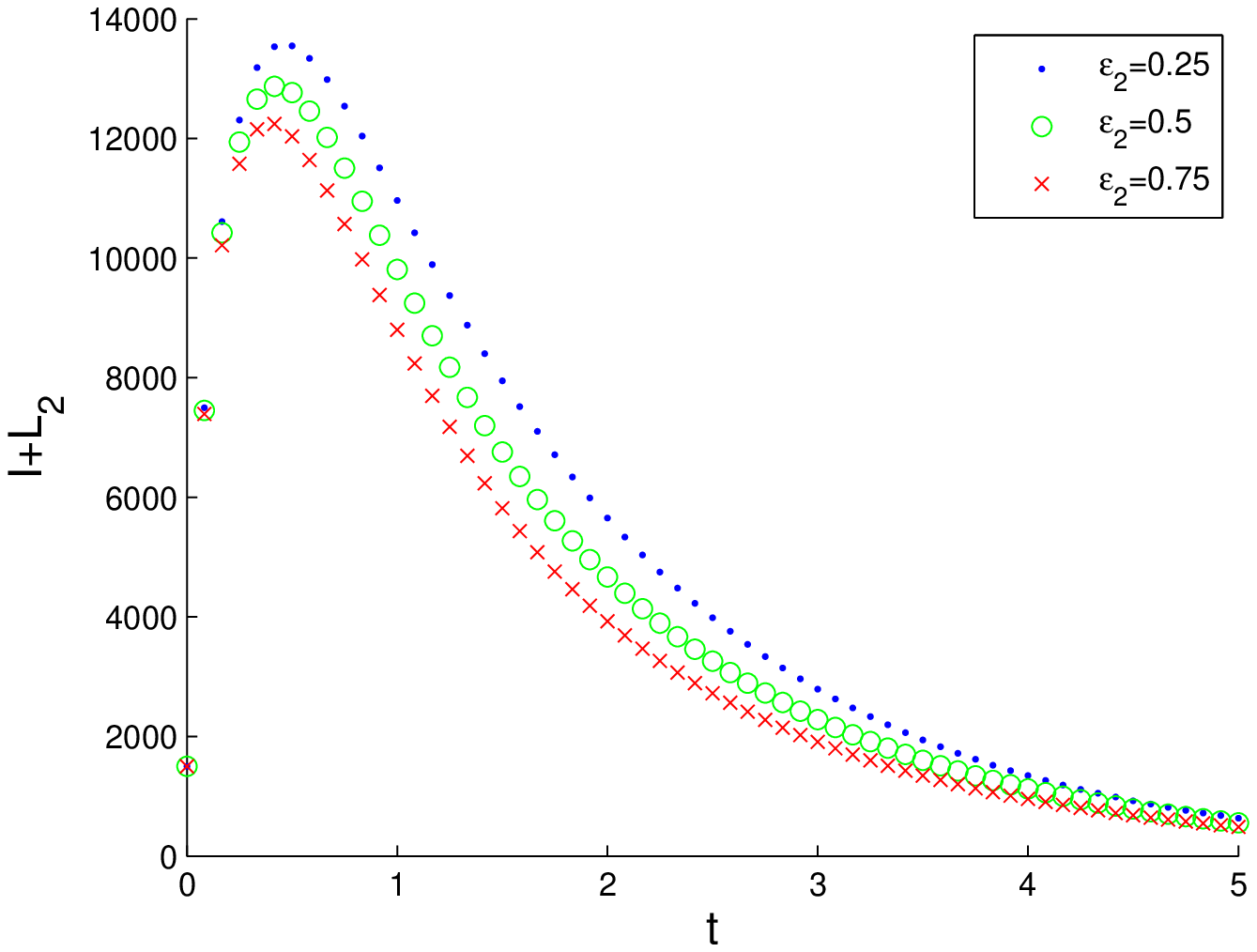,width=0.333\linewidth}\label{tb:sol:eps2s:e}}
\subfigure[$u_1$ for $f_2=5$]{\epsfig{file=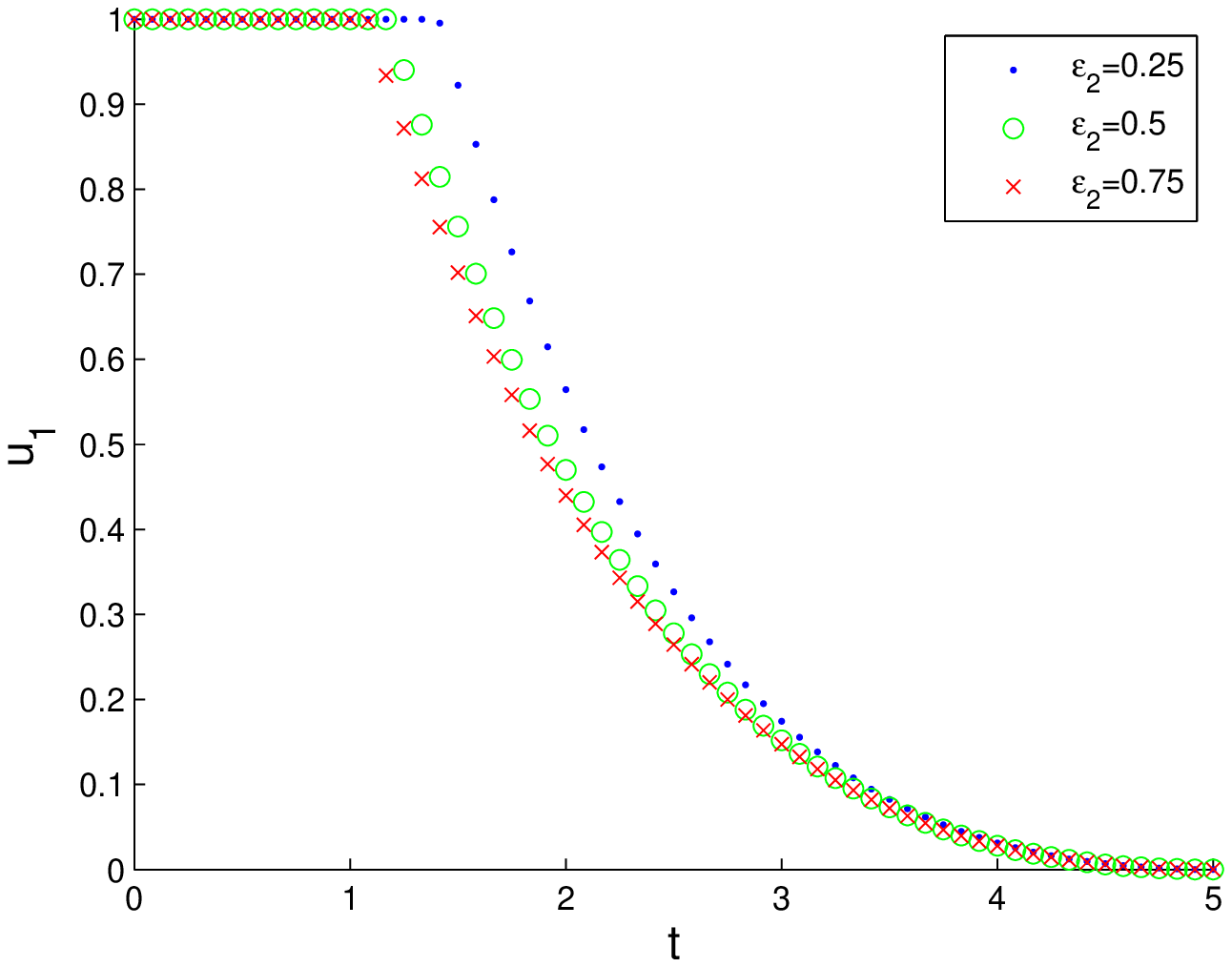,width=0.333\linewidth}\label{tb:sol:eps2s:f}}%
\subfigure[$u_2$ for $f_2=5$]{\epsfig{file=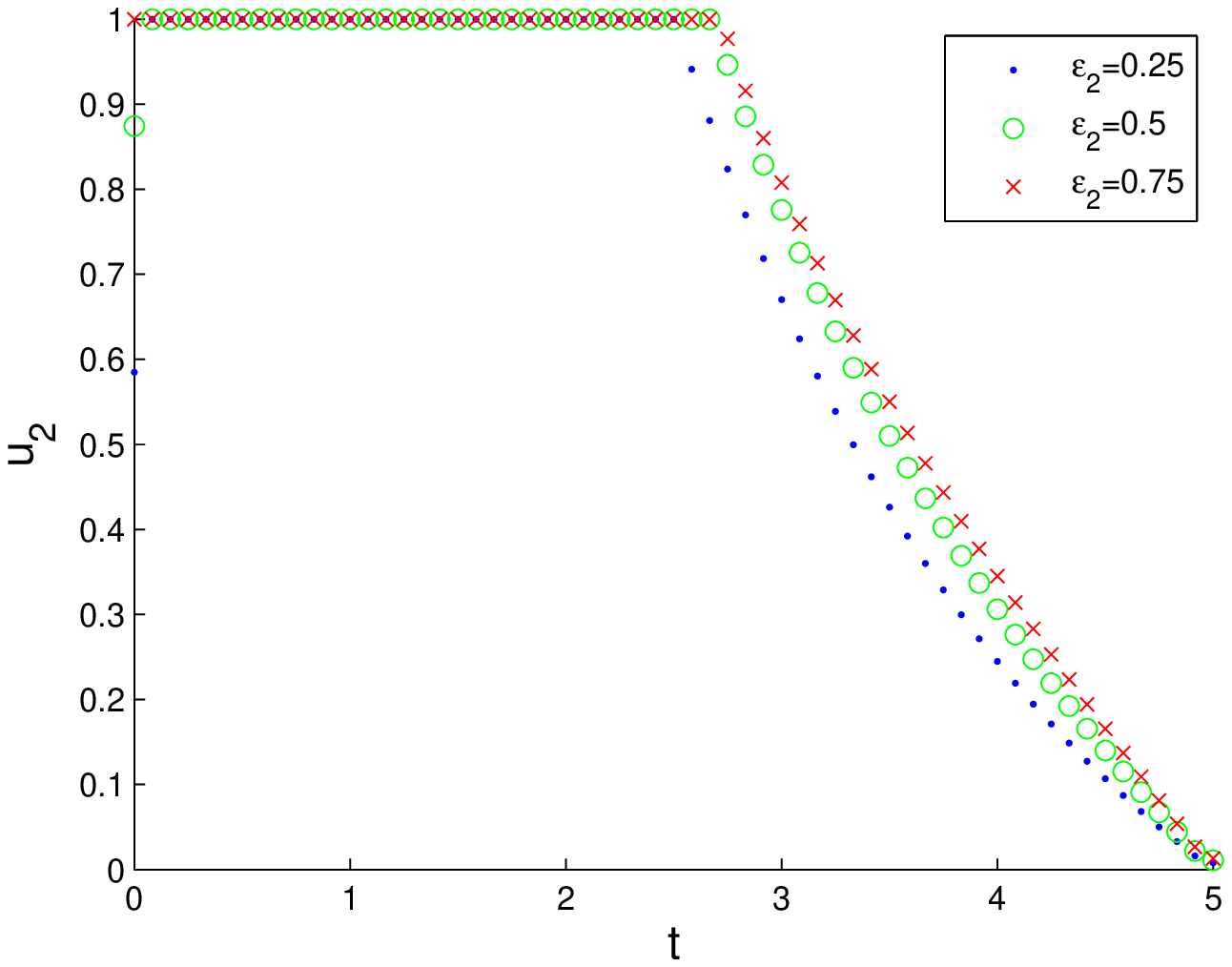,width=0.333\linewidth}\label{tb:sol:eps2s:g}}%
\subfigure[$I+L_2$ for $f_2=5$]{\epsfig{file=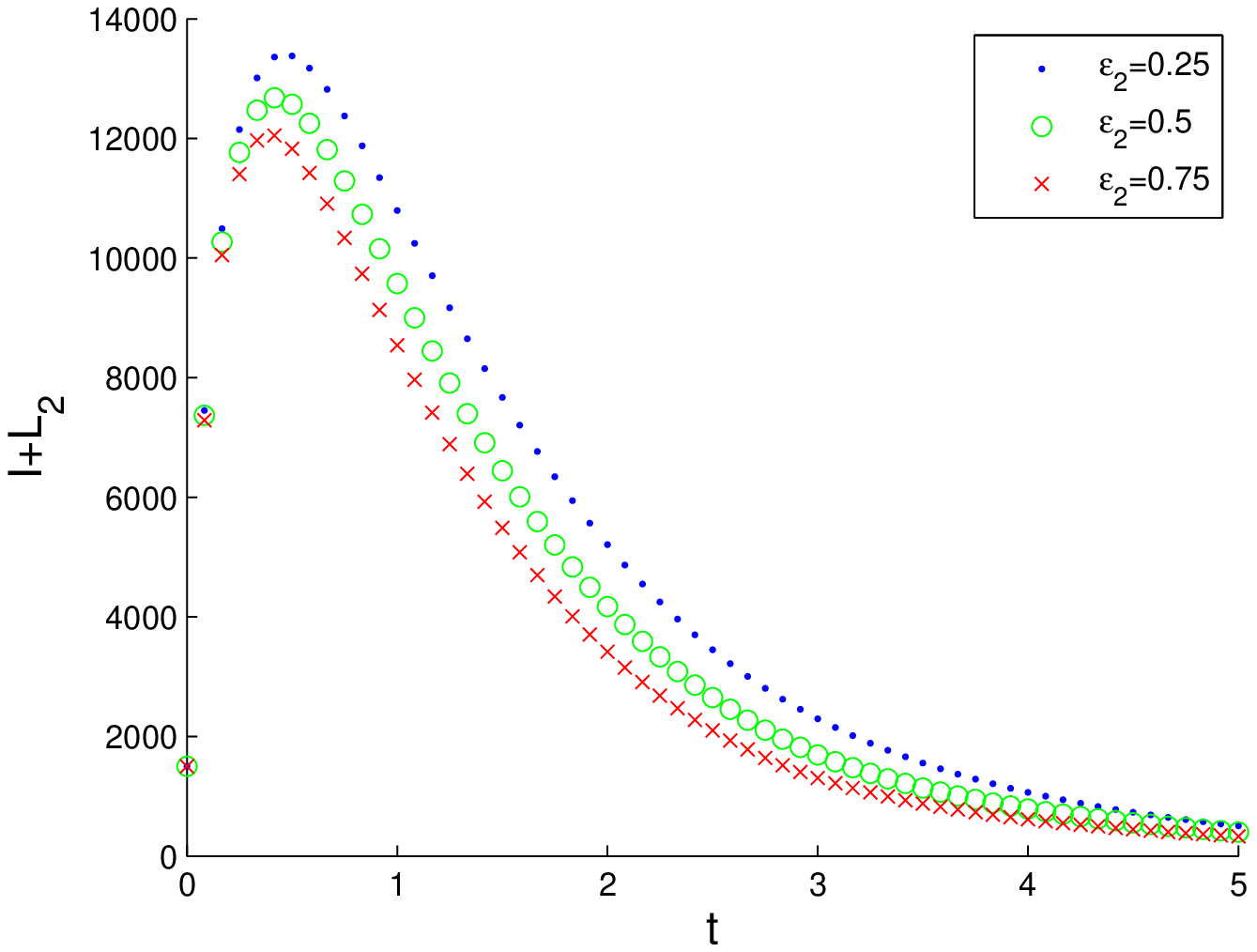,width=0.333\linewidth}\label{tb:sol:eps2s:h}}
\subfigure[$u_1$ for $f_2=7.5$]{\epsfig{file=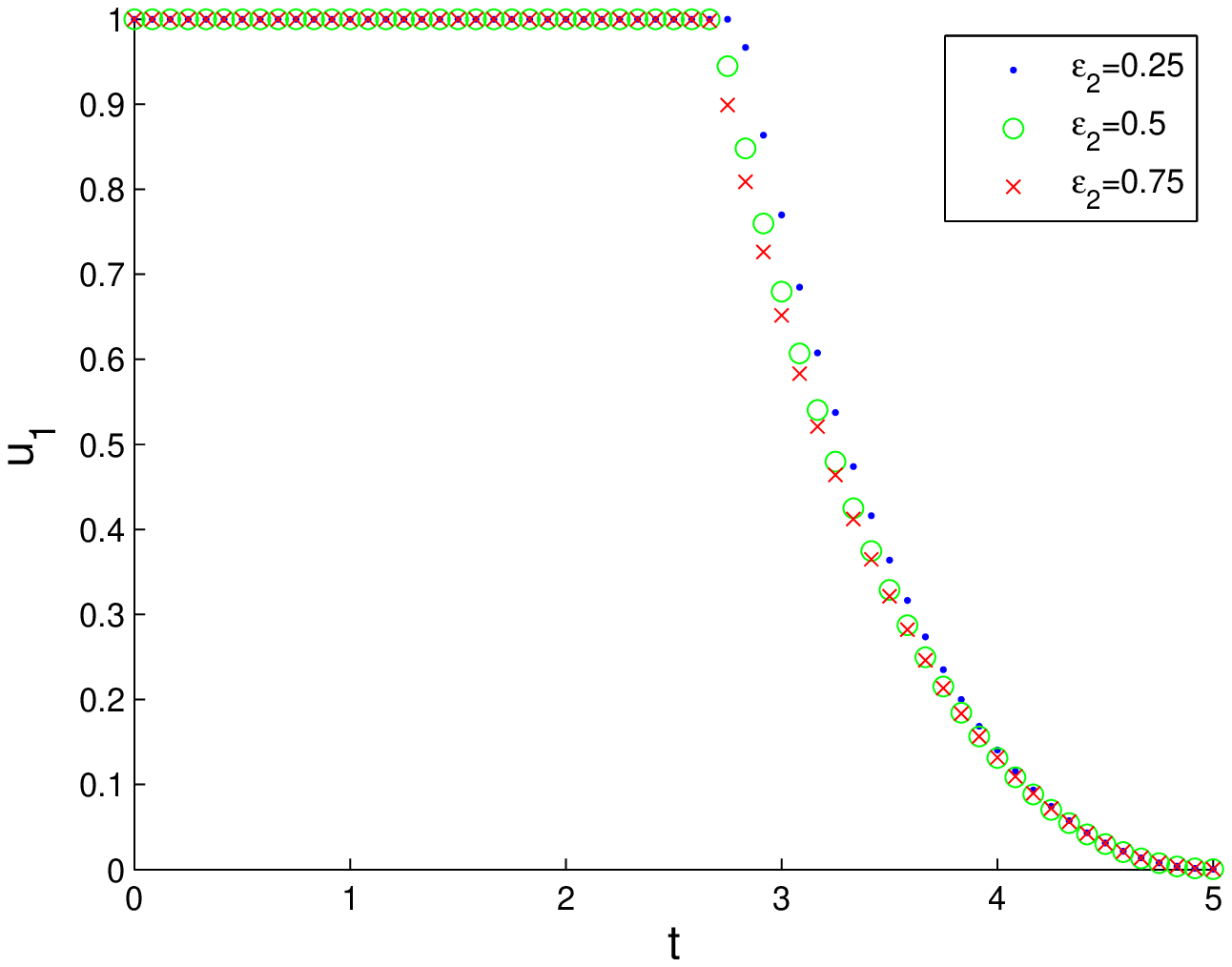,width=0.333\linewidth}\label{tb:sol:eps2s:i}}%
\subfigure[$u_2$ for $f_2=7.5$]{\epsfig{file=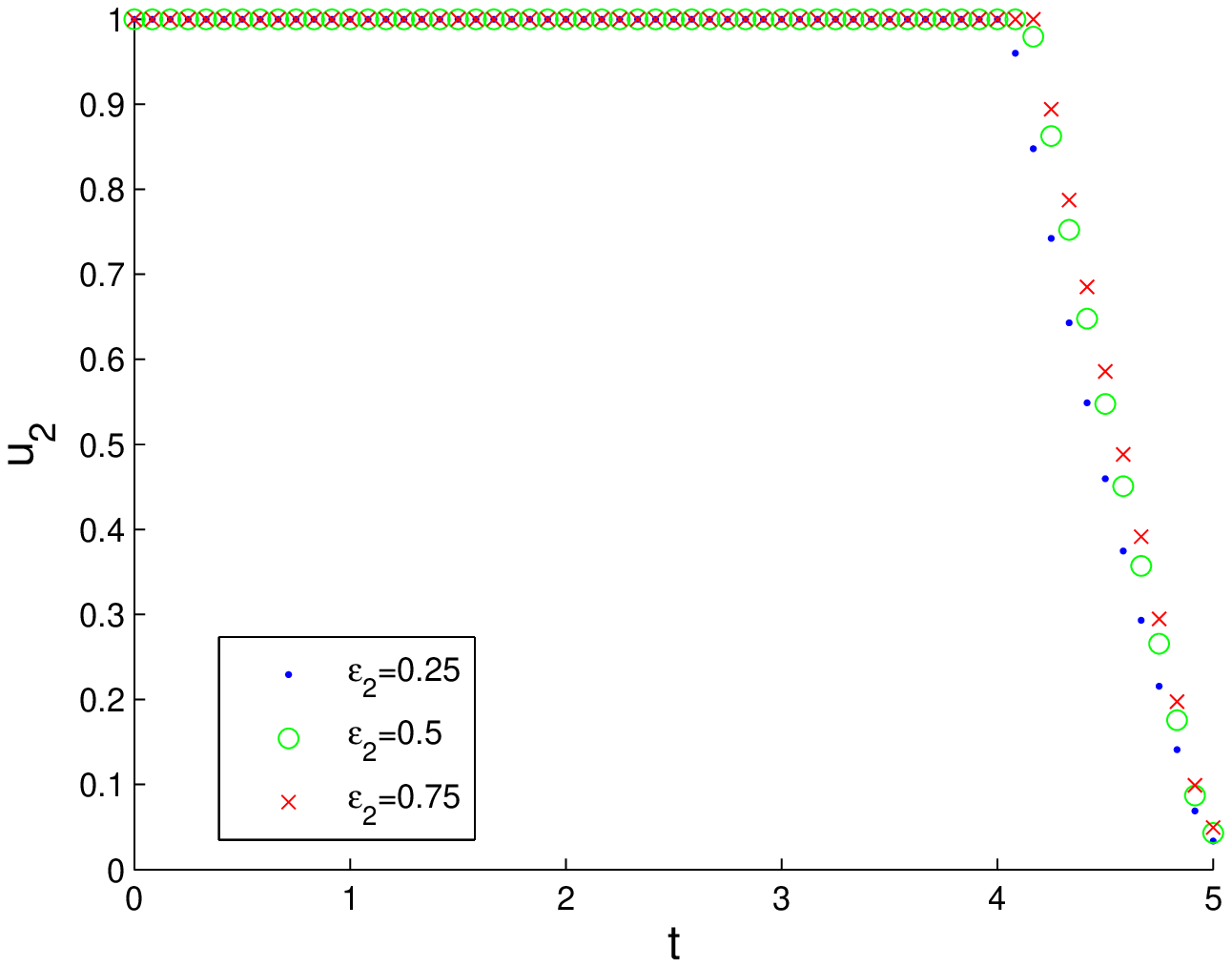,width=0.333\linewidth}\label{tb:sol:eps2s:j}}%
\subfigure[$I+L_2$ for $f_2=7.5$]{\epsfig{file=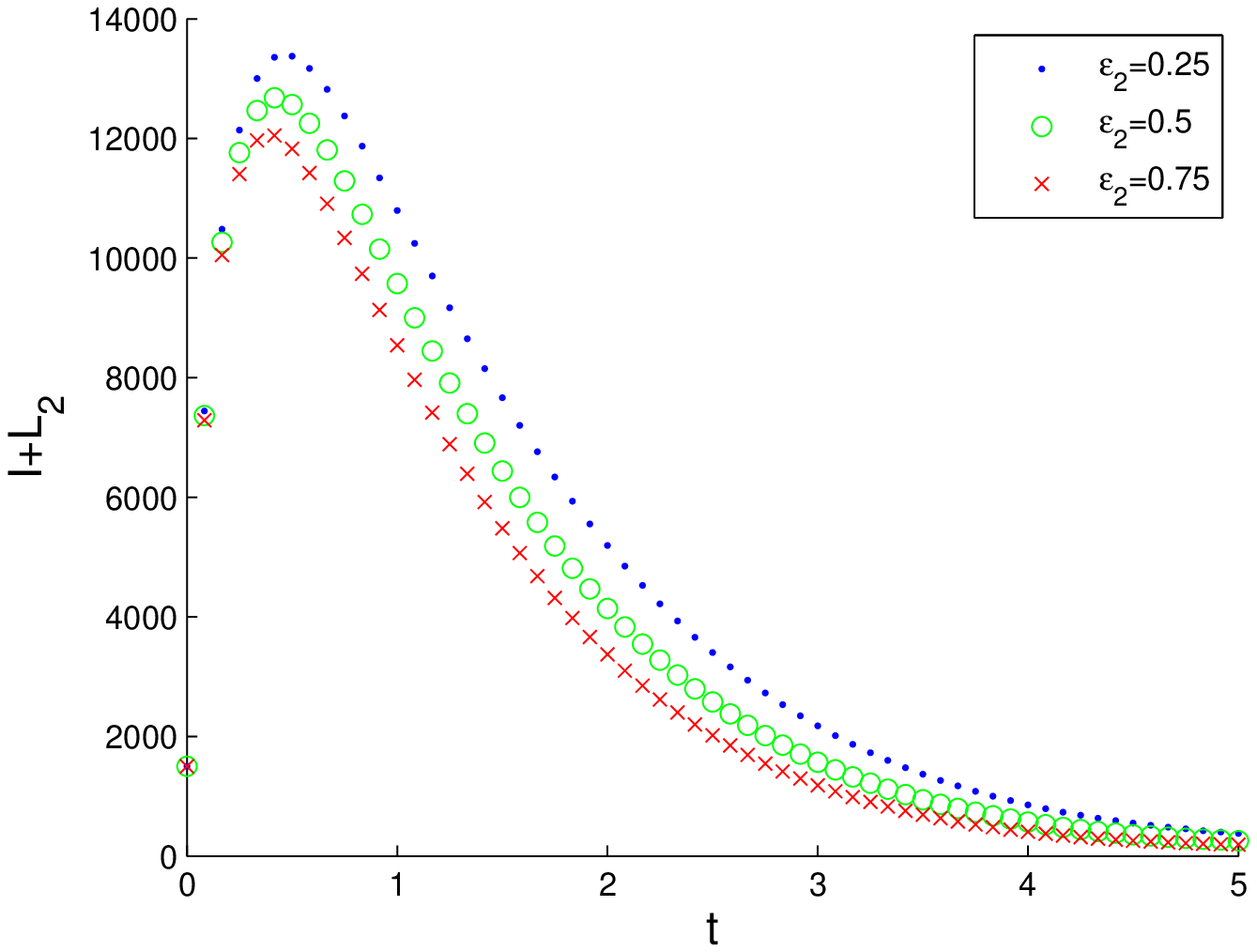,width=0.333\linewidth}\label{tb:sol:eps2s:k}}
\subfigure[$I$ for $f_2=10$]{\epsfig{file=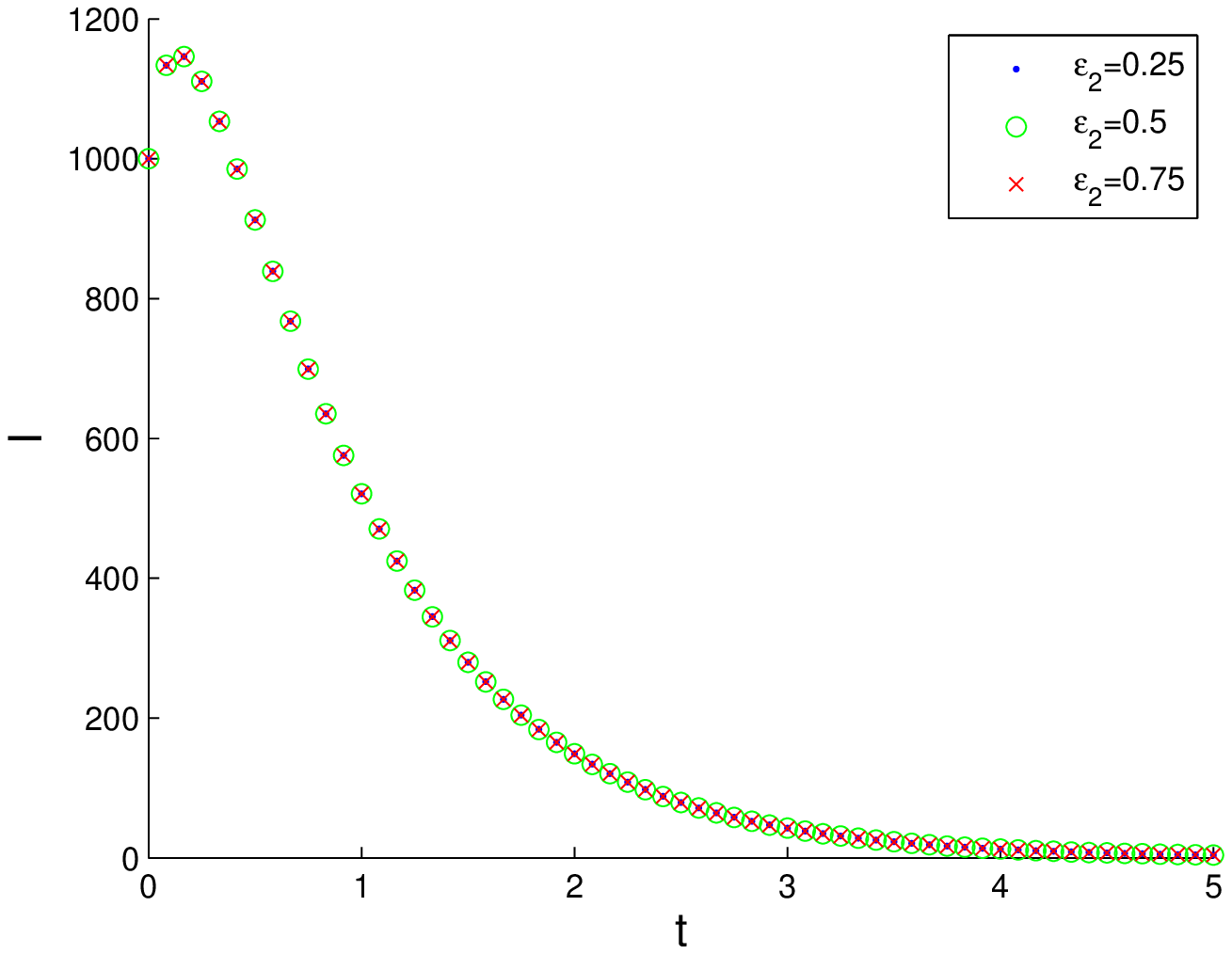,width=0.333\linewidth}\label{tb:sol:eps2s:l}}%
\subfigure[$L_2$ for $f_2=10$]{\epsfig{file=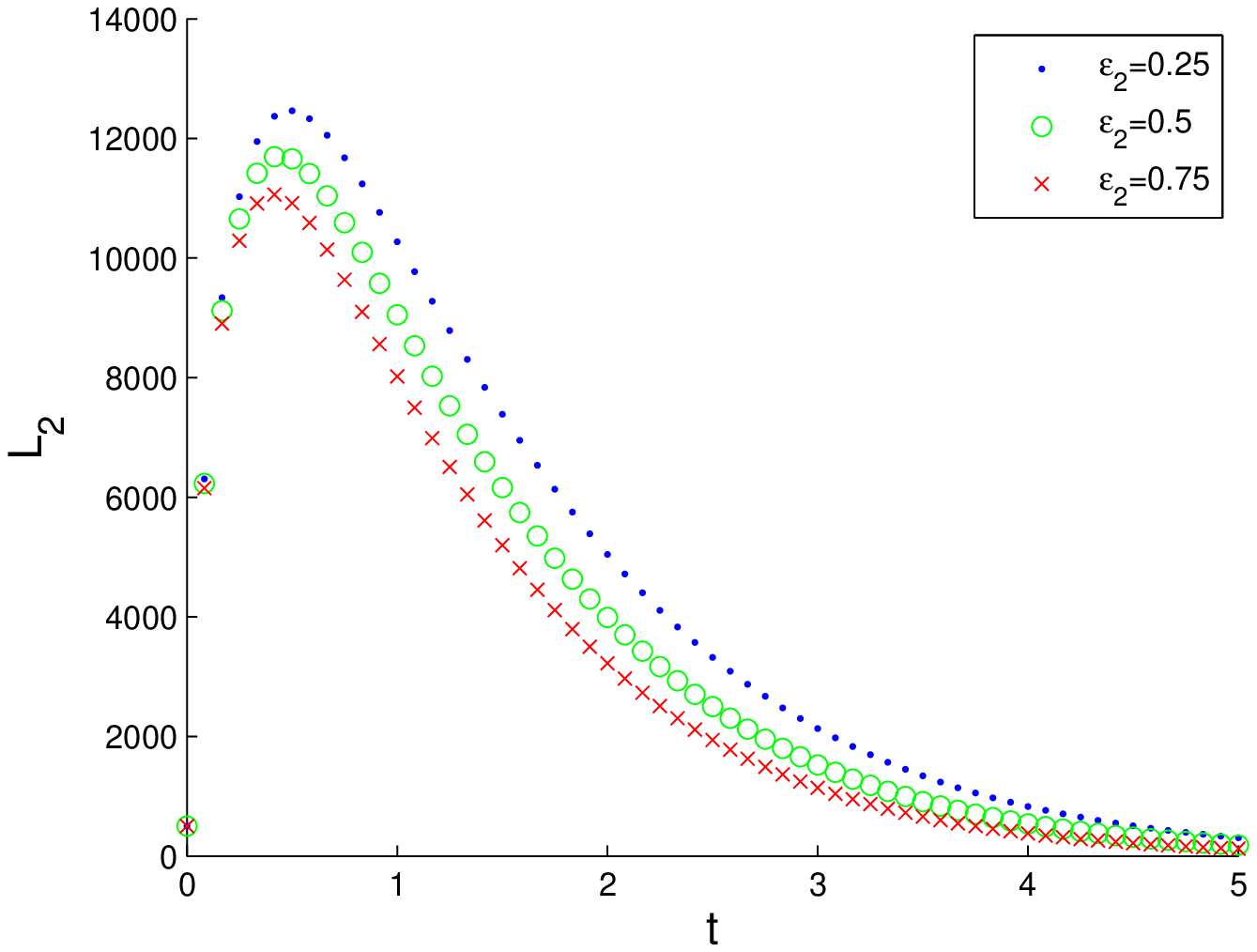,width=0.333\linewidth}\label{tb:sol:eps2s:m}}
\caption{Changes of controls, infectious and persistent latent individuals for different values of $\beta$.}
\label{tb:sol:eps2s}
\end{figure}

In the following, we discuss the obtained optimal solutions to the problem \eqref{mo:problem},
considering the variations of some model parameters separately. Since there is a set of optimal
solutions to the problem~\eqref{mo:problem}, considering all the solutions is somewhat cumbersome,
and we divide the obtained trade-off solutions into different parts, selecting a representative solution
to each part for analysis. For this purpose, we consider five cases: $f_2=0$, $f_2=2.5$, $f_2=5$,
$f_2=7.5$, and $f_2=10$. Each case represents an amount of available resources for treatment.
For each case, the best solution with respect to $f_1$ in the corresponding set of trade-off
solutions is considered. In particular, the first case ($f_2=0$) represents the situation
without controls ($u_1(\cdot),u_2(\cdot) \equiv 0$), reflecting the economical perspective.
The last case ($f_2=10$) represents the situation where the maximum controls are applied
($u_1(\cdot),u_2(\cdot) \equiv 1$), being the most preferable scenario from medical perspective.
Whereas three other cases ($f_2=2.5$, $f_2=5$, $f_2=7.5$) represent intermediate scenarios,
providing trade-offs between the number of affected by the tuberculosis and expenses for treatment.

Figure~\ref{tb:tradeoffs:betas} plots the trade-off solutions obtained for four different
values of the transition coefficient, $\beta$. From the figure, it can be seen that the higher
$\beta$ the higher number of infectious and persistent latent individuals is. Also,
the difference between the worst and the best case scenarios from medical perspective
becomes larger, if $\beta$ is increased.

Figure~\ref{tb:sol:betas} shows the changes of the number of infectious and persistent
latent individuals, as well as the controls, for solutions in different parts of the
trade-off curves with varying $\beta$. Figures~\ref{tb:sol:betas:a} and~\ref{tb:sol:betas:b}
show the changes of $I$ and $L_2$, respectively, without the controls ($f_2=0$).
On the other hand,  Figures~\ref{tb:sol:betas:l} and~\ref{tb:sol:betas:m} illustrate
the case where the maximum controls are applied ($u_1(\cdot),u_2(\cdot) \equiv 1$).
From these figures and Figures~\ref{tb:sol:betas:e}, \ref{tb:sol:betas:h} and \ref{tb:sol:betas:k},
it can be seen that if $\beta$ increases, then the number of infectious and persistent
latent individuals grows. From Figures~\ref{tb:sol:betas:c}, \ref{tb:sol:betas:f}
and \ref{tb:sol:betas:i}, one can see that $u_1$ is usually larger for higher values of $\beta$.
However, Figures~\ref{tb:sol:betas:d}, \ref{tb:sol:betas:g} and \ref{tb:sol:betas:j}
suggest that optimal values of $u_2$ are not always larger for higher values of $\beta$.
Thus, with an increasing of $\beta$, more effort must be put on the prevention of failure
of treatment in active infectious individuals. Moreover, the fraction of individuals
that are put under treatment must be decreased during the first part of the period,
and increased during the second part of the period when $\beta$ grows.

Figure~\ref{tb:tradeoffs:pops} displays the trade-off solutions obtained for three
different population sizes, $N$. We can see that larger population
sizes result in the higher numbers of infectious and persistent latent individuals.
Also, the range of optimal values of $f_1$ grows when $N$ rises.

Figure~\ref{tb:sol:pops} shows the changes of the controls as well as the fractions
of the number of infectious and persistent latent individuals for solutions
in different parts of the trade-off curves for different $N$. Although the number
of infectious and persistent latent individuals grows for larger $N$,
Figures~\ref{tb:sol:pops:a},~\ref{tb:sol:pops:b},~\ref{tb:sol:pops:l}
and~\ref{tb:sol:pops:m} show that the fractions of $I/N$ and $L_2/N$
do not vary when $(u_1(\cdot),u_2(\cdot)) \equiv 0$ and $(u_1(\cdot),u_2(\cdot)) \equiv 1$,
respectively. Similarly, the fractions $(I+L_2)/N$ remain unchanged for intermediate
solutions when $N$ is varied (Figures~\ref{tb:sol:pops:e},~\ref{tb:sol:pops:h}
and~\ref{tb:sol:pops:k}). From the plots for the controls in Figure~\ref{tb:sol:pops},
we can seen that the optimal control strategies do not depend on the population size.

Figure~\ref{tb:tradeoffs:eps1s} presents the trade-off solutions when the efficacy
of treatment of active tuberculosis individuals, $\epsilon_1$, is varying. The figure
shows that increasing $\epsilon_1$ allows to reduce the number of infectious
and persistent latent individuals.

The plots for infectious and persistent latent individuals presented in Figure~\ref{tb:sol:eps1s}
demonstrate that the higher efficacy of treatment of active individuals the lower the number
of infectious and persistent latent individuals is. Concerning the controls,
Figures~\ref{tb:sol:eps1s:c},~\ref{tb:sol:eps1s:f} and~\ref{tb:sol:eps1s:i}
show that the optimal control $u_1$ is higher when $\epsilon_1$ is increased.
On the contrary, the optimal control $u_2$ decreases when $\epsilon_1$ is increased
(Figures~\ref{tb:sol:eps1s:d},~\ref{tb:sol:eps1s:g} and~\ref{tb:sol:eps1s:j}).
This suggests that for efficiently reducing $I+L_2$ when increasing $\epsilon_1$
we must focus on policies associated with $u_1$, namely the supervision and
the support of active infectious individuals, and decrease the fraction
of individuals that are put under treatment.

Figure~\ref{tb:tradeoffs:eps2s} shows the trade-off solutions when the efficacy
of treatment of latent tuberculosis individuals, $\epsilon_2$, is varying.
Similarly to the case of varying $\epsilon_1$ shown in Figure~\ref{tb:tradeoffs:eps1s},
this figure suggests that increasing $\epsilon_2$ allows to reduce the number of infectious
and persistent latent individuals. Though, the larger decrease in $I+L_2$ can be achieved
by increasing $\epsilon_2$ compared to the same values of $\epsilon_1$.

Figures~\ref{tb:sol:eps2s:l} and~\ref{tb:sol:eps2s:m} show that by implementing the maximum controls,
$I$ has the same values for different $\epsilon_2$, whereas $L_2$ decreases when $\epsilon_2$
is increased. Also, $I+L_2$ becomes smaller when $\epsilon_2$ increases, as it can be seen
analyzing intermediate trade-off solutions in Figures~\ref{tb:sol:eps2s:e},~\ref{tb:sol:eps2s:h}
and~\ref{tb:sol:eps2s:k}. The plots for the controls presented in Figure~\ref{tb:sol:eps2s}
show that $u_1$ becomes smaller when $\epsilon_2$ increases, whereas $u_2$ grows when
$\epsilon_2$ is increased. This suggest that, to efficiently reduce $I+L_2$ with increasing
$\epsilon_2$, we must focus on implementation of $u_2$, which is related to the fraction
of persistent latent individuals that is put under treatment.


\section{Methods Comparison}
\label{sec:5}

In this section, we compare our approach to find the optimal control
strategies in a tuberculosis model, based on the $\epsilon$-constraint method
(Section~\ref{sec:results}), with two other scalarization methods
for solving multiobjective optimization problems.
The first one is based on the goal attainment method, which was applied for solving
multiobjective control problems in some recent studies \cite{KayaMaurer2014,LogistHouska2010}.
Another popular method is based on Chebyshev problem (or Chebyshev method).
Follows the description of the two methods.


\subsection{Goal attainment method}

The goal attainment method \cite{MiettinenBook}
reformulates the problem shown in \eqref{mop} as folows:
\begin{equation}
\label{method:gam}
\left\{
\begin{array}{rl}
\underset{\boldsymbol{x} \in \Omega, \alpha \geq 0}{\text{minimize:}}   & \alpha \\
\text{subject to:} & w_1(f_1(\boldsymbol{x}) - z^{\ast}_1) \leq \alpha \\
& \qquad \vdots \\
& w_m(f_m(\boldsymbol{x}) - z^{\ast}_m) \leq \alpha,\\
\end{array}
\right.
\end{equation}
where $\boldsymbol{z}^{\ast} = (z^{\ast}_1,\ldots,z^{\ast}_m)^{\text{T}}$
is a reference point, $\boldsymbol{w} = (w_1, \ldots, w_m )^{\text{T}}$
is a weight vector ($\sum\limits_{j=1}\limits^{m} w_j = 1$),
and $\boldsymbol{x} \in \Omega, \alpha \in \mathbb{R}_{+}$ are variables.
Solving the above problem for different weight vectors allows to obtain multiple
Pareto optimal solutions. In \cite{EichfelderBook} the problem \eqref{method:gam}
is referred to as Pascoletti--Serafini scalarization.


\subsection{Chebyshev method}

The Chebyshev method belongs to the class of weighted metric methods
\cite{MiettinenBook}, which minimize the distance between some reference point
and the feasible objective region. The weighted metric methods use the weighted $L_p$
metrics for measuring the distance from any solution to the reference point.
Chebyshev method is referred to the case with $p=\infty$, which can be formulated as
\begin{equation}
\label{method:chb}
\underset{\boldsymbol{x} \in \Omega}{\text{minimize:}}
\max \limits_{1\leq i\leq m}
\, \{ w_i (f_i(\boldsymbol{x})-z^{\ast}_i) \},
\end{equation}
where $\boldsymbol{z}^{\ast} = (z^{\ast}_1,\ldots,z^{\ast}_m)^{\text{T}}$
is a reference point, and $\boldsymbol{w} = (w_1, \ldots, w_m )^{\text{T}}$
is a weight vector ($\sum\limits_{j=1}\limits^{m} w_j = 1$). Similarly
to the goal attainment method, for finding multiple Pareto optimal solutions
the problem \eqref{method:chb} must be solved for different weight vectors.
The problem shown in \eqref{method:chb} was originally introduced
in Bowman \cite{Bowman1976}. In the literature, the name of the method
may vary due to different ways of spelling.


\subsection{Results}

\begin{figure}
\centering
\subfigure[$\epsilon$-Constraint method]{\epsfig{file=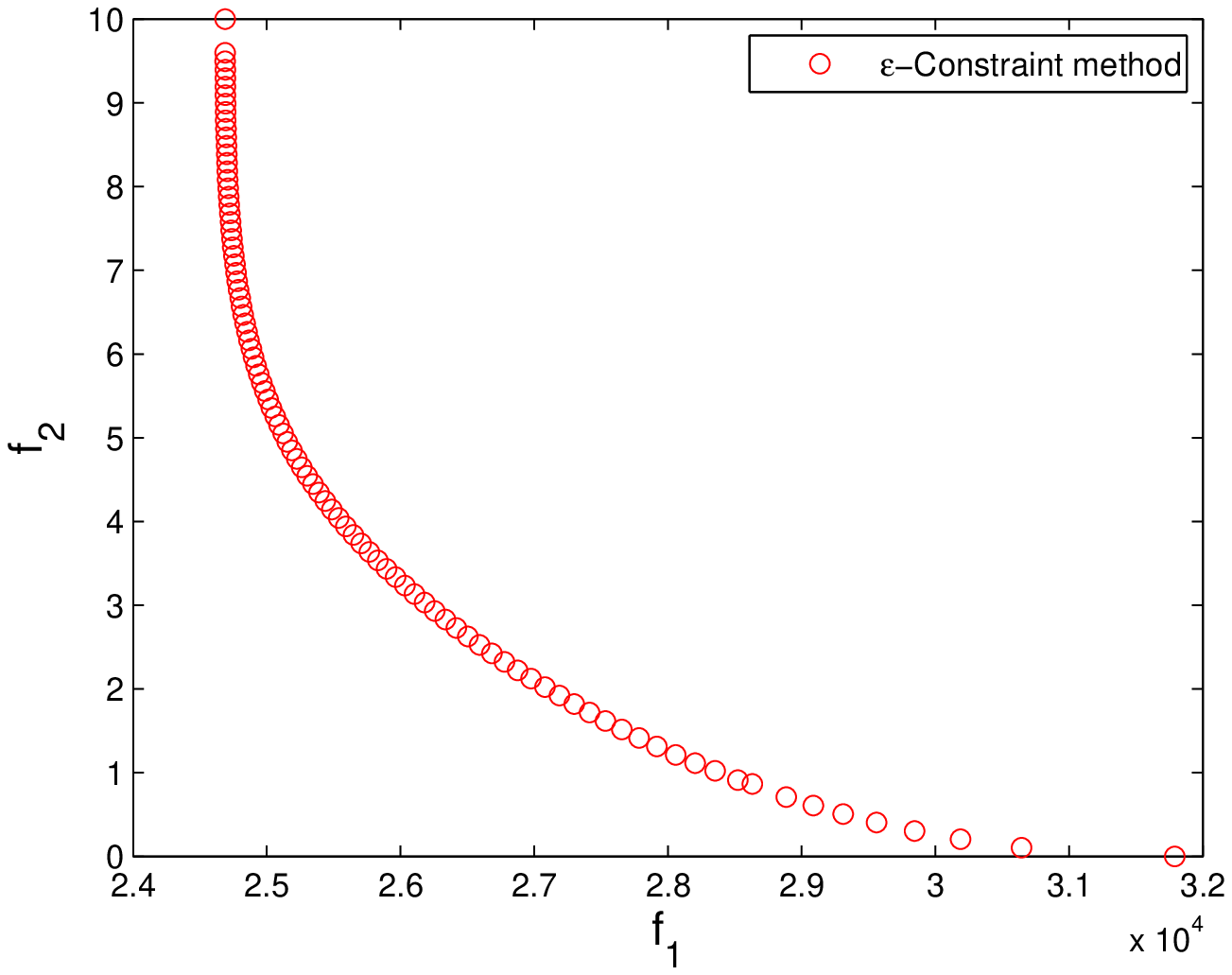,width=0.333\linewidth}\label{methods:eps}}%
\subfigure[Goal attainment method]{\epsfig{file=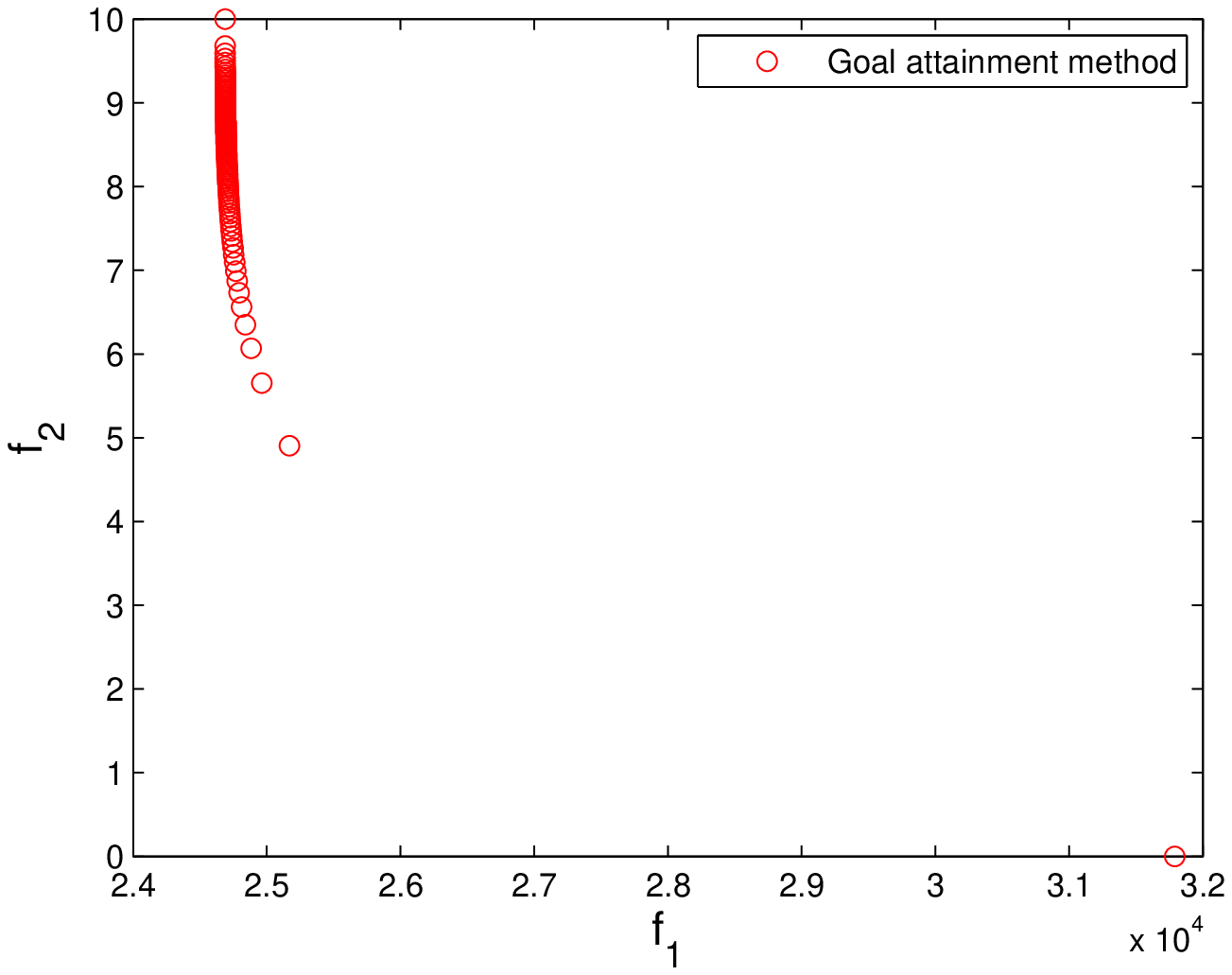,width=0.333\linewidth}\label{methods:gam}}%
\subfigure[Chebyshev method]{\epsfig{file=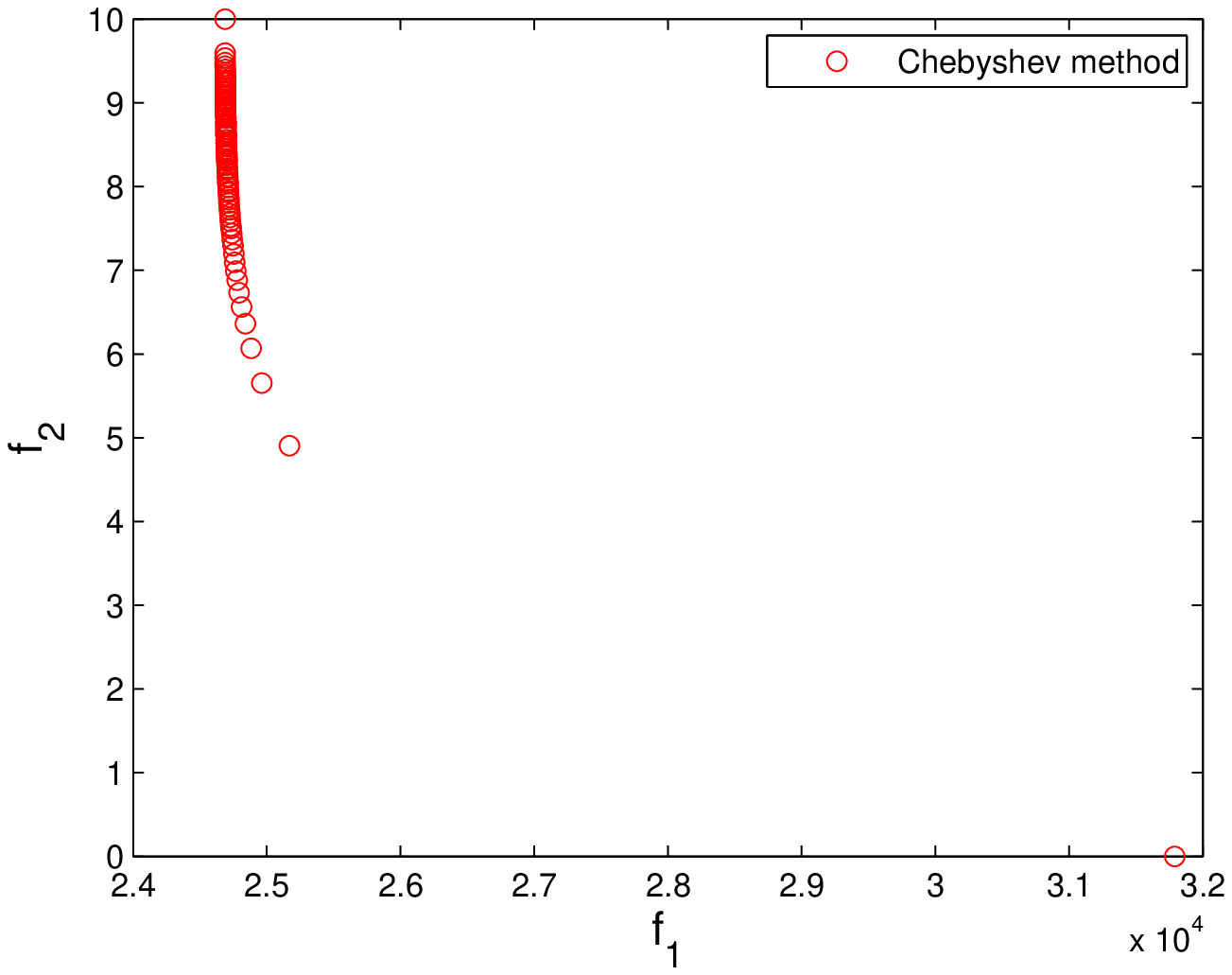,width=0.333\linewidth}\label{methods:chb}}%
\caption{Trade-off curves obtained by different methods.}\label{methods}
\end{figure}

For a set of 100 uniformly distributed weight vectors,
each method was executed on the tuberculosis model, having the remaining
settings as in Section~\ref{sec:setup}. Figure~\ref{methods} plots trade-off
solutions obtained afterwards, including the $\epsilon$-constraint method.
As it is seen, the $\epsilon$-constraint method produces a well-distributed
set of trade-of solutions, whereas the solutions obtained by the other methods
do not cover the whole Pareto optimal region. This is because the objective functions
are differently scaled. The use of scalarization schemes based on the weight vectors
makes the goal attainment method and Chebyshev method sensitive to the scale of objectives,
resulting in significant performance deterioration, in terms of the uniformness
of obtained solutions, in the case of disparately scaled objectives.
On the other hand, the $\epsilon$-constraint method divides the Pareto region
into a number of subregions with respect to the values of $\epsilon$,
minimizing the chosen objective. In this case, the relative scale of objectives
does not matter significatively, with the approach being able
to locate solutions in all parts of the Pareto optimal region.

To quantitatively compare the outcomes of the three methods,
we rely on the hypervolume \cite{ZitzlerThiele1998}, which has been utilized
extensively for the comparison of multiobjective algorithms. The hypervolume
uses the volume of the dominated portion of the objective space as a measure
for the quality. Table~\ref{tab:hv} shows the hypervolume values for the three methods,
which are computed after normalizing the objective values of the obtained solutions
and using the nadir point as a reference point. The results summarized in Table~\ref{tab:hv}
confirm our previous observation, with the $\epsilon$-constraint method producing the best results.
\begin{table}[t]
\centering
\small
\begin{tabular}{cccc} \hline
& $\epsilon$-Constraint method & Goal attainment method & Chebyshev method \\ \hline
Hypervolume & {\bf 0.82481} & 0.50034 & 0.50031 \\ 	\hline		
\end{tabular}
\caption{Hypervolume values for trade-off solutions obtained by the three methods (the higher the better).}
\label{tab:hv}
\end{table}


\section{Conclusions}
\label{sec:conclusions}

The incidence rates of tuberculosis have been declining since 2004 worldwide.
Mortality rates, at global level, fell down around 45\% between 1990 and 2012,
and if the current rate of decline is sustained, by 2015 the target of a 50\%
reduction can be achieved. The reduction of mortality and incidence rates
is due to prevention and treatment policies that have been applied in the last years.

In this paper we study a mathematical model for tuberculosis from the optimal control
point of view, using a multiobjective approach. The optimal control strategies are found
by simultaneously minimizing the number of individuals affected by the tuberculosis
and the cost of implementation of prevention and treatment policies. This approach avoids
the use of additional weight coefficients to formulate a single cost functional and reflects
the intrinsic nature of the problem.

The results obtained in this study clearly show
that a multiobjective approach is effective to finding optimal control strategies
in a mathematical model for tuberculosis. The obtained trade-off solutions reveal
different perspectives on the implementation of prevention and treatment policies.
Once a set of optimal solutions is calculated, the final decision on the control strategy
can be made taking into account the goals of public health care and the available resources for treatment.
We also investigate the optimal control strategies with varying model parameters.
It is observed that as the transmission coefficient  increases, the fraction of active
infectious and persistent latent individuals increases as well, corresponding to the
case where the disease may become endemic. Varying the population size, the optimal
control strategies remain unchanged. Increasing the efficacy of control policies allows
to reduce the number of active infectious and persistent latent individuals. When
the measure of efficacy for some control is increased, the main focus of efficiently
dealing with the disease must be on the policies associated with the corresponding
control for which the efficacy is improved.

Finally, we compared the approach proposed in our work with other
scalarization techniques for multiobjective optimization. The obtained results
show that the $\epsilon$-constraint method is an appropriate choice for finding
the optimal control strategies in the tuberculosis model.

As future work, it would be interesting to further investigate different values
for the model parameters and observe the variations on the optimal control strategies.
We also plan to consider the second objective as an $L^1$ functional,
with the control variable appearing linearly.


\subsection*{Acknowledgements}

Denysiuk would like to thank AdI -- Innovation Agency,
for the financial support awarded through POFC program,
for the R\&D project SustIMS -- Sustainable Infrastructure Management Systems
(FCOMP-01-0202-FEDER-023113) and to ISISE -- Institute for Sustainability
and Innovation in Structural Engineering (PEst-C/ECI/UI4029/2011 FCOM-01-0124-FEDER-022681).
Silva and Torres were supported by Portuguese funds through the
Center for Research and Development in Mathematics and Applications (CIDMA),
and The Portuguese Foundation for Science and Technology (FCT),
within project PEst-OE/MAT/UI4106/2014. Silva is also grateful to the FCT post-doc
fellowship SFRH/BPD/72061/2010; Torres to the FCT project PTDC/EEI-AUT/1450/2012,
co-financed by FEDER under \mbox{POFC-QREN} with COMPETE reference FCOMP-01-0124-FEDER-028894.
The authors would like to thank the Editor and two anonymous referees
for valuable comments and suggestions.



\end{document}